\documentclass[a4paper,11pt]{amsart}
\usepackage{amssymb,amsfonts,amsxtra,amscd,mathrsfs,stmaryrd}
\usepackage[all]{xy}
\usepackage{fullpage}
\theoremstyle{plain}
\newtheorem{theorem}{Theorem}[section]
\newtheorem{lemma}[theorem]{Lemma}
\newtheorem{cor}[theorem]{Corollary}
\newtheorem{prop}[theorem]{Proposition}
\theoremstyle{definition}
\newtheorem{defi}[theorem]{Definition}
\newtheorem{example}[theorem]{Example}
\theoremstyle{remark}
\newtheorem{rem}[theorem]{Remark}
\numberwithin{equation}{section}
\newcommand{\ai}{\ensuremath{A_\infty}}
\newcommand{\SA}{\ensuremath{\mathcal{SA}_\infty}}
\newcommand{\ISA}{\ensuremath{\mathcal{SA}_\infty}/\cong}
\newcommand{\SAZ}{\ensuremath{\mathcal{SA}_\infty^0}}
\newcommand{\FT}{\ensuremath{\mathcal{OC}}}
\newcommand{\forms}[1]{\ensuremath{\Omega^\bullet(#1)}}
\newcommand{\drham}[2][\bullet]{\ensuremath{DR^{#1}(#2)}}
\newcommand{\cdrham}[2][\bullet]{\ensuremath{\widehat{DR}{\vphantom{DR}}^{#1}(#2)}}
\newcommand{\hlie}{\ensuremath{\mathfrak{h}_{2n|m}}}
\newcommand{\glie}{\ensuremath{\mathfrak{g}_{2n|m}}}
\newcommand{\tglie}{\ensuremath{\tilde{\mathfrak{g}}_{2n|m}}}
\newcommand{\osp}{\ensuremath{\mathfrak{osp}_{2n|m}}}
\newcommand{\gc}[1][\bullet\bullet]{\ensuremath{\mathcal{G}_{#1}}}
\newcommand{\gh}[1][\bullet\bullet]{\ensuremath{H\mathcal{G}_{#1}}}
\newcommand{\cek}[1][\bullet]{\ensuremath{C_{#1}(\glie,\osp)}}
\newcommand{\stlim}[1][\bullet\bullet]{\ensuremath{C_{#1}(\mathfrak{g},\mathfrak{osp})}}
\newcommand{\hstlim}[1][\bullet\bullet]{\ensuremath{H_{#1}(\mathfrak{g},\mathfrak{osp})}}
\newcommand{\chd}[1]{\ensuremath{\mathscr{C}(#1)}}
\newcommand{\ochd}[1]{\ensuremath{\mathscr{OC}(#1)}}
\newcommand{\grpchd}[1]{\ensuremath{\Gamma_{k_1,\ldots,k_m}(#1)}}
\newcommand{\ctalg}[1]{\ensuremath{\widehat{T}(\Pi{#1}^*)}}
\newcommand{\gf}{\ensuremath{\mathbb{C}}}
\newcommand{\dilim}[2]{\ensuremath{\varinjlim_{#1} #2}}
\newcommand{\innprod}{\ensuremath{\langle -,- \rangle}}
\newcommand{\noproof}{\begin{flushright} \ensuremath{\square} \end{flushright}}
\DeclareMathOperator{\sgn}{sgn}
\DeclareMathOperator{\ad}{ad}
\DeclareMathOperator{\ord}{ord}
\DeclareMathOperator{\orb}{orb}
\DeclareMathOperator{\Hom}{Hom}
\DeclareMathOperator{\Der}{Der}
\DeclareMathOperator{\Aut}{Aut}
\DeclareMathOperator{\Vect}{Vect_\gf}

\begin{document}
\begin{abstract}
A standard combinatorial construction, due to Kontsevich, associates to any $\ai$-algebra with an invariant inner product, an inhomogeneous class in the cohomology of the moduli spaces of Riemann surfaces with marked points. We propose an alternative version of this construction based on noncommutative geometry and use it to prove that homotopy equivalent algebras give rise to the same cohomology classes. Along the way we re-prove Kontsevich's theorem relating graph homology to the homology of certain infinite-dimensional Lie algebras. An application to topological conformal field theories is given.
\end{abstract}
\title{Characteristic classes of $\ai$-algebras}
\author{Alastair Hamilton}
\address{Max-Planck-Institut f\"ur Mathematik, Vivatsgasse 7, 53111 Bonn, Germany.}
\email{hamilton@mpim-bonn.mpg.de}
\author{Andrey Lazarev}
\address{Mathematics Department, Bristol University, Bristol, England. BS8 1TW.}
\email{a.lazarev@bristol.ac.uk}
\keywords{$\ai$-algebra, graph homology, topological conformal field theory, noncommutative and symplectic geometry, Feynman calculus.}
\subjclass[2000]{14H10, 17B56, 17B66, 81T18, 81T40.}
\thanks{The work of the first author was supported by the Institut des Hautes \'Etudes Scientifiques and the Max-Planck-Institut f\"ur Mathematik, Bonn. The work of the second author was supported by a grant from the Engineering and Physical Sciences Research Council and by INTAS contract 03513251.}
\maketitle
\tableofcontents
\clearpage

\section{Introduction}

Graph complexes were introduced by Kontsevich in \cite{kontsympgeom} and \cite{kontfeynman}; the main ideas go back to R. Penner \cite{pen}, J. Harer \cite{har} and Culler-Vogtmann \cite{culv}. They were later generalised and further studied by Getzler and Kapranov under the name `Feynman transform', cf. \cite{GK}. These complexes are easy to define and carry a lot of information about the cohomology of mapping class groups, of groups of outer automorphisms of free groups and about invariants of smooth manifolds. On the other hand the complexity of the combinatorics of graphs has so far prevented explicit calculations and to our knowledge there are not even plausible conjectures on the structure of graph (co)homology.

In the absence of any structural or quantitative statements the best one can do is to look for regular ways to produce graph homology classes. The starting point is a homotopy
associative algebra with a type of Poincar\'e duality -- otherwise known as a \emph{symplectic} $\ai$-algebra. Such structures arise in the study of coherent sheaves on Calabi-Yau manifolds and also as cohomology algebras of closed manifolds. The standard construction, due to Kontsevich, associates to any symplectic $\ai$-algebra its \emph{partition function}; an inhomogeneous class in graph homology. An important technical remark is that we work with \emph{even} symplectic structures so that our `Poincare duality' is of the sort that happens to even-dimensional manifolds. The case of the odd symplectic structure admits an analogous treatment, albeit with a twisted version of the graph complex; this will be described in a separate work. Using the interpretation of ribbon graph homology as the cohomology of moduli spaces of Riemann surfaces, one can view this construction as producing  cohomology classes in these moduli spaces.

This construction is nontrivial -- recently Igusa \cite{igusa} and Mondello \cite{mondello} verified the claim made by Kontsevich in \cite{kontfeynman} that a certain simple family of $\ai$-algebras gives rise to the Miller-Morita-Mumford classes in the moduli of curves.

The natural question is then which $\ai$-algebras give rise to homologous cycles. We provide a partial answer to this question. We show that two $\ai$-algebras which are homotopy equivalent have the same characteristic class. For this we give another, combinatorics-free, construction of the partition function which takes values in the relative homology of a certain infinite-dimensional Lie algebra $\mathfrak{g}$ of `noncommutative Hamiltonians'. This is an inhomogeneous homology class which formally resembles the Chern character of a complex vector bundle. Based on this analogy, we call this homology class the \emph{characteristic class} of the $\ai$-algebra. This alternative construction is equivalent to the original one but we believe that it is of independent interest; for instance, it is only in this setting that the homotopy invariance of the construction becomes clear.

A crucial ingredient in our formulation of these constructions is Kontsevich's theorem which interprets graph homology as the relative homology of $\mathfrak{g}$. We give a new, hopefully conceptually clearer, proof of this theorem which makes explicit use of Feynman amplitudes of ribbon graphs.

Another application of our technology is concerned with TCFT's. Consider the topological category whose objects are unions of intervals (`open boundaries') and whose morphisms are given by moduli spaces of Riemann surfaces with a given embedding of intervals into their boundary. This is a monoidal category and an (open) conformal field theory (CFT) is essentially a monoidal functor from this category to the category of vector spaces. Replacing the moduli spaces of surfaces by their cellular chain complexes we get a differential graded (dg) category and the corresponding notion of an open TCFT. Costello \cite{Costello} proved that the data of an open TCFT is equivalent to that of a symplectic $\ai$-algebra (to facilitate the exposition we use a simplified version of his result not involving $D$-branes). The question concerning which open TCFT's correspond to homotopy equivalent algebras was left open. To answer this question we introduce the notion of homologically equivalent TCFT's and outline the proof that homotopy equivalent $\ai$-algebras give rise to TCFT's equivalent in this sense.

We work with ribbon graphs and ribbon graph complexes throughout, although all our results admit straightforward analogues, with the same proofs, for commutative (i.e. undecorated) graphs. It is very likely that our results could be carried through without much work for graphs decorated by a cyclic operad, or even a cyclic dg-operad. It is unclear, however, that this level of abstraction which would necessitate quite a bit of extra verbiage, is justified in the present context. We refer the reader to \cite{CV} and \cite{LV} for the treatment of general decorated graph complexes.

The paper is organised as follows. In section \ref{sec_noncomgeom} we recall the relevant basics from Kontsevich's noncommutative symplectic geometry, references for this material are \cite{hamlaz} and \cite{ginzsympgeom} as well as Kontsevich's original paper \cite{kontsympgeom} and the more recent \cite{KS}. In sections \ref{sec_graph} and \ref{sec_Feynman} we introduce ribbon graphs and graph (co)homology and give a proof of Kontsevich's theorem mentioned above. Section \ref{sec_Ainfinity} collects relevant results on $\ai$-algebras; it is based on \cite{hamlaz}. Section \ref{sec_classes} contains our main results on the characteristic classes of $\ai$-algebras. Finally, section \ref{sec_TCFT} outlines applications to TCFT's.

\subsection{Acknowledgement} The authors would like to express their appreciation to J. Chuang, K. Costello and  M. Movshev for useful discussions on graph complexes.

\subsection{Notation and conventions}

A vector superspace $V$ is a vector space with a natural decomposition $V=V_0\oplus V_1$ into an even part ($V_0$) and a odd part ($V_1$). Given a vector superspace $V$, we will denote the parity reversion of $V$ by $\Pi V$. The parity reversion $\Pi V$ is defined by the formula,
\[ (\Pi V)_0 := V_1, \quad (\Pi V)_1 := V_0. \]
The canonical bijection $V \to \Pi V$ will be denoted by $\Pi$. In this paper we will refer to all vector superspaces simply as `vector spaces'. Likewise, the term `Lie algebra' will be used in place of `Lie superalgebra'.

In this paper we will work over the field of complex numbers $\mathbb{C}$. The canonical $2n|m$-dimensional vector space will be denoted by $\gf^{2n|m}$ and written using the coordinates familiar to symplectic geometry.
\[ \gf^{2n|m}:=\langle p_1,\ldots,p_n;q_1,\ldots,q_n;x_1,\ldots,x_m \rangle. \]
$\gf^{2n|m}$ comes equipped with a canonical nondegenerate skew-symmetric \emph{even} bilinear form $\innprod$ given by the following formula:
\begin{equation} \label{eqn_canonical_form}
\begin{split}
\langle p_i,q_j \rangle = \langle x_i,x_j \rangle & = \delta_{ij}, \\
\langle x_i,p_j \rangle = \langle x_i,q_j \rangle = \langle p_i,p_j \rangle = \langle q_i,q_j \rangle & = 0.
\end{split}
\end{equation}
The Lie algebra of linear endomorphisms of $\gf^{2n|m}$ which annihilate this bilinear form will be denoted by $\osp$.

Let $V$ be a vector space; then we denote the free associative algebra on $V$, with and without unit, by
\[ T(V):=\bigoplus_{i=0}^\infty V^{\otimes i} \quad \text{and} \quad T^+(V):=\bigoplus_{i=1}^\infty V^{\otimes i} \]
respectively. They have a natural grading on them; we say an element $x\in V^{\otimes i}$ has homogeneous \emph{order} $i$. The subspace of $T(V)$ generated by commutators will be denoted by $[T(V),T(V)]$.

Any invertible homomorphism of algebras $\phi:T(V)\to T(W)$ will be called a \emph{diffeomorphism} and any derivation $\xi:T(V)\to T(V)$ will be called a \emph{vector field}. We
denote the Lie algebra of all such vector fields on $T(V)$ by $\Der(T(V))$.

We will denote the symmetric group on $n$ letters by $S_n$ and the cyclic group of order $n$ by $Z_n$. Given a vector space $V$, its tensor power $V^{\otimes n}$ has an action of the cyclic group $Z_n$ on it which is the restriction of the canonical action of $S_n$ to the subgroup $Z_n \cong \langle (n\,n-1\ldots 2\,1) \rangle \subset S_n$. If we define the algebra $\beth$ as
\[ \beth:= \prod_{n=1}^\infty\mathbb{Z}[S_n] \]
then the action of each $S_n$ on $V^{\otimes n}$ for $n\geq 1$ gives $T^+(V)$ the structure of a left $\beth$-module.

Let $z_n$ denote the generator of $Z_n$ corresponding to the cycle $(n\,n-1\ldots 2\,1)$. The norm operator $N_n \in \mathbb{Z}[S_n]$ is the element given by the formula,
\[ N_n:=1+z_n+z_n^2 +\ldots + z_n^{n-1}. \]
We define the anti-symmetriser $\epsilon_n \in \mathbb{Z}[S_n]$ to be the element given by the formula,
\[\epsilon_n:=\sum_{\sigma\in S_n} \sgn(\sigma)\sigma.\]
The operators $z,N,\epsilon \in \beth$ are given by the formulae;
\[ z:=\sum_{n=1}^\infty z_n, \quad N:= \sum_{n=1}^\infty N_n, \quad \epsilon:=\sum_{n=1}^\infty \epsilon_n.\]

Let $k$ be a positive integer. A partition $c$ of $\{1,\ldots 2k\}$ such that every $x \in c$ is a set consisting of precisely two elements will be called a \emph{chord diagram}. The set of all such chord diagrams will be denoted by $\chd{k}$. A chord diagram with an ordering (orientation) of each two point set in that chord diagram will be called an \emph{oriented chord diagram}. The set of all oriented chord diagrams will be denoted by $\ochd{k}$.

$S_{2k}$ acts on the left of $\ochd{k}$ as follows: given an oriented chord diagram $c:=(i_1,j_1),\ldots,(i_k,j_k)$ and a permutation $\sigma \in S_{2k}$,
\begin{equation} \label{eqn_chord_action}
\sigma\cdot c := (\sigma(i_1),\sigma(j_1)),\ldots,(\sigma(i_k),\sigma(j_k)).
\end{equation}
$S_{2k}$ acts on $\chd{k}$ in a similar fashion.

Let $\sigma \in S_m$ and $k_1,\ldots,k_m$ be positive integers such that $k_1+\ldots+k_m=N$, we define the permutation $\sigma_{(k_1,\ldots,k_m)} \in S_N$ by the following commutative diagram:
\begin{equation} \label{eqn_perm_embedding}
\xymatrix{ V^{\otimes k_1}\otimes\ldots\otimes V^{\otimes k_m} \ar@{=}[d] \ar^{\sigma}_{x_1\otimes\ldots\otimes x_m \mapsto x_{\sigma(1)}\otimes\ldots\otimes x_{\sigma(m)}}[rrrr] &&&& V^{\otimes k_{\sigma(1)}}\otimes\ldots\otimes V^{\otimes k_{\sigma(m)}} \ar@{=}[d] \\ V^{\otimes N} \ar^{\sigma_{(k_1,\ldots,k_m)}}_{y_1\otimes\ldots\otimes y_N \mapsto y_{\sigma_{(k_1,\ldots,k_m)}[1]}\otimes\ldots\otimes y_{\sigma_{(k_1,\ldots,k_m)}[N]}}[rrrr] &&&& V^{\otimes N} }
\end{equation}

\section{Noncommutative and symplectic geometry} \label{sec_noncomgeom}

In this section we recall the relevant background material on noncommutative and symplectic geometry as described in \cite{ginzsympgeom}, \cite{hamlaz}, and \cite{kontsympgeom}. This is used to construct a family of Lie (super)algebras consisting of noncommutative Hamiltonians. It is the Lie algebra homology of this family of Lie algebras that will give rise to graph homology and they play a fundamental role in the constructions described in section \ref{sec_classes}.

\subsection{Noncommutative geometry}

Here we recall the basic construction of the noncommutative forms in the simple case of a free associative algebra and prove some basic facts about this definition such as the formal analogue of the Poincar\'e lemma. First, the algebra of noncommutative forms $\Omega^\bullet$ is defined and used to construct the de Rham complex $DR^\bullet$. There then follows an explicit description of the low-dimensional components of this complex and a statement of the formal Poincar\'e lemma. In this section we will provide no proofs for the stated results. All the necessary details can be found in \cite{ginzsympgeom}, \cite{hamlaz} and \cite{kontsympgeom}.

We start by defining the module of 1-forms.

\begin{defi}
Let $V$ be a vector space. The module of noncommutative 1-forms $\Omega^1(V)$ is defined as
\[ \Omega^1(V):= T(V) \otimes T^+(V). \]
$\Omega^1(V)$ has the structure of a $T(V)$-bimodule via the actions
\begin{align*}
& a \cdot (x\otimes y):= ax \otimes y, \\
& (x \otimes y) \cdot a:= x \otimes ya - xy \otimes a;
\end{align*}
for $a,x \in T(V)$ and $y \in T^+(V)$.

Let $d:T(V) \to \Omega^1(V)$ be the map given by the formulae
\begin{displaymath}
\begin{array}{ll}
d(x):= 1 \otimes x, & x \in T^+(V); \\
d(x):= 0, & x \in \gf.
\end{array}
\end{displaymath}
The map $d$ thus defined is a derivation of degree zero.
\end{defi}

Next we introduce the algebra of forms by formally multiplying 1-forms together.

\begin{defi} \label{def_forms}
Let $V$ be a vector space and let $A:=T(V)$. The module of noncommutative forms $\forms{V}$ is defined as
\[ \forms{V}:=T_A\left[\Pi\Omega^1(V)\right]=A \oplus \bigoplus_{i=1}^\infty \underbrace{\Pi\Omega^1(V)\underset{A}{\otimes} \ldots \underset{A}{\otimes} \Pi\Omega^1(V)}_{i \text{ factors}}. \]

Since $\Omega^1(V)$ is an $A$-bimodule, $\forms{V}$ has the structure of an  associative algebra whose multiplication is the standard associative multiplication on the tensor algebra $T_A\left[\Pi\Omega^1(V)\right]$. The map $d:T(V) \to \Omega^1(V)$ lifts uniquely to a map $d:\forms{V} \to \forms{V}$ which gives $\forms{V}$ the structure of a differential graded algebra.
\end{defi}

\begin{rem} \label{rem_form_func}
This construction is functorial; i.e. given vector spaces $V$ and $W$ and a homomorphism of algebras $\phi:T(V) \to T(W)$ (in particular, any linear map from $V$ to $W$ gives rise to such a map), there is a unique map $\phi^*:\forms{V} \to \forms{W}$ extending $\phi$ to a homomorphism of differential graded algebras.
\end{rem}

\begin{rem} \label{rem_form_order}
The module of noncommutative forms has a grading by \emph{order} defined by setting the order of a 1-form $xdy$ to be $\ord(x)+\ord(y)-1$. This suffices to completely determine the order of any homogeneous element in $\forms{V}$; for instance, the $n$-form
\[ x_0dx_1\ldots dx_n \]
has homogeneous order $\ord(x_0)+\ord(x_1)+\ldots+\ord(x_n)-n$.
\end{rem}

It is possible to define noncommutative analogs of the Lie derivative and contraction operator.

\begin{defi} \label{def_operators}
Let $V$ be a vector space and let $\xi:T(V) \to T(V)$ be a vector field:
\begin{enumerate}
\item
We can define a vector field $L_\xi: \forms{V} \to \forms{V}$, called the Lie derivative, by the formulae:
\begin{align*}
L_\xi(x)&:=\xi(x), \\
L_\xi(dx)&:=(-1)^{|\xi|}d(\xi(x)); \\
\end{align*}
for any $x \in T(V)$.
\item
We can define a vector field $i_\xi:\forms{V} \to \forms{V}$, called the contraction operator, by the formulae:
\begin{align*}
i_\xi(x)&:=0, \\
i_\xi(dx)&:=\xi(x); \\
\end{align*}
for any $x \in T(V)$.
\end{enumerate}
\end{defi}

These maps can be shown to satisfy the following identities:

\begin{lemma} \label{lem_operator_identities}
Let $V$ be a vector space and let $\xi,\gamma:T(V) \to T(V)$ be vector fields, then
\begin{enumerate}
\item $L_\xi=[i_\xi,d]$.
\item $[L_\xi,i_\gamma]=i_{[\xi,\gamma]}.$
\item $L_{[\xi,\gamma]}=[L_\xi,L_\gamma].$
\item $[i_\xi,i_\gamma]=0.$
\item $[L_\xi,d]=0.$
\end{enumerate}
\end{lemma}

Now we can define the de Rham complex. Note that the de Rham complex is \emph{not} an algebra.

\begin{defi}
Let $V$ be a vector space. The de Rham complex $\drham{V}$ is defined as
\[ \drham{V}:=\frac{\forms{V}}{\left[\forms{V},\forms{V}\right]}. \]
The differential on $\drham{V}$ is induced by the differential on $\forms{V}$ defined in Definition \ref{def_forms} and is similarly denoted by $d$.
\end{defi}

\begin{rem}
The de Rham complex inherits a grading by order from the grading by order on $\forms{V}$ which was defined in Remark \ref{rem_form_order}.

It follows from Remark \ref{rem_form_func} that the construction of the de Rham complex is also functorial; that is to say that if $\phi:T(V)\to T(W)$ is a homomorphism of algebras then it induces a morphism of complexes $\phi^*:\drham{V} \to \drham{W}$.

Furthermore the Lie and contraction operators $L_\xi,i_\xi:\forms{V}\to\forms{V}$ defined by Definition \ref{def_operators} factor naturally through Lie and contraction operators $L_\xi,i_\xi:\drham{V}\to\drham{V}$, which by an abuse of notation, we denote by the same letters.

We say that an element $x\in\forms{V}$ is a homogeneous $n$-form if it is product of $n$ elements in $\Pi\Omega^1(V)$ (a $0$-form is an element of $T(V)$). Likewise, an element of $\drham{V}$ is an $n$-form if it can be represented by an $n$-form of $\forms{V}$. We define $\drham[n]{V}$ to be the $n$-fold parity reversion of the module of $n$-forms in $\drham{V}$, i.e. we have the identity
\[\drham{V}=\bigoplus_{n=0}^\infty\Pi^n\drham[n]{V} .\]
In particular $\drham[0]{V}$ is the abelianisation of $T(V)$,
\[ \drham[0]{V}:=\frac{T(V)}{\left[T(V),T(V)\right]}=\bigoplus_{i=0}^\infty\left(V^{\otimes i}\right)_{Z_i}. \]

A description of $\drham[1]{V}$ is provided by the following lemma:
\end{rem}

\begin{lemma} \label{lem_oneformiso}
Let $V$ be a vector space, then there is an isomorphism of vector spaces;
\[ \Theta: V\otimes T(V) \to \drham[1]{V}\]
given by the formula, $\Theta(x\otimes y):=dx\cdot y$.
\end{lemma}
\noproof

\begin{rem}
Identifying $V\otimes T(V)$ with $\bigoplus_{i=1}^\infty \left(V^{\otimes i}\right)$ gives a map
\[ \Theta:\bigoplus_{i=1}^\infty V^{\otimes i} \to \drham[1]{V}. \]
\end{rem}

The following lemma describes the map $d:\drham[0]{V} \to \drham[1]{V}$.

\begin{lemma} \label{lem_normdiff}
Let $V$ be a vector space, then the following diagram commutes:
\begin{displaymath}
\xymatrix{ \bigoplus_{i=1}^\infty V^{\otimes i} \ar@{=}[r] & V\otimes T(V) \ar[r]^\Theta & \drham[1]{V} \\
\bigoplus_{i=1}^\infty (V^{\otimes i})_{Z_i} \ar[u]^{N} \ar@{^{(}->}[r] & \bigoplus_{i=0}^\infty (V^{\otimes i})_{Z_i} \ar@{=}[r] & \drham[0]{V} \ar[u]^{d} \\ }
\end{displaymath}
\end{lemma}
\noproof

Given a vector space $V$, we define the de Rham cohomology of $V$ to be the cohomology of the complex $\drham{V}$. We define $H^i(\drham{V})$ as
\[ H^i\left(\drham{V}\right):=\frac{\left\{ x \in \drham[i]{V}:dx=0 \right\}}{d\left(\drham[i-1]{V}\right)}. \]
The cohomology of $V$ is described by the following lemma, which is the noncommutative analogue of the Poincar\'e lemma.

\begin{lemma} \label{lem_poincare}
Let $V$ be a vector space, then the de Rham complex has trivial homology, except in degree zero.
\begin{displaymath}
\begin{array}{lcl}
H^i\left(\drham{V}\right) & = & 0, \quad \text{for all } i \geq 1; \\
H^0\left(\drham{V}\right) & = & \gf. \\
\end{array}
\end{displaymath}
\end{lemma}
\noproof

\subsection{Symplectic geometry}

Here we recall the basic definitions and theorems of noncommutative \emph{symplectic} geometry, cf. \cite{ginzsympgeom}, \cite{hamlaz} and \cite{kontsympgeom}. Again we provide no proofs for any theorems in this section and refer the reader to the cited sources.

\begin{defi} \label{def_symplectic_form}
Let $V$ be a vector space and $\omega \in \drham[2]{V}$ be any 2-form. We say that $\omega$ is a \emph{symplectic form} if:
\begin{enumerate}
\item
it is a closed form, that is to say that $d\omega = 0$;
\item
it is nondegenerate, that is to say that the following map is bijective;
\begin{equation} \label{eqn_nondegenerate}
\begin{array}{ccc}
\Der(T(V)) & \to & \drham[1]{V}, \\
\xi & \mapsto & i_\xi(\omega).
\end{array}
\end{equation}
\end{enumerate}
\end{defi}

\begin{defi} \label{def_sympmaps}
Let $V$ and $W$ be vector spaces and let $\omega \in \drham[2]{V}$ and $\omega' \in \drham[2]{W}$ be symplectic forms:
\begin{enumerate}
\item
We say a vector field $\xi:T(V) \to T(V)$ is a \emph{symplectic vector field} if $L_\xi(\omega)=0$.
\item
We say a diffeomorphism $\phi:T(V) \to T(W)$ is a \emph{symplectomorphism} if $\phi^*(\omega)=\omega'$.
\end{enumerate}
\end{defi}

We have the following simple lemma relating symplectic vector fields to closed 1-forms.

\begin{lemma} \label{lem_sympvfield_form_iso}
Let $V$ be a vector space and $\omega \in \drham[2]{V}$ be a symplectic form. Then the map
\[ \Phi:\Der(T(V)) \to \drham[1]{V}\]
given by the formula $\Phi(\xi):=-i_\xi(\omega)$ induces a one-to-one correspondence between \emph{symplectic} vector fields and \emph{closed} 1-forms:
\[ \Phi: \{ \xi \in \Der(T(V)) : L_\xi(\omega) = 0 \} \overset{\cong}{\longrightarrow} \{ \alpha \in \drham[1]{V} : d\alpha = 0 \}. \]
\end{lemma}
\noproof

The next two results provide an interpretation of closed 2-forms.

\begin{prop} \label{prop_closed_twoforms}
Let $V$ be a vector space, then there is an isomorphism $\Upsilon$ mapping closed 2-forms to $\left[T(V),T(V)\right]$ which makes the following diagram commute:
\begin{displaymath}
\xymatrix{ \{ \omega \in \drham[2]{V} : d\omega = 0 \} \ar[rr]^-{\Upsilon}_-{\cong} && \left[T(V),T(V)\right] \\ \drham[1]{V} \ar[u]^{d} && \bigoplus_{i=1}^\infty V^{\otimes i} \ar[ll]_{\Theta}^{\cong} \ar[u]^{1-z} }
\end{displaymath}
\end{prop}
\noproof

Now let us restrict our attention to 2-forms \emph{of order zero}. Such 2-forms are naturally closed. We have the following lemma describing such 2-forms.

\begin{lemma} \label{lem_symp_innprod}
Let $U$ be a vector space. The map $\Upsilon$ defined by Proposition \ref{prop_closed_twoforms} induces an isomorphism between those 2-forms having order zero and skew-symmetric bilinear forms on $U$:
\[ \Upsilon:\{ \omega \in \drham[2]{U^*} : \ord(\omega) = 0 \} \overset{\cong}{\longrightarrow} (\Lambda^2 U)^*. \]
\end{lemma}
\noproof

Recall that a (skew)-symmetric bilinear form $\innprod$ on $U$ is called \emph{nondegenerate} if the map
\begin{displaymath}
\begin{array}{ccc}
U & \to & U^* \\
x & \mapsto & [a \mapsto \langle x,a\rangle]
\end{array}
\end{displaymath}
is bijective. In this case we call $\innprod$ an \emph{inner product} on $U$. We have the following proposition relating this notion to that of Definition \ref{def_symplectic_form}:

\begin{prop}
Let $U$ be a vector space and let $\omega \in \drham[2]{U^*}$ be a 2-form \emph{of order zero}. Let $\innprod:=\Upsilon(\omega)$ be the skew-symmetric bilinear form corresponding to $\omega$, then the 2-form $\omega$ is nondegenerate if and only if the bilinear form $\innprod$ is nondegenerate.
\end{prop}
\noproof

\subsection{Lie algebras of vector fields}

The purpose of this section is to introduce a series of Lie algebras $\glie$ and show that they define Lie algebras of symplectic vector fields. We verify the Jacobi identity for the bracket on $\glie$ and provide an explicit formula for the bracket in Lemma \ref{lem_bracket_formula}.

Recall from Lemma \ref{lem_sympvfield_form_iso} that given a vector space $V$ and a symplectic form $\omega \in \drham[2]{V}$, there is an isomorphism between symplectic vector fields and closed 1-forms
\[ \Phi: \{ \xi \in \Der(T(V)) : L_\xi(\omega) = 0 \} \overset{\cong}{\longrightarrow} \{ \alpha \in \drham[1]{V} : d\alpha = 0 \}, \]
given by the formula $\Phi(\xi):=-i_\xi(\omega)$. This means that to any function in $\drham[0]{V}$ we can associate a certain symplectic vector field;
\begin{displaymath}
\begin{array}{ccc}
\drham[0]{V} & \to & \{ \xi \in \Der(T(V)) : L_\xi(\omega) = 0 \}, \\
a & \mapsto & \Phi^{-1}(da).
\end{array}
\end{displaymath}
We will typically denote the resulting vector field by the corresponding Greek letter, for instance in this case we have $\alpha:=\Phi^{-1}(da)$.

Let $V:=\gf^{2n|m}$ be the canonical $2n|m$-dimensional vector space written in the standard coordinates
\[ \gf^{2n|m}:=\langle p_1,\ldots,p_n;q_1,\ldots,q_n;x_1,\ldots,x_m \rangle. \]
$V$ can be endowed with the canonical \emph{even} symplectic form
\begin{equation} \label{eqn_canonicalsymplecticform}
\omega:= \sum_{i=1}^n dp_i dq_i + \frac{1}{2}\sum_{i=1}^m dx_i dx_i.
\end{equation}
Any other symplectic vector space of the same dimension is isomorphic to this one. This symplectic vector space can be used to define a series of Lie algebras as follows:

\begin{defi}
Given integers $n,m\geq 0$ we can define a Lie algebra $\hlie$ with the underlying space
\[ \hlie:=\drham[0]{V}=\bigoplus_{i=0}^\infty\left(V^{\otimes i}\right)_{Z_i},\]
where $V:=\gf^{2n|m}$ is the canonical $2n|m$-dimensional symplectic vector space.

The Lie bracket on $\hlie$ is given by the formula;
\begin{equation} \label{eqn_bracket}
\{ a,b \}:= (-1)^a L_\alpha(b), \quad \text{for any } a,b \in \hlie.
\end{equation}

We define a new Lie algebra $\glie$ as the Lie subalgebra of $\hlie$ which has the underlying vector space
\[\glie:=\bigoplus_{i=2}^\infty\left(V^{\otimes i}\right)_{Z_i}.\]
The Lie algebra $\tglie$ is defined to be the Lie subalgebra of $\glie$ which has the underlying vector space
\[\tglie:=\bigoplus_{i=3}^\infty\left(V^{\otimes i}\right)_{Z_i}.\]

We denote the direct limits of these Lie algebras by
\[ \mathfrak{h}:=\dilim{n,m}{\hlie}, \quad \mathfrak{g}:=\dilim{n,m}{\glie} \quad \text{and} \quad \tilde{\mathfrak{g}}:=\dilim{n,m}{\tglie}. \]
\end{defi}

The fact that the bracket defined by equation \eqref{eqn_bracket} gives $\hlie$ the structure of a Lie algebra still needs to be verified. This verification is provided by the following proposition.

\begin{prop}
The bracket $\{ -,- \}$ defined on $\hlie$ by equation \eqref{eqn_bracket} is skew-symmetric and satisfies the Jacobi identity.
\end{prop}

\begin{proof}
We must use the identities of Lemma \ref{lem_operator_identities} to verify that the bracket $\{-,-\}$ satisfies the aforementioned identities. First we prove that the bracket is skew-symmetric.
\begin{displaymath}
\begin{split}
\{a,b\} & = (-1)^a L_\alpha(b) = (-1)^a i_\alpha (db), \\
        & = (-1)^{a+1} i_\alpha i_\beta(\omega) = (-1)^{a+1}(-1)^{(a-1)(b-1)} i_\beta i_\alpha(\omega), \\
        & = -(-1)^{ab}(-1)^b L_\beta(a) = -(-1)^{ab}\{b,a\}. \\
\end{split}
\end{displaymath}
Next we verify the Jacobi identity. This will follow from the following equation.
\begin{equation} \label{eqn_liemap}
\begin{split}
d\{a,b\} =&(-1)^a dL_\alpha(b) = L_\alpha(db) = -L_\alpha i_\beta(\omega), \\
=& -[L_\alpha,i_\beta](\omega) = -i_{[\alpha,\beta]}(\omega) = \Phi([\alpha,\beta]). \\
\end{split}
\end{equation}
From this equation, the Jacobi identity follows thusly.
\begin{displaymath}
\begin{split}
\{\{a,b\},c\} &= (-1)^{a+b}L_{[\alpha,\beta]}(c), \\
&= (-1)^{a+b}L_\alpha L_\beta(c) -(-1)^{ab+a+b}L_\beta L_\alpha(c), \\
&= \{a,\{b,c\}\} -(-1)^{ab}\{b,\{a,c\}\}.
\end{split}
\end{displaymath}
\end{proof}

\begin{prop} \label{prop_hamsympiso}
The map from the Lie algebra $\hlie$ which has underlying vector space $\drham[0]{V}$ to the Lie algebra of symplectic vector fields
\begin{displaymath}
\begin{array}{ccc}
\hlie & \to & \{ \xi \in \Der(T(V)) : L_\xi(\omega) = 0 \} \\
a & \mapsto & \alpha:=\Phi^{-1}(da)
\end{array}
\end{displaymath}
is a surjective homomorphism of Lie algebras with kernel $\gf$.
\end{prop}

\begin{proof}
The statement that the map is surjective and has kernel $\gf$ follows from Lemma \ref{lem_poincare} and equation \eqref{eqn_liemap} implies that it is a map of Lie algebras.
\end{proof}

\begin{rem}
Consider the Lie subalgebra $\mathfrak{k}_{2n|m}$ of $\glie$ which has the underlying vector space
\[\mathfrak{k}_{2n|m}:=\left(V^{\otimes 2}\right)_{Z_2}.\]
The Lie algebra $\mathfrak{k}_{2n|m}$ acts on $\hlie$ via the adjoint action, in particular it acts on $V$, $\glie$ and $\tglie$. It follows from Proposition \ref{prop_hamsympiso} that the following map is an isomorphism of Lie algebras
\begin{displaymath}
\begin{array}{ccc}
\mathfrak{k}_{2n|m} & \to & \osp, \\
a & \mapsto & \left[\{a,-\}: x \mapsto \{a,x\},\quad \text{for } x \in V\right]; \\
\end{array}
\end{displaymath}
where $\osp$ is the Lie algebra of \emph{linear} symplectic vector fields. The consequence of this is that $\glie$ splits as a semi-direct product of $\osp$ and $\tglie$;
\[ \glie = \osp\ltimes \tglie, \]
where $\osp$ acts canonically on $\tglie$ according to the Leibniz rule.
\end{rem}

\begin{lemma} \label{lem_bracket_formula}
Let $a:=a_1\ldots a_k,b:=b_1\ldots b_l \in \hlie$ where $a_i,b_j\in V:=\gf^{2n|m}$, then we have the following explicit formula for the bracket $\{-,-\}$ on $\hlie$:
\begin{equation} \label{eqn_bracket_formula}
\{a,b\}=\sum_{i=1}^k\sum_{j=1}^l (-1)^p \langle a_i,b_j\rangle(z_{k-1}^{i-1}\cdot a_1\ldots\hat{a_i}\ldots a_k)(z_{l-1}^{j-1}\cdot b_1\ldots\hat{b_j}\ldots b_l),
\end{equation}
where $p:=|a_i|(|a_{i+1}|+\ldots+|a_k|+|b_1|+\ldots+|b_{j-1}|)$ and $\innprod$ is the canonical nondegenerate skew-symmetric bilinear form on $V$ (cf. equation \eqref{eqn_canonical_form}).
\end{lemma}

\begin{proof}
First, we calculate the vector field $\alpha$ corresponding to $a$. We have,
\[ da=\sum_{i=1}^k (-1)^{a_i(a_1+\ldots+a_{i-1})}da_i \cdot (z_{k-1}^{i-1}\cdot a_1\ldots\hat{a_i}\ldots a_k) .\]
Let $A_i:=(-1)^{a_i(a_1+\ldots+a_{i-1})}z_{k-1}^{i-1}\cdot a_1\ldots\hat{a_i}\ldots a_k$. Using the identity
\[v=\sum_{j=1}^n\left[\langle v,q_j\rangle p_j - \langle v,p_j\rangle q_j\right] +\sum_{j=1}^m\langle v,x_j\rangle x_j, \quad \text{for all } v\in V;\]
we have that
\[ da = \sum_{j=1}^n\sum_{i=1}^k\left[\langle a_i,q_j \rangle dp_j\cdot A_i - \langle a_i,p_j \rangle dq_j\cdot A_i\right] + \sum_{j=1}^m\sum_{i=1}^k \langle a_i,x_j \rangle dx_j\cdot A_i.\]
Since
\[i_\alpha(\omega)=\sum_{j=1}^n\left[(-1)^a dq_j\cdot\alpha(p_j)-(-1)^a dp_j\cdot\alpha(q_j)\right]+\sum_{j=1}^m dx_j\cdot\alpha(x_j)\]
we conclude that for all $v\in V$
\[\alpha(v)=(-1)^{av+a+v}\sum_{i=1}^k (-1)^{a_i(a_1+\ldots+a_{i-1})}\langle a_i,v \rangle (z_{k-1}^{i-1}\cdot a_1\ldots\hat{a_i}\ldots a_k).\]
An immediate consequence of this formula is that
\[ \{u,v\}=\langle u,v \rangle; \quad \text{for all } u,v\in V.\]
Equation \eqref{eqn_bracket_formula} is  established by the following calculation:
\begin{displaymath}
\begin{split}
\{a,b\} &= (-1)^a L_\alpha(b) = (-1)^a\sum_{j=1}^l (-1)^{b_j(b_1+\ldots+b_{j-1})}\alpha(b_j)(z_{l-1}^{j-1}\cdot b_1\ldots\hat{b_j}\ldots b_l) \\
&= (-1)^a\sum_{i=1}^k\sum_{j=1}^l(-1)^q\langle a_i,b_j\rangle(z_{k-1}^{i-1}\cdot a_1\ldots\hat{a_i}\ldots a_k)(z_{l-1}^{j-1}\cdot b_1\ldots\hat{b_j}\ldots b_l) \\
&= \sum_{i=1}^k\sum_{j=1}^l (-1)^p \langle a_i,b_j\rangle (z_{k-1}^{i-1}\cdot a_1\ldots\hat{a_i}\ldots a_k)(z_{l-1}^{j-1}\cdot b_1\ldots\hat{b_j}\ldots b_l);
\end{split}
\end{displaymath}
where
\begin{align*}
q&:=|a||b_j|+|a|+|b_j|+|b_j|(|b_1|+\ldots+|b_{j-1}|)+|a_i|(|a_1|+\ldots+|a_{i-1}|), \\
p&:=|a_i|(|a_{i+1}|+\ldots+|a_k|+|b_1|+\ldots+|b_{j-1}|).
\end{align*}
The calculation with the signs on the last line follows by assuming that $|a_i|=|b_j|$, which is possible since the bilinear form $\innprod$ is even.
\end{proof}

\subsection{Lie algebra homology}

In this section we briefly recall the definition of Lie algebra homology and \emph{relative} Lie algebra homology in the $\mathbb{Z}/2\mathbb{Z}$-graded setting.

\begin{defi}
Let $\mathfrak{l}$ be a Lie algebra: the underlying space of the Chevalley-Eilenberg complex of $\mathfrak{l}$ is the exterior algebra $\Lambda_{\bullet}(\mathfrak{l})$ which is defined to be the quotient of $T(\mathfrak{l})$ by the ideal generated by the relation
\[ g \otimes h = -(-1)^{|g||h|} h \otimes g; \quad g,h \in \mathfrak{l}. \]
There is a natural grading on $\Lambda_{\bullet}(\mathfrak{l})$ where an element $g\in\mathfrak{l}$ has bidegree $(1,|g|)$ in $\Lambda_{\bullet\bullet}(\mathfrak{l})$ and total degree $|g|+1$; this implicitly defines a bigrading and total grading on the whole of $\Lambda_{\bullet}(\mathfrak{l})$.

The differential $d:\Lambda_{i}(\mathfrak{l}) \to \Lambda_{i-1}(\mathfrak{l})$ is defined by the following formula:
\[ d(g_1\wedge\ldots\wedge g_m):= \sum_{1\leq i < j \leq m} (-1)^{p(g)} [g_i,g_j]\wedge g_1\wedge\ldots\wedge \hat{g_i}\wedge\ldots\wedge \hat{g_j}\wedge\ldots\wedge g_m, \]
for $g_1,\ldots,g_m \in \mathfrak{l}$; where
\[ p(g):=|g_i|(|g_1|+\ldots+|g_{i-1}|)+|g_j|(|g_1|+\ldots+|g_{j-1}|)+|g_i||g_j|+i+j-1. \]
The differential $d$ has bidegree $(-1,0)$ in the above bigrading. We will denote the Chevalley-Eilenberg complex of $\mathfrak{l}$ by $C_{\bullet}(\mathfrak{l})$; the homology of $C_{\bullet}(\mathfrak{l})$ will be called the Lie algebra homology of $\mathfrak{l}$ and will be denoted by $H_{\bullet}(\mathfrak{l})$.
\end{defi}

\begin{rem}
Let $\mathfrak{l}$ be a Lie algebra; $\mathfrak{l}$ acts on $\Lambda(\mathfrak{l})$ via the adjoint action. This action commutes with the Chevalley-Eilenberg differential $d$; in fact it is nullhomotopic (cf. \cite{fuchscohom} or \cite{loday}).

As a consequence of this the space
\[ C_{\bullet}(\tglie)_{\osp} \]
of $\osp$-coinvariants of the Chevalley-Eilenberg complex of $\tglie$ forms a complex when equipped with the Chevalley-Eilenberg differential $d$. This is the \emph{relative} Chevalley-Eilenberg complex of $\glie$ modulo $\osp$ (cf. \cite{fuchscohom}) and is denoted by $\cek$. The homology of this complex is called the \emph{relative} homology of $\glie$ modulo $\osp$ and is denoted by $H_{\bullet}(\glie,\osp)$.

Let $\mathfrak{osp}$ denote the direct limit of the Lie algebras $\osp$, then the stable limit of the relative Chevalley-Eilenberg complexes is given by
\[ \stlim[\bullet]=\dilim{n,m}{\cek}.\]
We denote the homology of this stable limit by $\hstlim[\bullet]$.
\end{rem}

\section{Graph homology} \label{sec_graph}

In this section we will describe the combinatorial model for \emph{ribbon graphs} or \emph{fat graphs} that we use throughout the rest of the paper. We describe various operations on graphs such as contracting edges and expanding vertices and how these operations can be used to define a complex which is generated by isomorphism classes of oriented ribbon graphs.

\subsection{Basic definitions and operations on graphs}

In this section we set out the basic definition of a graph and its variants such as \emph{fully ordered graphs}, \emph{ribbon graphs} and \emph{oriented graphs}. We describe basic operations on these various types of graphs, such as contracting edges and expanding vertices, in terms of the combinatorial model for such graphs given by the subsequent definitions.

We will use the following set theoretic definition for the basic notion of a graph (see Igusa \cite{igusa}):

\begin{defi}
Choose a fixed infinite set $\Xi$ which is disjoint from its power set; a graph $\Gamma$ is a finite subset of $\Xi$ (the set of half-edges) together with:
\begin{enumerate}
\item
a partition $V(\Gamma)$ of $\Gamma$ ($V(\Gamma)$ is the set of vertices of $\Gamma$),
\item
a partition $E(\Gamma)$ of $\Gamma$ into sets having cardinality equal to two ($E(\Gamma)$ is the set of edges of $\Gamma$).
\end{enumerate}
We say that a vertex $v \in V$ has valency $n$ if $v$ has cardinality $n$. The elements of $v$ are called the \emph{incident half-edges} of $v$. An edge whose composite half-edges are incident to the same vertex will be called a loop. The set consisting of those edges of $\Gamma$ which are loops will be denoted by $L(\Gamma)$. In this paper we shall assume that all our graphs have vertices of valency $\geq 3$ (except in section \ref{sec_TCFT}, where we allow graphs with legs).
\end{defi}

The following definition of a fully ordered graph will be useful in section \ref{sec_Feynman} when we come to defining certain Feynman amplitudes.

\begin{defi} \label{def_graph}
A \emph{fully ordered graph} $\Gamma$ is a graph together with
\begin{enumerate}
\item \label{item_graphdummy0}
a linear ordering of the vertices,
\item
a linear ordering of the half-edges in every edge of $\Gamma$ (i.e. an orientation of each edge, as for a directed graph),
\item \label{item_graphdummy1}
a linear ordering of the incident half-edges for every vertex of $\Gamma$.
\end{enumerate}
We say that a fully ordered graph has type $(k_1,\ldots,k_m)$ if it has $m$ vertices and the $i$th vertex has valency $k_i$.
\end{defi}

Fully ordered graphs are rigid objects in that they have no nontrivial automorphisms. As such, their use in this paper is as a notational device, rather than as an interesting object of study.

\begin{example}
A fully ordered graph of type $(3,3,3,3)$:
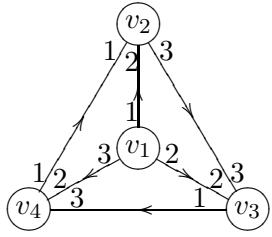
\begin{displaymath}
\begin{xy}
(0,0)="a" *[o]=<16pt>{v_1} *\frm{o} ;
"a"+a(90);"a" ; **{} ; "a";?/16mm/="b" *[o]=<16pt>{v_2} *\frm{o} ;
"a"+a(330);"a" ; **{} ; "a";?/16mm/="c" *[o]=<16pt>{v_3} *\frm{o} ;
"a"+a(210);"a" ; **{} ; "a";?/16mm/="d" *[o]=<16pt>{v_4} *\frm{o} ;
"b";"a" ; **{} ; ?(0)/8pt/="t1" ; "a";"b" ; **{} ; ?(0)/8pt/="t2" ; "t1";"t2" **\dir{-} ; ?*\dir{>} ; ?(0)/6pt/*_!/2.5pt/{1} ; ?(1)/-6pt/*_!/2.5pt/{2} ;
"c";"a" ; **{} ; ?(0)/8pt/="t1" ; "a";"c" ; **{} ; ?(0)/8pt/="t2" ; "t1";"t2" **\dir{-} ; ?*\dir{>} ; ?(0)/4pt/*_!/4.5pt/{2} ; ?(1)/-8pt/*_!/4.5pt/{2} ;
"d";"a" ; **{} ; ?(0)/8pt/="t1" ; "a";"d" ; **{} ; ?(0)/8pt/="t2" ; "t1";"t2" **\dir{-} ; ?*\dir{>} ; ?(0)/4pt/*_!/-4.5pt/{3} ; ?(1)/-8pt/*_!/-4.5pt/{2} ;
"c";"b" ; **{} ; ?(0)/8pt/="t1" ; "b";"c" ; **{} ; ?(0)/8pt/="t2" ; "t1";"t2" **\dir{-} ; ?*\dir{>} ; ?(0)/6pt/*_!/4pt/{3} ; ?(1)/-6pt/*_!/4pt/{3} ;
"d";"c" ; **{} ; ?(0)/8pt/="t1" ; "c";"d" ; **{} ; ?(0)/8pt/="t2" ; "t1";"t2" **\dir{-} ; ?*\dir{>} ; ?(0)/10pt/*_!/-4.5pt/{1} ; ?(1)/-10pt/*_!/-4.5pt/{3} ;
"b";"d" ; **{} ; ?(0)/8pt/="t1" ; "d";"b" ; **{} ; ?(0)/8pt/="t2" ; "t1";"t2" **\dir{-} ; ?*\dir{>} ; ?(0)/6pt/*_!/4pt/{1} ; ?(1)/-6pt/*_!/4pt/{1} ;
\end{xy}
\end{displaymath}
\end{example}

Let $\Gamma$ be a graph with $n$ edges and $m$ vertices which have valency $k_1,\ldots,k_m$ respectively. Let $\chi$ be the set of structures on $\Gamma$ consisting of items (\ref{item_graphdummy0}-\ref{item_graphdummy1}) in Definition \ref{def_graph}, then the group
\begin{equation} \label{eqn_perm_group}
S_m \times (\underbrace{S_2 \times \ldots \times S_2}_{n \text{ copies}})\times(S_{k_1}\times\ldots\times S_{k_m})
\end{equation}
acts freely and transitively on the right of $\chi$. The first factor changes the ordering of vertices, the central factors flip the orientation of the edges and the final factors change the ordering of the incident half-edges at each vertex.

Given an oriented chord diagram $c:=(i_1,j_1),\ldots,(i_k,j_k)$ and a sequence of positive integers $k_1,\ldots,k_m$ such that $k_1+\ldots+k_m=2k$ we can construct a fully ordered graph $\Gamma=\grpchd{c}$ of type $(k_1,\ldots,k_m)$ with half-edges $h_1,\ldots,h_{2k}$ as follows:
\begin{enumerate}
\item \label{item_grpchddummy0}
The vertices of $\Gamma$ are
\[ (h_1,\ldots,h_{k_1}),(h_{k_1+1},\ldots,h_{k_1+k_2}),\ldots,(h_{k_1+\ldots+k_{m-1}+1},\ldots,h_{k_1+\ldots+k_m}). \]
\item \label{item_grpchddummy1}
The edges of $\Gamma$ are
\[ (h_{i_1},h_{j_1}),\ldots,(h_{i_k},h_{j_k}). \]
\item
The set of all vertices is ordered as in \eqref{item_grpchddummy0} and the linear ordering of the incident half-edges for each vertex is also given by \eqref{item_grpchddummy0}. The edges are oriented as in \eqref{item_grpchddummy1}.
\end{enumerate}


Conversely, given a fully ordered graph $\Gamma$ of type $(k_1,\ldots,k_m)$ we can define an oriented chord diagram $c_\Gamma\in\ochd{k}$. This gives us the following lemma,
which is tautological in nature, but allows us to describe the combinatorics of graphs in terms of the actions of permutation groups on chord diagrams.

\begin{lemma} \label{lem_graph_chord}
The set of all (isomorphism classes of) fully ordered graphs of type $(k_1,\ldots,k_m)$ is in one-to-one correspondence with the set of all oriented chord diagrams.
\begin{displaymath}
\begin{array}{ccc}
\ochd{k} & \cong & \mathrm{Fully \ ordered \ graphs} \\
c & \rightharpoonup & \grpchd{c} \\
c_\Gamma & \leftharpoondown & \Gamma
\end{array}
\end{displaymath}

Furthermore, these maps are equivariant with respect to the actions of the permutation groups defined by \eqref{eqn_perm_group} and \eqref{eqn_chord_action}; that is to say that given positive integers $k_1,\ldots,k_m$ such that $k_1+\ldots+k_m=2k$ and given an oriented chord diagram $c:=(i_1,j_1),\ldots,(i_k,j_k)$ we have that:
\begin{enumerate}
\item
\[\grpchd{(i_r,j_r)\cdot c} = \grpchd{c}\cdot \tau_r,\]
where $\tau_r\in S_2$ reverses the orientation on the corresponding edge.
\item
For all $\mathbf{\sigma}\in S_{k_1}\times\ldots\times S_{k_m}\subset S_{2k}$;
\begin{displaymath}
\begin{array}{ccc}
\grpchd{\mathbf{\sigma}^{-1}\cdot c} & \cong & \grpchd{c}\cdot\mathbf{\sigma}, \\
h_i & \mapsto & h_{\sigma(i)}.
\end{array}
\end{displaymath}
\item
For all $\sigma\in S_m$;
\begin{displaymath}
\begin{array}{ccc}
\Gamma_{k_{\sigma(1)},\ldots,k_{\sigma(m)}}({\sigma_{(k_1,\ldots,k_m)}}^{-1}\cdot c) & \cong & \grpchd{c}\cdot\sigma, \\
h_i & \mapsto & h_{\sigma_{(k_1,\ldots,k_m)}[i]}.
\end{array}
\end{displaymath}
\end{enumerate}
\end{lemma}
\noproof

Once again, let $\Gamma$ be a graph with $n$ edges and $m$ vertices of valency $k_1,\ldots,k_m$ respectively. Now, consider the subgroup
\[G\lhd S_m \times (\underbrace{S_2 \times \ldots \times S_2}_{n \text{ copies}})\]
consisting of permutations having total sign one, i.e. $G$ is the kernel of the map
\begin{displaymath}
\begin{array}{ccc}
S_m \times S_2 \times \ldots \times S_2 & \to & \{ 1,-1 \} \\
(\sigma,\tau_1,\ldots,\tau_n) & \mapsto & \sgn\sigma \sgn\tau_1 \ldots \sgn\tau_n \\
\end{array}
\end{displaymath}
Define $K$ to be the subgroup of $S_{k_1}\times\ldots\times S_{k_m}$ consisting of cyclical permutations:
\[K:=Z_{k_1}\times\ldots\times Z_{k_m}.\]

Recall that $\chi$ is the set of structures on $\Gamma$ consisting of items (\ref{item_graphdummy0}-\ref{item_graphdummy1}) in Definition \ref{def_graph}. This terminology allows us to make the following definition of an \emph{oriented ribbon graph}:

\begin{defi} \label{def_ribbon_graph}
An \emph{oriented ribbon graph} $\Gamma$ is a graph together with an element in $\chi/(G\times K)$; that is to say that it is a graph together with a pair of orbits, the first corresponding to the group $G$ and the last corresponding to the group $K$. The orbit $\sigma$ corresponding to the group $G$ is called the \emph{orientation} on $\Gamma$, in particular there are only two possible choices for $\sigma$. The orbit $\pi$ corresponding to $K$ consists of a cyclic ordering of the incident half-edges at every vertex of $\Gamma$ and defines the \emph{ribbon structure} on $\Gamma$.

A \emph{ribbon graph} $\Gamma$ is a graph together with just the orbit $\pi$ corresponding to the group $K$ and without the orbit corresponding to the orientation.
\end{defi}

It follows from the definition that every \emph{oriented ribbon graph} is represented by some \emph{fully ordered graph}. When we want to emphasise the orientation on an oriented ribbon graph, we will denote the oriented ribbon graph by $(\Gamma,\sigma)$, where $\Gamma$ is a ribbon graph and $\sigma$ is the orientation on $\Gamma$; usually it will be denoted simply by $\Gamma$.

Having introduced the basic notion of a graph and its variants, we now describe their corresponding morphisms. We restrict our attention to isomorphisms alone.

\begin{defi}
Let $\Gamma$ and $\Gamma'$ be two graphs, we say $\Gamma$ is isomorphic to $\Gamma'$ if there exists a bijective map $f:\Gamma \to \Gamma'$ between the half-edges of $\Gamma$ and $\Gamma'$ such that:
\begin{enumerate}
\item
the image of a vertex of $\Gamma$ under $f$ is a vertex of $\Gamma'$,
\item
the image of an edge of $\Gamma$ under $f$ is an edge of $\Gamma'$.
\end{enumerate}

If $\Gamma$ and $\Gamma'$ are \emph{oriented ribbon graphs} then we say that they are isomorphic if there is a map $f$ as above which preserves the orientations and the ribbon structures on $\Gamma$ and $\Gamma'$. More explicitly, given orientations $\omega$ and $\omega'$ on $\Gamma$ and $\Gamma'$ respectively and a map $f$ as above, there is a naturally induced orientation $f(\omega)$ on $\Gamma'$ defined in an obvious way; we say that $f$ preserves the orientations of $\Gamma$ and $\Gamma'$ if $f(\omega)=\omega'$. Likewise, given ribbon structures $\pi$ and $\pi'$ on $\Gamma$ and $\Gamma'$ respectively and a map $f$ as above, there is a naturally induced ribbon structure $f(\pi)$ on $\Gamma'$ defined by transferring the cyclic orderings on the incident half-edges at each vertex of $\Gamma$ through the map $f$. We say that $f$ preserves the ribbon structures on $\Gamma$ and $\Gamma'$ if $f(\pi)=\pi'$.

Given a \emph{ribbon graph} $\Gamma$ we say that a map $f:\Gamma\to\Gamma$ is a \emph{ribbon graph automorphism} if it is an isomorphism of graphs which preserves the ribbon structure on $\Gamma$. If $\Gamma$ is an \emph{oriented ribbon graph} then we say that $f$ is an \emph{oriented ribbon graph automorphism} if it also preserves the orientation on $\Gamma$. The group of all oriented ribbon graph automorphisms will be denoted by $\Aut(\Gamma)$.
\end{defi}

In the remainder of this section we describe various operations on the previously defined types of graphs. These operations will be used to construct the differential in the graph complex in the next section. The first operation is that of contracting an edge.

\begin{defi} \label{def_edge_contract}
Let $\Gamma$ be an oriented ribbon graph and $e=(a,b)$ be one of its edges, which we assume is not a loop. The half-edges $a$ and $b$ will be incident to vertices denoted by
\[v=(h_1,\ldots,h_k,a) \text{ and } v'=(h'_1,\ldots h'_l,b) \]
respectively. We may assume that the orientation on $\Gamma$ is represented by an ordering of the vertices of the form
\[ (v,v',v_3,\ldots,v_m) \]
and where the edge $e$ is oriented in the order $e=(a,b)$. We define a new oriented ribbon graph $\Gamma/e$ as follows:
\begin{enumerate}
\item
The set of half-edges comprising $\Gamma/e$ are the half-edges of $\Gamma$ minus $a$ and $b$.
\item \label{item_edgecontractdummy0}
The set of vertices of $\Gamma/e$ is
\[ v_0,v_3,\ldots,v_m; \]
where $v_0:=(h_1,\ldots,h_k,h'_1,\ldots,h'_l)$.
\item
The set of edges of $\Gamma/e$ is the set of edges of $\Gamma$ minus the edge $e$.
\item
The orientation on $\Gamma/e$ is given by ordering the vertices as above in \eqref{item_edgecontractdummy0}; the edges of $\Gamma/e$ are oriented in the same way as the edges of $\Gamma$. It is clear that this definition produces a well defined orientation modulo the actions of the relevant permutation groups.
\item
The ribbon structure on $\Gamma/e$ is defined by setting the cyclic orderings of the incident half-edges for $v_3,\ldots,v_m$ to be the same as those for $\Gamma$ and ordering the incident half-edges at the first vertex $v_0$ as above in \eqref{item_edgecontractdummy0}.
\end{enumerate}
\end{defi}

\begin{example}
Contracting an edge in a ribbon graph.
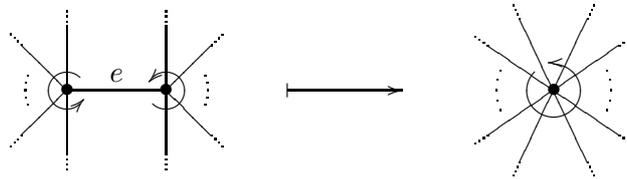
\begin{displaymath}
\begin{xy}
(10,0);(-5,0) ; **\dir{-} ; ?(1)*\dir{>} ; ?(0)*\dir{|} ;
(-4,0)-(30,0)="a" ; "a"+(13,0)="b" ;
"a"*\cir<2.5mm>{dl_dr} ; "a"+a(315);"a" ; **{} ; ?(0)/2.5mm/;"a" ; **{} ; ?(1)/-0.25pt/*_!/2pt/\dir{<} ;
"b"*\cir<2.5mm>{ur_ul} ; "b"+a(135);"b" ; **{} ; ?(0)/2.5mm/;"b" ; **{} ; ?(1)/-0.25pt/*_!/2pt/\dir{<} ;
"a";"b" ; **\dir{-} ; ?+(0,2)*{e} ;
"a"+a(90);"a" ; **{} ; ?/8mm/="t";"a" ; **\dir{-} ; "a"+a(90);"a" ; **{} ; ?/10mm/;"t" ; **\dir{.} ;
"a"+a(135);"a" ; **{} ; ?/8mm/="t";"a" ; **\dir{-} ; "a"+a(135);"a" ; **{} ; ?/10mm/;"t" ; **\dir{.} ;
"a"+a(225);"a" ; **{} ; ?/8mm/="t";"a" ; **\dir{-} ; "a"+a(225);"a" ; **{} ; ?/10mm/;"t" ; **\dir{.} ;
"a"+a(270);"a" ; **{} ; ?/8mm/="t";"a" ; **\dir{-} ; "a"+a(270);"a" ; **{} ; ?/10mm/;"t" ; **\dir{.} ;
"a"+a(180);"a" ; **{} ; ?/5mm/*\dir{.} ; "a"+a(190);"a" ; **{} ; ?/5mm/*\dir{.} ; "a"+a(200);"a" ; **{} ; ?/5mm/*\dir{.} ; "a"+a(170);"a" ; **{} ; ?/5mm/*\dir{.} ; "a"+a(160);"a" ; **{} ; ?/5mm/*\dir{.} ;
"a"*{\bullet} ;
"b"+a(90);"b" ; **{} ; ?/8mm/="t";"b" ; **\dir{-} ; "b"+a(90);"b" ; **{} ; ?/10mm/;"t" ; **\dir{.} ;
"b"+a(45);"b" ; **{} ; ?/8mm/="t";"b" ; **\dir{-} ; "b"+a(45);"b" ; **{} ; ?/10mm/;"t" ; **\dir{.} ;
"b"+a(270);"b" ; **{} ; ?/8mm/="t";"b" ; **\dir{-} ; "b"+a(270);"b" ; **{} ; ?/10mm/;"t" ; **\dir{.} ;
"b"+a(315);"b" ; **{} ; ?/8mm/="t";"b" ; **\dir{-} ; "b"+a(315);"b" ; **{} ; ?/10mm/;"t" ; **\dir{.} ;
"b"-a(180);"b" ; **{} ; ?/5mm/*\dir{.} ; "b"-a(190);"b" ; **{} ; ?/5mm/*\dir{.} ; "b"-a(200);"b" ; **{} ; ?/5mm/*\dir{.} ; "b"-a(170);"b" ; **{} ; ?/5mm/*\dir{.} ; "b"-a(160);"b" ; **{} ; ?/5mm/*\dir{.} ;
"b"*{\bullet} ;
(0,0)+(30,0)="a" ;
"a"+a(35);"a" ; **{} ; ?/10mm/="t";"a" ; **\dir{-} ; "a"+a(35);"a" ; **{} ; ?/12mm/;"t" ; **\dir{.} ;
"a"+a(65);"a" ; **{} ; ?/10mm/="t";"a" ; **\dir{-} ; "a"+a(65);"a" ; **{} ; ?/12mm/;"t" ; **\dir{.} ;
"a"+a(115);"a" ; **{} ; ?/10mm/="t";"a" ; **\dir{-} ; "a"+a(115);"a" ; **{} ; ?/12mm/;"t" ; **\dir{.} ;
"a"+a(145);"a" ; **{} ; ?/10mm/="t";"a" ; **\dir{-} ; "a"+a(145);"a" ; **{} ; ?/12mm/;"t" ; **\dir{.} ;
"a"+a(215);"a" ; **{} ; ?/10mm/="t";"a" ; **\dir{-} ; "a"+a(215);"a" ; **{} ; ?/12mm/;"t" ; **\dir{.} ;
"a"+a(245);"a" ; **{} ; ?/10mm/="t";"a" ; **\dir{-} ; "a"+a(245);"a" ; **{} ; ?/12mm/;"t" ; **\dir{.} ;
"a"+a(295);"a" ; **{} ; ?/10mm/="t";"a" ; **\dir{-} ; "a"+a(295);"a" ; **{} ; ?/12mm/;"t" ; **\dir{.} ;
"a"+a(325);"a" ; **{} ; ?/10mm/="t";"a" ; **\dir{-} ; "a"+a(325);"a" ; **{} ; ?/12mm/;"t" ; **\dir{.} ;
"a"+a(180);"a" ; **{} ; ?/7mm/*\dir{.} ; "a"+a(190);"a" ; **{} ; ?/7mm/*\dir{.} ; "a"+a(200);"a" ; **{} ; ?/7mm/*\dir{.} ; "a"+a(170);"a" ; **{} ; ?/7mm/*\dir{.} ; "a"+a(160);"a" ; **{} ; ?/7mm/*\dir{.} ;
"a"-a(180);"a" ; **{} ; ?/7mm/*\dir{.} ; "a"-a(190);"a" ; **{} ; ?/7mm/*\dir{.} ; "a"-a(200);"a" ; **{} ; ?/7mm/*\dir{.} ; "a"-a(170);"a" ; **{} ; ?/7mm/*\dir{.} ; "a"-a(160);"a" ; **{} ; ?/7mm/*\dir{.} ;
"a"*\cir<3.7mm>{dl^l} ; "a"+a(90);"a" ; **{} ; ?(0)/3.7mm/;"a" ; **{} ; ?(1)/-0.25pt/*_!/2pt/\dir{<} ;
"a"*{\bullet} ;
\end{xy}
\end{displaymath}
\end{example}

The next operation we describe is how to expand vertices. In order to expand vertices it is necessary to introduce the notion of an \emph{ideal edge}.

\begin{defi}
Let $\Gamma$ be an oriented ribbon graph and $v\in V(\Gamma)$ be a vertex of $\Gamma$. An \emph{ideal edge} $i$ of the vertex $v$ consists of the following data:
\begin{enumerate}
\item
a partition of the incident half-edges of $v$ into two disjoint sets $a$ and $b$ which both have cardinality$\geq 2$,
\item
linear orderings on $a$ and $b$ such that the juxtaposition of these linear orderings gives the cyclic ordering on $v$.
\end{enumerate}

We say that $i$ is an \emph{ideal edge} of $\Gamma$ if it is an ideal edge of some vertex of $\Gamma$ and denote the set of all ideal edges of $\Gamma$ by $I(\Gamma)$.
\end{defi}

\begin{rem}
Note that if $f:\Gamma\to\Gamma'$ is an isomorphism of oriented ribbon graphs then to any ideal edge $i\in I(\Gamma)$ there naturally corresponds an ideal edge $f(i)\in I(\Gamma')$.
\end{rem}

\begin{defi} \label{def_edge_expand}
Let $\Gamma$ be an oriented ribbon graph with $m$ vertices. We may assume that the orientation on $\Gamma$ is represented by an ordering of the vertices of the form
\[ (v_1,v_2,\ldots,v_m). \]
Let
\[ i=\{(h_1,\ldots,h_k),(h'_1,\ldots,h'_l)\} \]
be an ideal edge of the vertex $v_1$ of $\Gamma$. We define a new oriented ribbon graph $\Gamma\backslash i$ as follows:
\begin{enumerate}
\item
The set of half-edges comprising $\Gamma\backslash i$ are the half-edges of $\Gamma$ plus two new half-edges $a$ and $b$.
\item
The set of edges of $\Gamma\backslash i$ is the set of edges of $\Gamma$ plus the edge $e:=(a,b)$.
\item \label{item_edgeexpanddummy0}
The set of vertices of $\Gamma\backslash i$ is
\[ v,v',v_2,\ldots,v_m; \]
where
\[ v:=(h_1,\ldots,h_k,a) \quad \text{and} \quad v':=(h'_1,\ldots,h'_l,b).\]
\item
The orientation on $\Gamma\backslash i$ is given by ordering the vertices as above in \eqref{item_edgeexpanddummy0}; the edges of $\Gamma\backslash i$ are oriented in the same way as the edges of $\Gamma$ with the new edge being oriented as $e:=(a,b)$. It is clear that this definition produces a well defined orientation modulo the actions of the relevant permutation groups.
\item
The ribbon structure on $\Gamma\backslash i$ is defined by setting the cyclic orderings of the incident half-edges for $v_2,\ldots,v_m$ to be the same as those for $\Gamma$ and ordering the incident half-edges of $v$ and $v'$ as above in \eqref{item_edgeexpanddummy0}.
\end{enumerate}
\end{defi}

\begin{rem} \label{rem_contract_expand}
Ideal edge expansion is the inverse operation to edge contraction. Let $\Gamma$ be an oriented ribbon graph, then any edge $e$ of $\Gamma$ canonically determines an ideal edge $i_e$ of $\Gamma/e$; moreover there is a canonical isomorphism between $(\Gamma/e)\backslash i_e$ and $\Gamma$. Conversely, any ideal edge $i$ of $\Gamma$ canonically determines an edge $e_i$ of $\Gamma\backslash i$ such that $(\Gamma\backslash i)/e_i=\Gamma$. This is illustrated by the following diagram:
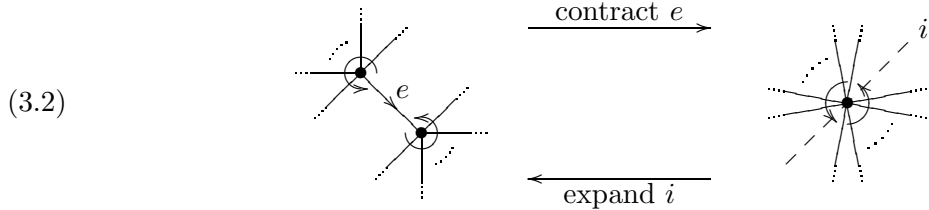
\begin{equation} \label{eqn_contract_expand}
\begin{xy}
(-12,10);(12,10) ; **\dir{-} ; ?(0)*\dir{<} ; ?+(0,2.2)*{\text{contract} \ e} ;
(12,-10);(-12,-10) ; **\dir{-} ; ?(0)*\dir{<} ; ?-(0,2.2)*{\text{expand} \ i} ;
(-4,4)-(30,0)="a" ; "a"+(8,-8)="b" ;
"a"*\cir<2.2mm>{l_d} ; "a"+(0.9,-2.1)*\dir{>} ;
"b"*\cir<2.2mm>{r_u} ; "b"+(-0.9,2.1)*\dir{<} ;
"a";"b" ; **\dir{-} ; ?*!/3pt/\dir{<} +(1.5,1.5)*{e} ;
"a"+a(45);"a" ; **{} ; ?/6mm/="t";"a" ; **\dir{-} ; "a"+a(45);"a" ; **{} ; ?/8mm/;"t" ; **\dir{.} ;
"a"+a(90);"a" ; **{} ; ?/6mm/="t";"a" ; **\dir{-} ; "a"+a(90);"a" ; **{} ; ?/8mm/;"t" ; **\dir{.} ;
"a"+a(180);"a" ; **{} ; ?/6mm/="t";"a" ; **\dir{-} ; "a"+a(180);"a" ; **{} ; ?/8mm/;"t" ; **\dir{.} ;
"a"+a(225);"a" ; **{} ; ?/6mm/="t";"a" ; **\dir{-} ; "a"+a(225);"a" ; **{} ; ?/8mm/;"t" ; **\dir{.} ;
"a"+a(135);"a" ; **{} ; ?/4mm/*\dir{.} ; "a"+a(145);"a" ; **{} ; ?/4mm/*\dir{.} ; "a"+a(155);"a" ; **{} ; ?/4mm/*\dir{.} ; "a"+a(125);"a" ; **{} ; ?/4mm/*\dir{.} ; "a"+a(115);"a" ; **{} ; ?/4mm/*\dir{.} ;
"a"*{\bullet} ;
"b"+a(45);"b" ; **{} ; ?/6mm/="t";"b" ; **\dir{-} ; "b"+a(45);"b" ; **{} ; ?/8mm/;"t" ; **\dir{.} ;
"b"+a(0);"b" ; **{} ; ?/6mm/="t";"b" ; **\dir{-} ; "b"+a(0);"b" ; **{} ; ?/8mm/;"t" ; **\dir{.} ;
"b"+a(225);"b" ; **{} ; ?/6mm/="t";"b" ; **\dir{-} ; "b"+a(225);"b" ; **{} ; ?/8mm/;"t" ; **\dir{.} ;
"b"+a(270);"b" ; **{} ; ?/6mm/="t";"b" ; **\dir{-} ; "b"+a(270);"b" ; **{} ; ?/8mm/;"t" ; **\dir{.} ;
"b"-a(135);"b" ; **{} ; ?/4mm/*\dir{.} ; "b"-a(145);"b" ; **{} ; ?/4mm/*\dir{.} ; "b"-a(155);"b" ; **{} ; ?/4mm/*\dir{.} ; "b"-a(125);"b" ; **{} ; ?/4mm/*\dir{.} ; "b"-a(115);"b" ; **{} ; ?/4mm/*\dir{.} ;
"b"*{\bullet} ;
(0,0)+(30,0)="a" ;
"a"+a(80);"a" ; **{} ; ?/8mm/="t";"a" ; **\dir{-} ; "a"+a(80);"a" ; **{} ; ?/10mm/;"t" ; **\dir{.} ;
"a"+a(100);"a" ; **{} ; ?/8mm/="t";"a" ; **\dir{-} ; "a"+a(100);"a" ; **{} ; ?/10mm/;"t" ; **\dir{.} ;
"a"+a(170);"a" ; **{} ; ?/8mm/="t";"a" ; **\dir{-} ; "a"+a(170);"a" ; **{} ; ?/10mm/;"t" ; **\dir{.} ;
"a"+a(190);"a" ; **{} ; ?/8mm/="t";"a" ; **\dir{-} ; "a"+a(190);"a" ; **{} ; ?/10mm/;"t" ; **\dir{.} ;
"a"+a(260);"a" ; **{} ; ?/8mm/="t";"a" ; **\dir{-} ; "a"+a(260);"a" ; **{} ; ?/10mm/;"t" ; **\dir{.} ;
"a"+a(280);"a" ; **{} ; ?/8mm/="t";"a" ; **\dir{-} ; "a"+a(280);"a" ; **{} ; ?/10mm/;"t" ; **\dir{.} ;
"a"+a(350);"a" ; **{} ; ?/8mm/="t";"a" ; **\dir{-} ; "a"+a(350);"a" ; **{} ; ?/10mm/;"t" ; **\dir{.} ;
"a"+a(10);"a" ; **{} ; ?/8mm/="t";"a" ; **\dir{-} ; "a"+a(10);"a" ; **{} ; ?/10mm/;"t" ; **\dir{.} ;
"a"+a(135);"a" ; **{} ; ?/6mm/*\dir{.} ; "a"+a(145);"a" ; **{} ; ?/6mm/*\dir{.} ; "a"+a(155);"a" ; **{} ; ?/6mm/*\dir{.} ; "a"+a(125);"a" ; **{} ; ?/6mm/*\dir{.} ; "a"+a(115);"a" ; **{} ; ?/6mm/*\dir{.} ;
"a"-a(135);"a" ; **{} ; ?/6mm/*\dir{.} ; "a"-a(145);"a" ; **{} ; ?/6mm/*\dir{.} ; "a"-a(155);"a" ; **{} ; ?/6mm/*\dir{.} ; "a"-a(125);"a" ; **{} ; ?/6mm/*\dir{.} ; "a"-a(115);"a" ; **{} ; ?/6mm/*\dir{.} ;
"a"*\cir<3mm>{l^dr} ; "a"*\cir<3mm>{r^ul} ; "a"+a(225);"a" ; **{} ; ?/2.3mm/="t" ; "t"-a(135);"t" ; **{} ; ?/0.3mm/*\dir{>} ; "a"+a(45);"a" ; **{} ; ?/2.3mm/="t" ; "t"-a(315);"t" ; **{} ; ?/0.3mm/*\dir{>} ;
"a"-(8,8);"a"+(8,8) ; **\dir{--} ; ?(0)+(2,2)*{i} ;
"a"*{\bullet} ;
\end{xy}
\end{equation}

This correspondence can be stated more precisely. Let $\mathcal{C}_1$ be the category whose objects are ordered pairs $(\Gamma,e)$, where $\Gamma$ is an oriented ribbon graph and $e$ is an edge of $\Gamma$ which is not a loop. The morphisms
\[ f:(\Gamma,e) \to (\Gamma',e') \]
in the category $\mathcal{C}_1$ are isomorphisms $f:\Gamma\to\Gamma'$ such that $f(e)=e'$.

Let $\mathcal{C}_2$ be the category whose objects are ordered pairs $(\Gamma,i)$, where $\Gamma$ is an oriented ribbon graph and $i$ is an ideal edge of $\Gamma$. The morphisms
\[ f:(\Gamma,i) \to (\Gamma',i') \]
in the category $\mathcal{C}_2$ are isomorphisms $f:\Gamma\to\Gamma'$ such that $f(i)=i'$.

Diagram \eqref{eqn_contract_expand} describes an equivalence of the categories $\mathcal{C}_1$ and $\mathcal{C}_2$:
\begin{displaymath}
\begin{array}{ccc}
\mathcal{C}_1 & \approx & \mathcal{C}_2 \\
(\Gamma,e) & \rightharpoonup & (\Gamma/e,i_e) \\
(\Gamma\backslash i,e_i) & \leftharpoondown & (\Gamma,i) \\
\end{array}
\end{displaymath}
\end{rem}

\subsection{Graph homology and cohomology}

In this section we define the graph complex and graph homology. The differential in this graph complex is given by contracting edges. We also describe the complex that computes graph cohomology by introducing a differential that expands vertices.

Recall that in Definition \ref{def_ribbon_graph} it was noted that for any ribbon graph $\Gamma$ there were two possible choices of orientations on it, therefore given an orientation $\sigma$ on $\Gamma$ we may define $-\sigma$ to be the other choice of orientation on $\Gamma$. We are now ready to define the graph complex:

\begin{defi}
The underlying space of the graph complex $\gc$ is the free $\gf$-vector space generated by isomorphism classes of \emph{oriented ribbon graphs}, modulo the following relation:
\[ (\Gamma,-\sigma)=-(\Gamma,\sigma) \]
for any ribbon graph $\Gamma$ with orientation $\sigma$. The differential $\partial$ on $\gc$ is given by the following formula:
\begin{equation} \label{eqn_hdiff}
\partial(\Gamma):=\sum_{e \in E(\Gamma)-L(\Gamma)} \Gamma/e
\end{equation}
where $\Gamma$ is any oriented ribbon graph and the sum is taken over all the edges of $\Gamma$ which are not loops.

There is a natural bigrading on $\gc$ where $\gc[ij]$ is the vector space generated by isomorphism classes of oriented ribbon graphs with $i$ vertices and $j$ edges. In this bigrading the differential has bidegree $(-1,-1)$. The homology of this complex is called \emph{graph homology} and denoted by $\gh$.
\end{defi}

\begin{rem}
Since contracting edges in the opposite order produces opposite orientations, it follows that $\partial^2=0$; i.e. $\partial$ is a differential, cf. \cite{hamgraph}. Note that if a ribbon graph $\Gamma$ has an orientation reversing automorphism, then the oriented ribbon graph $(\Gamma,\sigma)$ is zero in $\gc$ for any choice of orientation $\sigma$ on $\Gamma$. For this reason we could usually assume that any ribbon graph automorphism is orientation preserving.
\end{rem}

\begin{rem}
The graph complex splits up as a direct sum of subcomplexes which are indexed by the Euler characteristic of the corresponding graph,
\[\gc=\bigoplus_{\chi=0}^{-\infty}\left[\bigoplus_{i-j=\chi}\gc[ij]\right].\]
Moreover, since all the vertices of any graph are at least trivalent, each of these subcomplexes has finite dimension.
\end{rem}

It is also possible to define a cohomological differential which expands vertices. At the end of this section we show that the associated complex computes graph cohomology.

\begin{defi}
We can define a cohomological differential $\delta$ of bidegree $(1,1)$ on the vector space $\gc$. First, we re-weight the graphs by introducing the notation $\langle\Gamma\rangle:=|\Aut(\Gamma)|\Gamma$. The differential $\delta$ is then defined by the formula
\[ \delta\langle\Gamma\rangle:=\sum_{i\in I(\Gamma)}\langle\Gamma\backslash i\rangle.\]
The fact that $\delta^2=0$ will follow from Proposition \ref{prop_adjointness}.
\end{defi}

\begin{defi} \label{def_graphinn}
There is a natural bilinear form $\innprod$ on $\gc$ coming from the natural basis of $\gc$ generated by ribbon graphs. Let $(\Gamma,\sigma)$ and $(\Gamma',\sigma')$ be oriented ribbon graphs and assume that they have no orientation reversing automorphisms, then $\innprod$ is given by the formula:
\[ \langle(\Gamma,\sigma),(\Gamma',\sigma')\rangle:=\left\{\begin{array}{cl} 1, & \text{if } (\Gamma,\sigma)\cong(\Gamma',\sigma'); \\ -1, & \text{if } (\Gamma,\sigma)\cong(\Gamma',-\sigma'); \\ 0, & \text{otherwise.} \end{array}\right. \]
In particular, the above expression is nonzero if and only if the ribbon graphs $\Gamma$ and $\Gamma'$ are isomorphic.
\end{defi}

\begin{rem}
Although the graph complex is infinite dimensional, each of its bigraded components $\gc[ij]$ is finite dimensional. Moreover the bilinear form $\innprod$ is given by a series of bilinear forms
\[ \innprod_{ij}:\gc[ij]\otimes\gc[ij]\to\gf \]
which are clearly nondegenerate.
\end{rem}

The following proposition shows that the differential $\delta$ computes graph cohomology.

\begin{prop} \label{prop_adjointness}
The differential $\partial$ on $\gc$ is adjoint to the differential $\delta$ under the bilinear form $\innprod$. More explicitly, for any two oriented ribbon graphs $\Gamma$ and $\Gamma'$ the following equation holds:
\begin{equation} \label{eqn_adjointness}
\langle\partial\Gamma,\Gamma'\rangle=\langle\Gamma,\delta\Gamma'\rangle.
\end{equation}
\end{prop}

\begin{proof}
From the definitions of $\partial$ and $\delta$ we see that \eqref{eqn_adjointness} is equivalent to
\[\sum_{e\in E(\Gamma)-L(\Gamma)}|\Aut(\Gamma')|\langle\Gamma/e,\Gamma'\rangle = \sum_{i\in I(\Gamma)}|\Aut(\Gamma'\backslash i)|\langle\Gamma,\Gamma'\backslash i\rangle.\]

Let $\Gamma_1$ and $\Gamma_2$ be any two ribbon graphs and assume that they have no orientation reversing automorphisms. Let $\mathcal{I}(\Gamma_1,\Gamma_2)$ denote the set of all ribbon graph isomorphisms (not necessarily orientation preserving) between $\Gamma_1$ and $\Gamma_2$, then $\Aut(\Gamma_1)$ and $\Aut(\Gamma_2)$ act freely and transitively on $\mathcal{I}(\Gamma_1,\Gamma_2)$. Since $\langle(\Gamma_1,\sigma_1),(\Gamma_2,\sigma_2)\rangle$ is nonzero if and only if $\Gamma_1$ and $\Gamma_2$ are isomorphic, it follows that \eqref{eqn_adjointness} is equivalent to
\[\sum_{e\in E(\Gamma)-L(\Gamma)} \sum_{f\in\mathcal{I}(\Gamma/e,\Gamma')}\langle\Gamma/e,\Gamma'\rangle = \sum_{i\in I(\Gamma')}\sum_{f\in\mathcal{I}(\Gamma,\Gamma'\backslash i)}\langle\Gamma,\Gamma'\backslash i\rangle.\]

The result now follows by constructing a bijection between the terms in the sum on the left hand side and the terms in the sum on the right hand side. This bijection is precisely that which is provided by the equivalence of categories described in Remark \ref{rem_contract_expand}.

More explicitly, let $e$ be an edge of $\Gamma$ and $f:\Gamma/e\to\Gamma'$ be an isomorphism; then $i':=f(i_e)$ is an ideal edge of $\Gamma'$ and $f$ extends to an isomorphism from $\Gamma=(\Gamma/e)\backslash i_e$ to $\Gamma'\backslash i'$. Since \eqref{eqn_contract_expand} describes an equivalence of categories, it follows that this map provides a bijection between the terms in the sum.
\end{proof}

\begin{cor}
Let $H^\delta_{ij}$ be the $ij$th component of the cohomology of the complex
\[ \left(\gc; \ \delta:\gc\to\gc[\bullet+1,\bullet+1]\right) \]
and let $H_{\partial^*}^{ij}$ be the $ij$th component of the cohomology of the complex
\[ \left(\Hom(\gc,\gf); \ \partial^*:\Hom(\gc[\bullet\bullet],\gf)\to\Hom(\gc[\bullet+1,\bullet+1],\gf)\right), \]
then $H^\delta_{ij}\cong H_{\partial^*}^{ij}$.
\end{cor}
\noproof

\section{Feynman calculus} \label{sec_Feynman}

This section relates the combinatorial structure present in the graph complex to the algebraic structure present in the relative Chevalley-Eilenberg complex of the Lie algebra $\mathfrak{g}$. This is achieved by making explicit use of Feynman diagrams and section \ref{sec_kontthm} contains a proof of Kontsevich's theorem exploiting this approach. The section concludes with a discussion of the Hopf algebra structures present in both these complexes.

\subsection{Feynman amplitudes} \label{sec_Feynman_amplitudes}

In this section we define certain Feynman amplitudes which we describe in terms of invariants of the Lie algebra $\osp$. This will allow us to apply the fundamental theorems of invariant theory for this Lie algebra in order to establish an isomorphism between graph and Lie algebra homology.

Given $\sigma\in S_k$, let $\tilde{\sigma}\in S_{2k}$ be the unique permutation that makes the following diagram commute for any vector space $V$,
\[ \xymatrix{V^{\otimes 2k} \ar[r]^{\tilde{\sigma}}& V^{\otimes 2k} \\ (V\otimes V)^{\otimes k} \ar@{=}[u] \ar[r]^{\sigma} & (V\otimes V)^{\otimes k} \ar@{=}[u]}\]
This gives an embedding
\begin{displaymath}
\begin{array}{ccc}
S_k & \hookrightarrow & S_{2k}, \\
\sigma & \mapsto & \tilde{\sigma}.
\end{array}
\end{displaymath}
We will denote the set of right cosets of this embedding by $S_{k}\backslash S_{2k}$. We now make the following definition.

\begin{defi}
Any oriented chord diagram $c \in \ochd{k}$ canonically determines a right coset $\sigma_c\in S_{k}\backslash S_{2k}$. If we write the oriented chord diagram $c$ as
\[c:= (i_1,j_1),\ldots,(i_k,j_k);\]
then the permutation $\sigma_c$ is given by
\begin{displaymath}
\begin{array}{cccccccc}
& i_1 & j_1 & i_2 & j_2 & & i_k & j_k \\
\sigma_c:= & \downarrow & \downarrow & \downarrow & \downarrow & \cdots & \downarrow & \downarrow \\
& 1 & 2 & 3 & 4 & & 2k-1 & 2k \\
\end{array}.
\end{displaymath}
\end{defi}

This definition satisfies the following identity which is easily verified; for any $\tau \in S_{2k}$ and any oriented chord diagram $c\in\ochd{k}$,
\begin{equation} \label{eqn_perm_identity}
\sigma_{\tau\cdot c}=\sigma_c\cdot\tau^{-1}.
\end{equation}

Let $V:=\gf^{2n|m}$ and let $\kappa:V^{\otimes 2k} \to \gf$ be the map defined by the formula
\[ \kappa(x_1\otimes\ldots\otimes x_{2k}):= \langle x_1,x_2 \rangle\langle x_3,x_4 \rangle\ldots\langle x_{2k-1},x_{2k}\rangle, \]
where $\innprod$ is the canonical form given by \eqref{eqn_canonical_form}. It follows from the definition that $\kappa$ is $\osp$-invariant. Using this map and the previous definition we can describe the fundamental invariants of the Lie algebra $\osp$.

\begin{defi}
Given an oriented chord diagram $c\in \ochd{k}$ we define the map $\beta_c:V^{\otimes 2k}\to\gf$ by the formula:
\[\beta_c(x):=\kappa(\sigma_c\cdot x), \quad \text{for } x\in V^{\otimes 2k}.\]
\end{defi}

\begin{rem}
Let $c:=(i_1,j_1),\ldots,(i_k,j_k)$ be an oriented chord diagram and $x_1,\ldots,x_{2k}\in V$, then a slightly more explicit formula for $\beta_c$ is (modulo 2)
\[ \beta_c(x_1\otimes\ldots\otimes x_{2k}):=\pm \langle x_{i_1},x_{j_1} \rangle\langle x_{i_2},x_{j_2} \rangle\ldots\langle x_{i_k},x_{j_k}\rangle.\]

Note that since the canonical form $\innprod$ is \emph{even} it follows that $\beta_c$ is well defined, i.e. it is independent of the choice of representative for $\sigma_c\in S_{k}\backslash S_{2k}$. It follows from \eqref{eqn_perm_identity} that for any oriented chord diagram $c \in \ochd{k}$ and any $\tau\in S_{2k}$,
\begin{equation} \label{eqn_inv_identity}
\beta_{\tau\cdot c}(x)=\beta_c(\tau^{-1}\cdot x).
\end{equation}
\end{rem}

Now we have all the terminology in place that we need in order to define Feynman amplitudes for fully ordered graphs.

\begin{defi}
Let $V:=\gf^{2n|m}$ and let $\Gamma$ be a fully ordered graph of type $(k_1,\ldots,k_i)$, we define the \emph{Feynman amplitude} map $\tilde{F}_\Gamma:T(T^+(V))\to\gf$ by the formula:
\begin{equation} \label{eqn_feynamp}
\tilde{F}_\Gamma(x):=\left\{\begin{array}{ll}\beta_{c_\Gamma}(x), & x\in V^{\otimes k_1}\otimes\ldots\otimes V^{\otimes k_i};\\ 0, & \text{for all other summands.}\end{array}\right.
\end{equation}
\end{defi}

\begin{rem}
This definition is best explained graphically. Let $x:=x_1\otimes\ldots\otimes x_i$ be a tensor in $V^{\otimes k_1}\otimes\ldots\otimes V^{\otimes k_i}$, where $x_j:=x_{j1}\otimes\ldots\otimes x_{jk_j}$ is a tensor in $V^{\otimes k_j}$. Given any fully ordered graph $\Gamma$ we can produce a number $\tilde{F}_\Gamma(x)$ by placing the tensor $x_j$ at the $j$th vertex of $\Gamma$; the edges of the graph provide a way to contract these tensors using the canonical form $\innprod$ and produce a number. This is illustrated by the following diagram.
\end{rem}

\begin{example} \label{exm_feynman_amplitude}
The Feynman amplitude of a graph of type $(4,3,3,3,3)$.
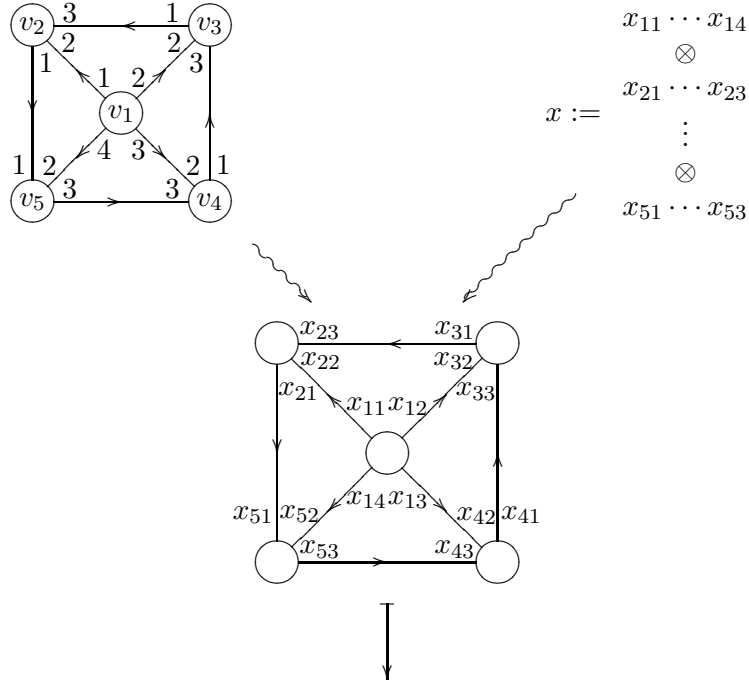
\begin{displaymath}
\begin{xy}
(-35,0)="a" ; "a"*[o]=<16pt>{v_1} *\frm{o} ;
"a"+a(135);"a" ; **{} ; "a";?/16mm/="b" *[o]=<16pt>{v_2} *\frm{o} ;
"a"+a(45);"a" ; **{} ; "a";?/16mm/="c" *[o]=<16pt>{v_3} *\frm{o} ;
"a"+a(315);"a" ; **{} ; "a";?/16mm/="d" *[o]=<16pt>{v_4} *\frm{o} ;
"a"+a(225);"a" ; **{} ; "a";?/16mm/="e" *[o]=<16pt>{v_5} *\frm{o} ;
"b";"a" ; **{} ; ?(0)/8pt/="t1" ; "a";"b" ; **{} ; ?(0)/8pt/="t2" ; "t1";"t2" **\dir{-} ; ?*\dir{>} ; ?(0)/6pt/*_!/-5pt/{1} ; ?(1)/-6pt/*_!/-5pt/{2} ;
"c";"a" ; **{} ; ?(0)/8pt/="t1" ; "a";"c" ; **{} ; ?(0)/8pt/="t2" ; "t1";"t2" **\dir{-} ; ?*\dir{>} ; ?(0)/6pt/*_!/5pt/{2} ; ?(1)/-6pt/*_!/5pt/{2} ;
"d";"a" ; **{} ; ?(0)/8pt/="t1" ; "a";"d" ; **{} ; ?(0)/8pt/="t2" ; "t1";"t2" **\dir{-} ; ?*\dir{>} ; ?(0)/6pt/*_!/-5pt/{3} ; ?(1)/-6pt/*_!/5pt/{2} ;
"e";"a" ; **{} ; ?(0)/8pt/="t1" ; "a";"e" ; **{} ; ?(0)/8pt/="t2" ; "t1";"t2" **\dir{-} ; ?*\dir{>} ; ?(0)/6pt/*_!/5pt/{4} ; ?(1)/-6pt/*_!/-5pt/{2} ;
"d";"e" ; **{} ; ?(0)/8pt/="t1" ; "e";"d" ; **{} ; ?(0)/8pt/="t2" ; "t1";"t2" **\dir{-} ; ?*\dir{>} ; ?(0)/6pt/*_!/5pt/{3} ; ?(1)/-6pt/*_!/5pt/{3} ;
"c";"d" ; **{} ; ?(0)/8pt/="t1" ; "d";"c" ; **{} ; ?(0)/8pt/="t2" ; "t1";"t2" **\dir{-} ; ?*\dir{>} ; ?(0)/6pt/*_!/-5pt/{1} ; ?(1)/-6pt/*_!/5pt/{3} ;
"b";"c" ; **{} ; ?(0)/8pt/="t1" ; "c";"b" ; **{} ; ?(0)/8pt/="t2" ; "t1";"t2" **\dir{-} ; ?*\dir{>} ; ?(0)/6pt/*_!/-5pt/{1} ; ?(1)/-6pt/*_!/-5pt/{3} ;
"e";"b" ; **{} ; ?(0)/8pt/="t1" ; "b";"e" ; **{} ; ?(0)/8pt/="t2" ; "t1";"t2" **\dir{-} ; ?*\dir{>} ; ?(0)/6pt/*_!/5pt/{1} ; ?(1)/-6pt/*_!/-5pt/{1} ;
(35,0)="x" ; "x"*{x:=\begin{array}{c} x_{11}\cdots x_{14} \\ \otimes \\ x_{21}\cdots x_{23} \\ \vdots \\ \otimes \\ x_{51}\cdots x_{53}\end{array}} ;
"d"+(5,-5);(-10,-25) ; **\dir{~} ; ?(0)*\dir{<} ;
"x"+(-10,-10);(10,-25) ; **\dir{~} ; ?(0)*\dir{<} ;
(0,-45)="a" ; "a"*[o]=<16pt>{} *\frm{o} ;
"a"+a(135);"a" ; **{} ; "a";?/20mm/="b" *[o]=<16pt>{} *\frm{o} ;
"a"+a(45);"a" ; **{} ; "a";?/20mm/="c" *[o]=<16pt>{} *\frm{o} ;
"a"+a(315);"a" ; **{} ; "a";?/20mm/="d" *[o]=<16pt>{} *\frm{o} ;
"a"+a(225);"a" ; **{} ; "a";?/20mm/="e" *[o]=<16pt>{} *\frm{o} ;
"b";"a" ; **{} ; ?(0)/8pt/="t1" ; "a";"b" ; **{} ; ?(0)/8pt/="t2" ; "t1";"t2" **\dir{-} ; ?*!/-2pt/\dir{>} ; ?(0)/10pt/*_!/-7pt/{x_{11}} ; ?(1)/-8pt/*_!/-7pt/{x_{22}} ;
"c";"a" ; **{} ; ?(0)/8pt/="t1" ; "a";"c" ; **{} ; ?(0)/8pt/="t2" ; "t1";"t2" **\dir{-} ; ?*!/-2pt/\dir{>} ; ?(0)/10pt/*_!/7pt/{x_{12}} ; ?(1)/-8pt/*_!/7pt/{x_{32}} ;
"d";"a" ; **{} ; ?(0)/8pt/="t1" ; "a";"d" ; **{} ; ?(0)/8pt/="t2" ; "t1";"t2" **\dir{-} ; ?*!/-2pt/\dir{>} ; ?(0)/10pt/*_!/-7pt/{x_{13}} ; ?(1)/-10pt/*_!/7pt/{x_{42}} ;
"e";"a" ; **{} ; ?(0)/8pt/="t1" ; "a";"e" ; **{} ; ?(0)/8pt/="t2" ; "t1";"t2" **\dir{-} ; ?*!/-2pt/\dir{>} ; ?(0)/10pt/*_!/7pt/{x_{14}} ; ?(1)/-11pt/*_!/-7pt/{x_{52}} ;
"d";"e" ; **{} ; ?(0)/8pt/="t1" ; "e";"d" ; **{} ; ?(0)/8pt/="t2" ; "t1";"t2" **\dir{-} ; ?*\dir{>} ; ?(0)/8pt/*_!/5pt/{x_{53}} ; ?(1)/-8pt/*_!/5pt/{x_{43}} ;
"c";"d" ; **{} ; ?(0)/8pt/="t1" ; "d";"c" ; **{} ; ?(0)/8pt/="t2" ; "t1";"t2" **\dir{-} ; ?*\dir{>} ; ?(0)/10pt/*_!/-9pt/{x_{41}} ; ?(1)/-10pt/*_!/8pt/{x_{33}} ;
"b";"c" ; **{} ; ?(0)/8pt/="t1" ; "c";"b" ; **{} ; ?(0)/8pt/="t2" ; "t1";"t2" **\dir{-} ; ?*\dir{>} ; ?(0)/8pt/*_!/-5pt/{x_{31}} ; ?(1)/-8pt/*_!/-5pt/{x_{23}} ;
"e";"b" ; **{} ; ?(0)/8pt/="t1" ; "b";"e" ; **{} ; ?(0)/8pt/="t2" ; "t1";"t2" **\dir{-} ; ?*\dir{>} ; ?(0)/10pt/*_!/8pt/{x_{21}} ; ?(1)/-10pt/*_!/-9pt/{x_{51}} ;
"a"-(0,20);"a"-(0,30) **\dir{-} ; ?(0)*\dir{|} ; ?(1)*\dir{>} ;
"a"-(0,35)*{(\pm)\langle x_{11},x_{22} \rangle\langle x_{12},x_{32} \rangle\langle x_{13},x_{42} \rangle\langle x_{14},x_{52} \rangle\langle x_{21},x_{51} \rangle\langle x_{53},x_{43} \rangle\langle x_{41},x_{33} \rangle\langle x_{31},x_{23} \rangle.}
\end{xy}
\end{displaymath}
\end{example}

\begin{rem} \label{rem_feynman_equivariance}
It follows from Lemma \ref{lem_graph_chord} and equation \eqref{eqn_inv_identity} that the Feynman amplitudes are equivariant with respect to the action of the relevant permutations groups, that is to say that given a fully ordered graph $\Gamma$ of type $(k_1,\ldots,k_m)$ with $k$ edges and a tensor $x \in V^{\otimes k_1}\otimes\ldots\otimes V^{\otimes k_m}$:
\begin{enumerate}
\item
For all $\tau_1,\ldots,\tau_k\in S_2$,
\[ \tilde{F}_{\Gamma\cdot(\tau_1,\ldots,\tau_k)}(x) = \sgn(\tau_1)\cdots\sgn(\tau_k)\tilde{F}_\Gamma(x). \]
\item
For all $\sigma\in S_{k_1}\times\ldots\times S_{k_m}$,
\[ \tilde{F}_{\Gamma\cdot\sigma}(x) = \tilde{F}_\Gamma(\sigma\cdot x).\]
\item
For all $\sigma\in S_m$,
\[ \tilde{F}_{\Gamma\cdot\sigma}(\sigma^{-1}\cdot x) = \tilde{F}_\Gamma(x).\]
\end{enumerate}
\end{rem}

This allows us to give the following definition of the Feynman amplitude for an oriented ribbon graph. Recall that the underlying space of the Lie algebra $\tglie$ is $\bigoplus_{i=3}^\infty (V^{\otimes i})_{Z_i}$, where $V:=\gf^{2n|m}$. This means that we have a map
\[\epsilon\circ T(N):\cek\to T(T^+(V))_{\osp};\]
given an $m$-chain $g_1\wedge\cdots\wedge g_m$, the map $\epsilon\circ T(N)$ is defined by the formula,
\[\epsilon\circ T(N):g_1\wedge\cdots\wedge g_m\mapsto \sum_{\sigma\in S_m}\sgn(\sigma)\sigma\cdot\left[(N\cdot g_1)\otimes\ldots\otimes(N\cdot g_m)\right];\]
where $N$ is the norm operator.

\begin{defi}
Let $\Gamma$ be an oriented ribbon graph and let $\tilde{\Gamma}$ be some fully ordered graph representing it, we define the Feynman amplitude map $F_\Gamma:\cek\to\gf$ by the formula:
\[F_\Gamma(x):=\tilde{F}_{\tilde{\Gamma}}(\epsilon\circ T(N)[x]).\]
\end{defi}

\begin{rem}
It follows from Remark \ref{rem_feynman_equivariance} that the definition of the Feynman amplitude is independent of the choice of representative $\tilde{\Gamma}$. Furthermore, it follows that reversing the orientation on an oriented ribbon graph gives
\[ F_{(\Gamma,-\sigma)}(x) = -F_{(\Gamma,\sigma)}(x).\]
\end{rem}

\subsection{Kontsevich's theorem} \label{sec_kontthm}

In this section we use the Feynman amplitudes defined in the previous section to formulate and prove a super-analogue of Kontsevich's theorem on graph homology. This theorem tells us that the homology of the graph complex may be computed as the relative Lie algebra homology of the Lie algebra $\mathfrak{g}$.

The Feynman amplitudes described in section \ref{sec_Feynman_amplitudes} can be used to define a pairing
\[ \langle\langle-,-\rangle\rangle:\cek\otimes\gc\to\gf \]
between Chevalley-Eilenberg chains and graphs. It is given by the formula,
\begin{equation} \label{eqn_grphpair}
\langle\langle x,\Gamma \rangle\rangle:=\frac{F_\Gamma(x)}{|\Aut(\Gamma)|}.
\end{equation}
The following important theorem says that this pairing is a map of complexes.

\begin{theorem} \label{thm_kontsevich}
The Chevalley-Eilenberg differential $d$ is adjoint to the differential $\delta$ in the graph complex; i.e. for any oriented ribbon graph $\Gamma$ and for all $x \in \cek$,
\begin{equation} \label{eqn_kontsevichdummy0}
\langle\langle dx,\Gamma\rangle\rangle = \langle\langle x,\delta\Gamma\rangle\rangle.
\end{equation}
\end{theorem}

\begin{proof}
Choose a fully ordered graph $\tilde{\Gamma}$ representing $\Gamma$ and write the vertices of $\tilde{\Gamma}$ as
\[ (v_1,v_2,\ldots,v_i). \]
Suppose that $i\in I(\Gamma)$ is an ideal edge of the vertex $v_j$ of $\Gamma$, then we can canonically define a fully ordered graph $\tilde{\Gamma}\backslash i$ which represents $\Gamma\backslash i$. Write the vertex $v_j$ as
\[v_j:=(h_1,\ldots,h_k)\]
and the ideal edge of $v_j$ as $i:=u\cup u'$, where $h_1\in u$; then the vertices of $\tilde{\Gamma}\backslash i$ are ordered as
\[ (v_1,\ldots,v_{j-1},u\sqcup\{a\},u'\sqcup\{b\},v_{j+1},\ldots,v_i) \]
so that the vertex $u\sqcup\{a\}$ containing $h_1$ appears first. The orientation on the edge $e:=\{a,b\}$ depends respectively on whether $j$ is even or odd.

Let $V:=\gf^{2n|m}$ and let $x_1,\ldots,x_l$ be homogeneous tensors in $T_{\geq 3}(V):=\bigoplus_{r=3}^\infty V^{\otimes r}$,
\[ x_i:=x_{i1}\cdots x_{ik_i} \in V^{\otimes k_i}; \]
then the vector $x:=x_1\otimes\ldots\otimes x_l$ represents an $l$-chain of $\cek$. Let us define the vector $x_i\langle r\rangle$ by
\[x_i\langle r\rangle:=x_{i1}\cdots \widehat{x_{ir}}\cdots x_{ik_i}.\]
We now make the following calculations:
\begin{equation} \label{eqn_kontsevichdummy1}
\begin{split}
|\Aut(\Gamma)|\langle\langle x,\delta\Gamma\rangle\rangle &=\sum_{i\in I(\Gamma)}F_{\Gamma\backslash i}(x) \\
&= \sum_{i\in I(\Gamma)}\sum_{\sigma\in S_l}\sum_{\begin{subarray}{c} 1\leq r_1\leq k_1 \\ \vdots \\ 1\leq r_l\leq k_l \end{subarray}}\sgn(\sigma)\tilde{F}_{\tilde{\Gamma}\backslash i}\left[\sigma\cdot(z_{k_1}^{r_1}\cdot x_1\otimes\ldots\otimes z_{k_l}^{r_l}\cdot x_l)\right].
\end{split}
\end{equation}
Using the explicit formula for the Lie bracket given by Lemma \ref{lem_bracket_formula} we get
\begin{equation} \label{eqn_kontsevichdummy2}
\begin{split}
&|\Aut(\Gamma)| \langle\langle dx,\Gamma\rangle\rangle = \sum_{1\leq i < j \leq l} (-1)^p F_\Gamma\left[\{x_i,x_j\}\wedge x_1\wedge\ldots\wedge \hat{x_i}\wedge\ldots\wedge \hat{x_j}\wedge\ldots\wedge x_l\right] \\
&= \sum_{\begin{subarray}{c}1\leq i < j \leq l \\ \sigma \in S_{l-1} \\ 0\leq r \leq k_i+k_j-3 \end{subarray}}\sum_{\begin{subarray}{c} 1\leq r_1\leq k_1 \\ \vdots \\ 1\leq r_l\leq k_l \end{subarray}} (-1)^q\sgn(\sigma)\langle x_{ir_i},x_{jr_j}\rangle \tilde{F}_{\tilde{\Gamma}}\left[\sigma\cdot\left(\begin{array}{c} z_{k_i+k_j-2}^r\cdot\left[\left(z_{k_i-1}^{r_i-1}\cdot x_i\langle r_i\rangle\right)\left(z_{k_j-1}^{r_j-1}\cdot x_j\langle r_j\rangle\right)\right]\otimes \\ (z_{k_1}^{r_1}\cdot x_1)\otimes\ldots\otimes \hat{x_i}\otimes\ldots\otimes \hat{x_j}\otimes\ldots\otimes (z_{k_l}^{r_l}\cdot x_l)\end{array}\right)\right];
\end{split}
\end{equation}
where
\begin{displaymath}
\begin{split}
p&:=|x_i|(|x_1|+\ldots+|x_{i-1}|)+|x_j|(|x_1|+\ldots+|x_{j-1}|)+|x_i||x_j|+i+j-1 \quad \text{and} \\
q&:=p+|x_{ir_i}|(|x_{i,r_i+1}|+\ldots+|x_{ik_i}|+|x_{j1}|+\ldots+|x_{j,r_j-1}|).
\end{split}
\end{displaymath}

In order to establish equation \eqref{eqn_kontsevichdummy0} we must construct a bijection between the terms in the sum \eqref{eqn_kontsevichdummy1} and the terms in the sum \eqref{eqn_kontsevichdummy2}:

Given positive integers $i<j$, let $\tau_{ij}$ be the permutation
\[\tau_{ij}:=(1\,2\ldots i+1)(1\,2\ldots j).\]
Given a permutation $\sigma\in S_{l-1}$, define $\tilde{\sigma}\in S_l$ to be the permutation given by the formula,
\[\tilde{\sigma}:=(\sigma(1)\,\sigma(1)-1\,\ldots 2\,1)(1\, 2\ldots l)\sigma (1\, 2\ldots l)^{-1}.\]

We map the terms in the sum \eqref{eqn_kontsevichdummy2} to the terms in the sum \eqref{eqn_kontsevichdummy1} as follows: given a permutation $\sigma\in S_{l-1}$ and positive integers $i$, $j$ and $r$ such that $1\leq i<j\leq l$ and $1\leq r \leq k_i+k_j-2$, we must construct a permutation $\sigma'\in S_l$ and an ideal edge $\iota$ of $\Gamma$. We may assume that the fully ordered graph $\tilde{\Gamma}\cdot\sigma$ has type $(k_i+k_j-2,k_1,\ldots,\hat{k_i},\ldots,\hat{k_j},\ldots,k_l)$, otherwise the corresponding term in the sum \eqref{eqn_kontsevichdummy2} is zero.

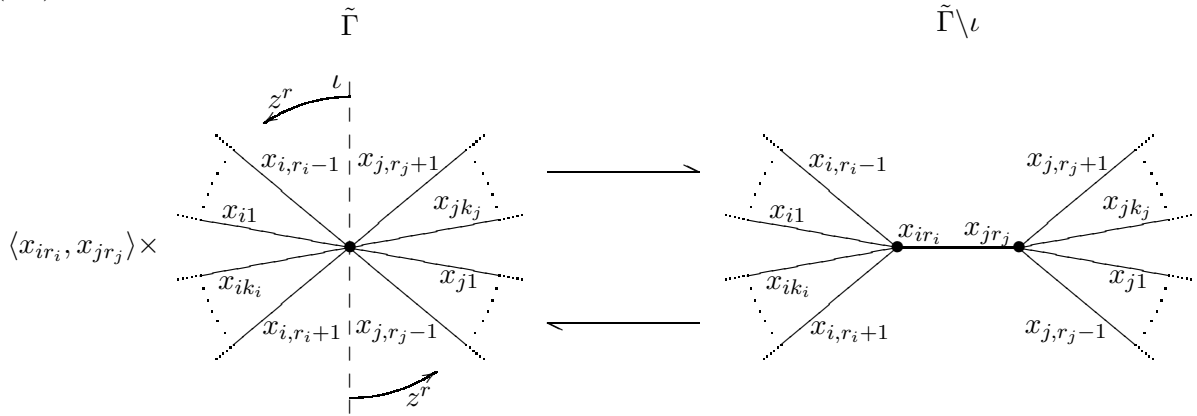
\begin{equation} \label{eqn_feynpair}
\begin{xy}
(-40,0)="a" ;
"a"+a(140);"a" ; **{} ; ?(0)/20mm/="t";"a" ; **\dir{-} ; ?(0)/12mm/*_!/-12pt/{x_{i,r_i-1}} ; "t"+a(140);"t" ; **{} ; ?(0)/3mm/;"t" ; **\dir{.} ;
"a"+a(170);"a" ; **{} ; ?(0)/20mm/="t";"a" ; **\dir{-} ; ?(0)/15mm/*_!/-5pt/{x_{i1}} ; "t"+a(170);"t" ; **{} ; ?(0)/3mm/;"t" ; **\dir{.} ;
"a"+a(190);"a" ; **{} ; ?(0)/20mm/="t";"a" ; **\dir{-} ; ?(0)/15mm/*_!/7pt/{x_{ik_i}} ; "t"+a(190);"t" ; **{} ; ?(0)/3mm/;"t" ; **\dir{.} ;
"a"+a(220);"a" ; **{} ; ?(0)/20mm/="t";"a" ; **\dir{-} ; ?(0)/12mm/*_!/13pt/{x_{i,r_i+1}} ; "t"+a(220);"t" ; **{} ; ?(0)/3mm/;"t" ; **\dir{.} ;
"a"+a(145);"a" ; **{} ; ?(0)/20mm/*\dir{.} ; "a"+a(150);"a" ; **{} ; ?(0)/20mm/*\dir{.} ; "a"+a(155);"a" ; **{} ; ?(0)/20mm/*\dir{.} ; "a"+a(160);"a" ; **{} ; ?(0)/20mm/*\dir{.} ; "a"+a(165);"a" ; **{} ; ?(0)/20mm/*\dir{.} ;
"a"+a(195);"a" ; **{} ; ?(0)/20mm/*\dir{.} ; "a"+a(200);"a" ; **{} ; ?(0)/20mm/*\dir{.} ; "a"+a(205);"a" ; **{} ; ?(0)/20mm/*\dir{.} ; "a"+a(210);"a" ; **{} ; ?(0)/20mm/*\dir{.} ; "a"+a(215);"a" ; **{} ; ?(0)/20mm/*\dir{.} ;
"a"-a(140);"a" ; **{} ; ?(0)/20mm/="t";"a" ; **\dir{-} ; ?(0)/12mm/*_!/-13pt/{x_{j,r_j-1}} ; "t"-a(140);"t" ; **{} ; ?(0)/3mm/;"t" ; **\dir{.} ;
"a"-a(170);"a" ; **{} ; ?(0)/20mm/="t";"a" ; **\dir{-} ; ?(0)/15mm/*_!/-5pt/{x_{j1}} ; "t"-a(170);"t" ; **{} ; ?(0)/3mm/;"t" ; **\dir{.} ;
"a"-a(190);"a" ; **{} ; ?(0)/20mm/="t";"a" ; **\dir{-} ; ?(0)/15mm/*_!/7pt/{x_{jk_j}} ; "t"-a(190);"t" ; **{} ; ?(0)/3mm/;"t" ; **\dir{.} ;
"a"-a(220);"a" ; **{} ; ?(0)/20mm/="t";"a" ; **\dir{-} ; ?(0)/12mm/*_!/12pt/{x_{j,r_j+1}} ; "t"-a(220);"t" ; **{} ; ?(0)/3mm/;"t" ; **\dir{.} ;
"a"-a(145);"a" ; **{} ; ?(0)/20mm/*\dir{.} ; "a"-a(150);"a" ; **{} ; ?(0)/20mm/*\dir{.} ; "a"-a(155);"a" ; **{} ; ?(0)/20mm/*\dir{.} ; "a"-a(160);"a" ; **{} ; ?(0)/20mm/*\dir{.} ; "a"-a(165);"a" ; **{} ; ?(0)/20mm/*\dir{.} ;
"a"-a(195);"a" ; **{} ; ?(0)/20mm/*\dir{.} ; "a"-a(200);"a" ; **{} ; ?(0)/20mm/*\dir{.} ; "a"-a(205);"a" ; **{} ; ?(0)/20mm/*\dir{.} ; "a"-a(210);"a" ; **{} ; ?(0)/20mm/*\dir{.} ; "a"-a(215);"a" ; **{} ; ?(0)/20mm/*\dir{.} ;
"a"+(0,22);"a"-(0,22) ; **\dir{--} ; ?(1)*_!/5pt/{\iota} ;
"a"+a(90);"a" ; **{} ; ?(0)/20mm/="t1" ; "a"+a(95);"a" ; **{} ; ?(0)/20mm/="t2" ; "a"+a(100);"a" ; **{} ; ?(0)/20mm/="t3" ; "a"+a(105);"a" ; **{} ; ?(0)/20mm/="t4" ; "a"+a(110);"a" ; **{} ; ?(0)/20mm/="t5" ; "a"+a(115);"a" ; **{} ; ?(0)/20mm/="t6" ; "a"+a(120);"a" ; **{} ; ?(0)/20mm/="t7" ; "a"+a(125);"a" ; **{} ; ?(0)/20mm/="t8" ;
"t8";"t1" ; **\crv{"t1"&"t2"&"t3"&"t4"&"t5"&"t6"&"t7"&"t8"} ; ?(1)*\dir{>} ; ?(0.7)*_!/-5pt/{z^r} ;
"a"-a(90);"a" ; **{} ; ?(0)/20mm/="t1" ; "a"-a(95);"a" ; **{} ; ?(0)/20mm/="t2" ; "a"-a(100);"a" ; **{} ; ?(0)/20mm/="t3" ; "a"-a(105);"a" ; **{} ; ?(0)/20mm/="t4" ; "a"-a(110);"a" ; **{} ; ?(0)/20mm/="t5" ; "a"-a(115);"a" ; **{} ; ?(0)/20mm/="t6" ; "a"-a(120);"a" ; **{} ; ?(0)/20mm/="t7" ; "a"-a(125);"a" ; **{} ; ?(0)/20mm/="t8" ;
"t8";"t1" ; **\crv{"t1"&"t2"&"t3"&"t4"&"t5"&"t6"&"t7"&"t8"} ; ?(1)*\dir{>} ; ?(0.7)*_!/-5pt/{z^r} ;
"a"*{\bullet} ; "a"-(35,0)*{\langle x_{ir_i},x_{jr_j}\rangle\times} ; "a"+(0,30)*{\tilde{\Gamma}} ;
(40,0)="b" ; "b"-(8,0)="b1" ; "b"+(8,0)="b2" ;
"b2";"b1" ; **\dir{-} ; ?(0)/3mm/*_!/5pt/{x_{ir_i}} ; ?(1)/-4mm/*_!/5pt/{x_{jr_j}} ;
"b1"+a(140);"b1" ; **{} ; ?(0)/20mm/="t";"b1" ; **\dir{-} ; ?(0)/12mm/*_!/-12pt/{x_{i,r_i-1}} ; "t"+a(140);"t" ; **{} ; ?(0)/3mm/;"t" ; **\dir{.} ;
"b1"+a(170);"b1" ; **{} ; ?(0)/20mm/="t";"b1" ; **\dir{-} ; ?(0)/15mm/*_!/-5pt/{x_{i1}} ; "t"+a(170);"t" ; **{} ; ?(0)/3mm/;"t" ; **\dir{.} ;
"b1"+a(190);"b1" ; **{} ; ?(0)/20mm/="t";"b1" ; **\dir{-} ; ?(0)/15mm/*_!/7pt/{x_{ik_i}} ; "t"+a(190);"t" ; **{} ; ?(0)/3mm/;"t" ; **\dir{.} ;
"b1"+a(220);"b1" ; **{} ; ?(0)/20mm/="t";"b1" ; **\dir{-} ; ?(0)/12mm/*_!/13pt/{x_{i,r_i+1}} ; "t"+a(220);"t" ; **{} ; ?(0)/3mm/;"t" ; **\dir{.} ;
"b1"+a(145);"b1" ; **{} ; ?(0)/20mm/*\dir{.} ; "b1"+a(150);"b1" ; **{} ; ?(0)/20mm/*\dir{.} ; "b1"+a(155);"b1" ; **{} ; ?(0)/20mm/*\dir{.} ; "b1"+a(160);"b1" ; **{} ; ?(0)/20mm/*\dir{.} ; "b1"+a(165);"b1" ; **{} ; ?(0)/20mm/*\dir{.} ;
"b1"+a(195);"b1" ; **{} ; ?(0)/20mm/*\dir{.} ; "b1"+a(200);"b1" ; **{} ; ?(0)/20mm/*\dir{.} ; "b1"+a(205);"b1" ; **{} ; ?(0)/20mm/*\dir{.} ; "b1"+a(210);"b1" ; **{} ; ?(0)/20mm/*\dir{.} ; "b1"+a(215);"b1" ; **{} ; ?(0)/20mm/*\dir{.} ;
"b2"-a(140);"b2" ; **{} ; ?(0)/20mm/="t";"b2" ; **\dir{-} ; ?(0)/12mm/*_!/-13pt/{x_{j,r_j-1}} ; "t"-a(140);"t" ; **{} ; ?(0)/3mm/;"t" ; **\dir{.} ;
"b2"-a(170);"b2" ; **{} ; ?(0)/20mm/="t";"b2" ; **\dir{-} ; ?(0)/15mm/*_!/-5pt/{x_{j1}} ; "t"-a(170);"t" ; **{} ; ?(0)/3mm/;"t" ; **\dir{.} ;
"b2"-a(190);"b2" ; **{} ; ?(0)/20mm/="t";"b2" ; **\dir{-} ; ?(0)/15mm/*_!/7pt/{x_{jk_j}} ; "t"-a(190);"t" ; **{} ; ?(0)/3mm/;"t" ; **\dir{.} ;
"b2"-a(220);"b2" ; **{} ; ?(0)/20mm/="t";"b2" ; **\dir{-} ; ?(0)/12mm/*_!/12pt/{x_{j,r_j+1}} ; "t"-a(220);"t" ; **{} ; ?(0)/3mm/;"t" ; **\dir{.} ;
"b2"-a(145);"b2" ; **{} ; ?(0)/20mm/*\dir{.} ; "b2"-a(150);"b2" ; **{} ; ?(0)/20mm/*\dir{.} ; "b2"-a(155);"b2" ; **{} ; ?(0)/20mm/*\dir{.} ; "b2"-a(160);"b2" ; **{} ; ?(0)/20mm/*\dir{.} ; "b2"-a(165);"b2" ; **{} ; ?(0)/20mm/*\dir{.} ;
"b2"-a(195);"b2" ; **{} ; ?(0)/20mm/*\dir{.} ; "b2"-a(200);"b2" ; **{} ; ?(0)/20mm/*\dir{.} ; "b2"-a(205);"b2" ; **{} ; ?(0)/20mm/*\dir{.} ; "b2"-a(210);"b2" ; **{} ; ?(0)/20mm/*\dir{.} ; "b2"-a(215);"b2" ; **{} ; ?(0)/20mm/*\dir{.} ;
"b1"*{\bullet} ; "b2"*{\bullet} ; "b"+(0,30)*{\tilde{\Gamma}\backslash \iota} ;
(6,10);(-14,10) ; **\dir{-} ; (6,10);(4,11) ; **\dir{-} ;
(6,-10);(-14,-10) ; **\dir{-} ; (-14,-10);(-12,-11) ; **\dir{-} ;
\end{xy}
\end{equation}

Let $v:=(h_1,\ldots,h_{k_i+k_j-2})$ be the $\sigma(1)$th vertex of $\tilde{\Gamma}$, then we define an ideal edge $\iota$ of the vertex $v$ of $\Gamma$ by,
\[\iota:=\left\{ \big(h_{z^r(1)},\ldots,h_{z^r(k_i-1)}\big) , \big(h_{z^r(k_i)},\ldots,h_{z^r(k_i+k_j-2)}\big) \right\};\]
where $z:=z_{k_i+k_j-2}$.

We define the permutation $\sigma'\in S_l$ by the formula:
\begin{equation} \label{eqn_kontsevichdummy3}
\sigma':=\left\{\begin{array}{ll} \tilde{\sigma}\tau_{ij}, & 0\leq r \leq k_i-2; \\ \tilde{\sigma}\tau_{ij}(i\,j), & k_i-1\leq r \leq k_i+k_j-3. \end{array}\right.
\end{equation}

It is straightforward to check that the terms in the sum \eqref{eqn_kontsevichdummy2} are mapped bijectively to the terms in the sum \eqref{eqn_kontsevichdummy1} under the above correspondence.
\end{proof}

\begin{rem}
The permutation $\sigma'$ defined by equation \eqref{eqn_kontsevichdummy3} is just the permutation which takes the letters $i$ and $j$ to the letters $\sigma(1)$ and $\sigma(1)+1$ and moves the remaining letters according to the permutation $\sigma$. Note that the appearance of the permutation $(i\, j)$ in \eqref{eqn_kontsevichdummy3} is due, essentially, to the possible symmetry arising in diagram \eqref{eqn_feynpair} when the ideal edge splits the vertex into two \emph{equal} parts.
\end{rem}

Given a chord diagram $c \in \chd{k}$ written as
\[ c:= \{i_1,j_1\},\ldots,\{i_k,j_k\}, \]
we can canonically define an \emph{oriented} chord diagram $\hat{c}\in\ochd{k}$ as follows: we may assume that $i_r < j_r$ for all $r$, then $\hat{c}$ is defined as,
\[ \hat{c}:= (i_1,j_1),\ldots,(i_k,j_k). \]

We can define a map $I:\cek\to\gc$ which formally resembles a kind of integration in which certain contributions to the integral are indexed by oriented ribbon graphs:

\begin{defi} \label{def_integral}
Let
\[x:= (x_{11}\cdots x_{1k_1})\otimes\ldots\otimes(x_{l1}\cdots x_{lk_l})\]
represent a Chevalley-Eilenberg chain, where $x_{ij}\in V:=\gf^{2n|m}$ and $k_1+\ldots+k_l=2k$. The map $I:\cek\to\gc$ is defined by the formula,
\[ I(x):= \sum_{c\in\chd{k}} \beta_{\hat{c}}(x)\grpchd{\hat{c}}.\]
\end{defi}

\begin{rem}
It follows from Lemma \ref{lem_graph_chord} and equation \eqref{eqn_inv_identity} that $I$ is a well defined map modulo the symmetry relations present in $\cek$. Furthermore, it will follow from theorems \ref{thm_kontsevich} and \ref{thm_Feynman_spinvariants} that $I$ is a map of complexes.
\end{rem}

\begin{theorem} \label{thm_Feynman_spinvariants}
The following diagram commutes:
\[\xymatrix{\cek \ar^{x\mapsto\langle\langle x,- \rangle\rangle}[rr] \ar_{I}[rd] && \Hom(\gc,\gf) \\ & \gc \ar_{\Gamma\mapsto\langle\Gamma,-\rangle}[ru]}\]
\end{theorem}

\begin{proof}
Let $\Gamma'$ be an oriented ribbon graph and let $x:= x_1\otimes\ldots\otimes x_l$ represent a Chevalley-Eilenberg chain, where $x_i\in V^{\otimes k_i}$ and $k_1+\ldots+k_l=2k$. We must show that
\begin{equation} \label{eqn_Feynman_spinvariants}
\langle\langle x,\Gamma' \rangle\rangle = \langle I(x),\Gamma' \rangle.
\end{equation}
Let $\tilde{\Gamma}'$ be a fully ordered graph representing $\Gamma'$. We may assume that $\tilde{\Gamma}'$ has type $(k_1,\ldots,k_l)$, otherwise both sides of \eqref{eqn_Feynman_spinvariants} are zero. It follows from the relevant definitions that \eqref{eqn_Feynman_spinvariants} is equivalent to
\[\sum_{\sigma\in S_l}\sum_{\begin{subarray}{c} \gamma_1\in Z_{k_1} \\ \vdots \\ \gamma_l\in Z_{k_l}\end{subarray}} \sgn(\sigma) \tilde{F}_{\tilde{\Gamma}'}\left[\sigma\cdot(\gamma_1\cdot x_1\otimes\ldots\otimes \gamma_l\cdot x_l)\right] = |\Aut(\Gamma')|\sum_{c\in\chd{k}}\beta_{\hat{c}}(x)\langle\Gamma_{k_1,\ldots,k_l}(\hat{c}),\Gamma'\rangle.\]
We may assume that all the permutations $\sigma\in S_l$ in the sum on the left hand side satisfy
\begin{equation} \label{eqn_Feynman_spinvariantsdummy0}
(k_{\sigma(1)},\ldots,k_{\sigma(l)}) = (k_1,\ldots,k_l),
\end{equation}
otherwise the terms in the sum are zero. Likewise, we may assume that all the chord diagrams $c\in\chd{k}$ in the sum on the right hand side are such that the ribbon graphs $\grpchd{\hat{c}}$ and $\Gamma'$ are isomorphic (disregarding orientation).

Let us denote the subgroup of $S_l$ consisting of permutations satisfying \eqref{eqn_Feynman_spinvariantsdummy0} by $\tilde{S}_l$. We define the group $H$ as the wreath product
\[H:=\tilde{S}_l\ltimes\left(Z_{k_1}\times\ldots\times Z_{k_l}\right).\]
Recall (cf. figure \eqref{eqn_perm_embedding}) that $S_l$ is embedded inside $S_{2k}$;
\begin{displaymath}
\begin{array}{ccc}
S_l & \hookrightarrow & S_{2k}, \\
\sigma & \mapsto & \sigma_{k_1,\ldots,k_l}.
\end{array}
\end{displaymath}
Obviously $Z_{k_1}\times\ldots\times Z_{k_l}$ is also embedded in $S_{2k}$ and this gives us a natural embedding $H\hookrightarrow S_{2k}$. Furthermore, there is an action of $H$ on the left of $V^{\otimes k_1}\otimes\ldots\otimes V^{\otimes k_l}$ which is defined in an obvious manner and this action is compatible with the embedding of $H$ into $S_{2k}$.

If we write the vertices of the fully ordered graph $\tilde{\Gamma}'$ as
\[(h_1,\ldots,h_{k_1})\quad(h_{k_1+1},\ldots,h_{k_1+k_2})\quad\ldots\quad(h_{k_1+\ldots+k_{l-1}+1},\ldots,h_{k_1+\ldots+k_l})\]
then every permutation $\sigma\in S_{2k}$ gives rise to a permutation of the half-edges of $\tilde{\Gamma}'$,
\begin{equation} \label{eqn_Feynman_spinvariantsdummy1}
h_i\mapsto h_{\sigma(i)}.
\end{equation}

Let $K$ be the subgroup of $H$ consisting of all permutations $(\sigma,\gamma)$ such that
\[(\sigma,\gamma)\cdot \bar{c}_{\tilde{\Gamma}'} = \bar{c}_{\tilde{\Gamma}'},\]
where $\bar{c}_{\tilde{\Gamma}'}$ is the \emph{unoriented} chord diagram corresponding to $\tilde{\Gamma}'$; that is to say that $K$ is the stabiliser of the point $\bar{c}_{\tilde{\Gamma}'}\in\chd{k}$. Since we may assume that every automorphism of $\Gamma'$ preserves the orientation on $\Gamma'$, it follows from Lemma \ref{lem_graph_chord} that \eqref{eqn_Feynman_spinvariantsdummy1} gives an embedding
\[K\to\Aut(\Gamma').\]

\begin{displaymath}
\begin{xy}
(0,0)="a" ;
"a"+a(90);"a" ; **{} ; ?(0)/25mm/="a1" ; ?(0)/13mm/*{p} ;
"a"+a(18);"a" ; **{} ; ?(0)/25mm/="a2" ; ?(0)/13mm/*{q} ;
"a"+a(306);"a" ; **{} ; ?(0)/25mm/="a3" ; ?(0)/13mm/*{r} ;
"a"+a(234);"a" ; **{} ; ?(0)/25mm/="a4" ; ?(0)/13mm/*{s} ;
"a"+a(162);"a" ; **{} ; ?(0)/25mm/="a5" ; ?(0)/13mm/*{t} ;
"a1"+a(270);"a1" ; **{} ; ?(0)/6mm/="t";"a1" ; **\dir{-} ; ?(1)*_!/4pt/{1} ; "t"+a(270);"t" ; **{} ; ?(0)/1mm/;"t" ; **\dir{.} ;
"a1"+a(225);"a1" ; **{} ; ?(0)/6mm/="t";"a1" ; **\dir{-} ; ?(1)*_!/5pt/{2} ; "t"+a(225);"t" ; **{} ; ?(0)/1mm/;"t" ; **\dir{.} ;
"a1"+a(210);"a1" ; **{} ; ?(0)/4mm/*\dir{.} ; "a1"+a(190);"a1" ; **{} ; ?(0)/4mm/*\dir{.} ; "a1"+a(170);"a1" ; **{} ; ?(0)/4mm/*\dir{.} ;
"a1"+a(155);"a1" ; **{} ; ?(0)/6mm/="t";"a1" ; **\dir{-} ; ?(1)*_!/5pt/{i_p} ; "t"+a(155);"t" ; **{} ; ?(0)/1mm/;"t" ; **\dir{.} ;
"a1"+a(140);"a1" ; **{} ; ?(0)/4mm/*\dir{.} ; "a1"+a(120);"a1" ; **{} ; ?(0)/4mm/*\dir{.} ; "a1"+a(100);"a1" ; **{} ; ?(0)/4mm/*\dir{.} ; "a1"+a(80);"a1" ; **{} ; ?(0)/4mm/*\dir{.} ;
"a1"+a(140);"a1" ; **{} ; ?(0)/6mm/="t1" ; "a1"+a(120);"a1" ; **{} ; ?(0)/6mm/="t2" ; "a1"+a(100);"a1" ; **{} ; ?(0)/6mm/="t3" ; "a1"+a(80);"a1" ; **{} ; ?(0)/6mm/="t4" ;
"t4";"t1" ; **\crv{"t1"&"t2"&"t3"&"t4"} ; ?(1)*\dir{>} ; ?*_!/5pt/{\gamma_p} ;
"a1"+a(45);"a1" ; **{} ; ?(0)/6mm/="t";"a1" ; **\dir{-} ; ?(1)*_!/5pt/{j_p} ; "t"+a(45);"t" ; **{} ; ?(0)/1mm/;"t" ; **\dir{.} ;
"a1"+a(30);"a1" ; **{} ; ?(0)/4mm/*\dir{.} ; "a1"+a(10);"a1" ; **{} ; ?(0)/4mm/*\dir{.} ;
"a1"+a(315);"a1" ; **{} ; ?(0)/6mm/="t";"a1" ; **\dir{-} ; ?(1)*_!/6pt/{k_p} ; "t"+a(315);"t" ; **{} ; ?(0)/1mm/;"t" ; **\dir{.} ;
"a1"*{\bullet} ; "a1"*\xycircle(10,10){} ;
"a2"+a(198);"a2" ; **{} ; ?(0)/6mm/="t";"a2" ; **\dir{-} ; ?(1)*_!/5pt/{1} ; "t"+a(198);"t" ; **{} ; ?(0)/1mm/;"t" ; **\dir{.} ;
"a2"+a(153);"a2" ; **{} ; ?(0)/6mm/="t";"a2" ; **\dir{-} ; ?(1)*_!/5pt/{2} ; "t"+a(153);"t" ; **{} ; ?(0)/1mm/;"t" ; **\dir{.} ;
"a2"+a(138);"a2" ; **{} ; ?(0)/4mm/*\dir{.} ; "a2"+a(118);"a2" ; **{} ; ?(0)/4mm/*\dir{.} ; "a2"+a(98);"a2" ; **{} ; ?(0)/4mm/*\dir{.} ;
"a2"+a(83);"a2" ; **{} ; ?(0)/6mm/="t";"a2" ; **\dir{-} ; ?(1)*_!/4pt/{i_q} ; "t"+a(83);"t" ; **{} ; ?(0)/1mm/;"t" ; **\dir{.} ;
"a2"+a(68);"a2" ; **{} ; ?(0)/4mm/*\dir{.} ; "a2"+a(48);"a2" ; **{} ; ?(0)/4mm/*\dir{.} ; "a2"+a(28);"a2" ; **{} ; ?(0)/4mm/*\dir{.} ; "a2"+a(8);"a2" ; **{} ; ?(0)/4mm/*\dir{.} ;
"a2"+a(68);"a2" ; **{} ; ?(0)/6mm/="t1" ; "a2"+a(48);"a2" ; **{} ; ?(0)/6mm/="t2" ; "a2"+a(28);"a2" ; **{} ; ?(0)/6mm/="t3" ; "a2"+a(8);"a2" ; **{} ; ?(0)/6mm/="t4" ;
"t4";"t1" ; **\crv{"t1"&"t2"&"t3"&"t4"} ; ?(1)*\dir{>} ; ?*_!/5pt/{\gamma_q} ;
"a2"+a(333);"a2" ; **{} ; ?(0)/6mm/="t";"a2" ; **\dir{-} ; ?(1)/2pt/*_!/7pt/{j_q} ; "t"+a(333);"t" ; **{} ; ?(0)/1mm/;"t" ; **\dir{.} ;
"a2"+a(318);"a2" ; **{} ; ?(0)/4mm/*\dir{.} ; "a2"+a(298);"a2" ; **{} ; ?(0)/4mm/*\dir{.} ;
"a2"+a(243);"a2" ; **{} ; ?(0)/6mm/="t";"a2" ; **\dir{-} ; ?(1)*_!/7pt/{k_q} ; "t"+a(243);"t" ; **{} ; ?(0)/1mm/;"t" ; **\dir{.} ;
"a2"*{\bullet} ; "a2"*\xycircle(10,10){} ;
"a3"+a(126);"a3" ; **{} ; ?(0)/6mm/="t";"a3" ; **\dir{-} ; ?(1)*_!/4pt/{1} ; "t"+a(126);"t" ; **{} ; ?(0)/1mm/;"t" ; **\dir{.} ;
"a3"+a(81);"a3" ; **{} ; ?(0)/6mm/="t";"a3" ; **\dir{-} ; ?(1)*_!/4pt/{2} ; "t"+a(81);"t" ; **{} ; ?(0)/1mm/;"t" ; **\dir{.} ;
"a3"+a(66);"a3" ; **{} ; ?(0)/4mm/*\dir{.} ; "a3"+a(46);"a3" ; **{} ; ?(0)/4mm/*\dir{.} ; "a3"+a(26);"a3" ; **{} ; ?(0)/4mm/*\dir{.} ;
"a3"+a(11);"a3" ; **{} ; ?(0)/6mm/="t";"a3" ; **\dir{-} ; ?(1)*_!/5pt/{i_r} ; "t"+a(11);"t" ; **{} ; ?(0)/1mm/;"t" ; **\dir{.} ;
"a3"+a(356);"a3" ; **{} ; ?(0)/4mm/*\dir{.} ; "a3"+a(336);"a3" ; **{} ; ?(0)/4mm/*\dir{.} ; "a3"+a(316);"a3" ; **{} ; ?(0)/4mm/*\dir{.} ; "a3"+a(296);"a3" ; **{} ; ?(0)/4mm/*\dir{.} ;
"a3"+a(356);"a3" ; **{} ; ?(0)/6mm/="t1" ; "a3"+a(336);"a3" ; **{} ; ?(0)/6mm/="t2" ; "a3"+a(316);"a3" ; **{} ; ?(0)/6mm/="t3" ; "a3"+a(296);"a3" ; **{} ; ?(0)/6mm/="t4" ;
"t4";"t1" ; **\crv{"t1"&"t2"&"t3"&"t4"} ; ?(1)*\dir{>} ; ?*_!/6pt/{\gamma_r} ;
"a3"+a(261);"a3" ; **{} ; ?(0)/6mm/="t";"a3" ; **\dir{-} ; ?(1)*_!/5pt/{j_r} ; "t"+a(261);"t" ; **{} ; ?(0)/1mm/;"t" ; **\dir{.} ;
"a3"+a(246);"a3" ; **{} ; ?(0)/4mm/*\dir{.} ; "a3"+a(226);"a3" ; **{} ; ?(0)/4mm/*\dir{.} ;
"a3"+a(171);"a3" ; **{} ; ?(0)/6mm/="t";"a3" ; **\dir{-} ; ?(1)*_!/6pt/{k_r} ; "t"+a(171);"t" ; **{} ; ?(0)/1mm/;"t" ; **\dir{.} ;
"a3"*{\bullet} ; "a3"*\xycircle(10,10){} ;
"a4"+a(54);"a4" ; **{} ; ?(0)/6mm/="t";"a4" ; **\dir{-} ; ?(1)*_!/5pt/{1} ; "t"+a(54);"t" ; **{} ; ?(0)/1mm/;"t" ; **\dir{.} ;
"a4"+a(9);"a4" ; **{} ; ?(0)/6mm/="t";"a4" ; **\dir{-} ; ?(1)*_!/5pt/{2} ; "t"+a(9);"t" ; **{} ; ?(0)/1mm/;"t" ; **\dir{.} ;
"a4"+a(354);"a4" ; **{} ; ?(0)/4mm/*\dir{.} ; "a4"+a(334);"a4" ; **{} ; ?(0)/4mm/*\dir{.} ; "a4"+a(314);"a4" ; **{} ; ?(0)/4mm/*\dir{.} ;
"a4"+a(299);"a4" ; **{} ; ?(0)/6mm/="t";"a4" ; **\dir{-} ; ?(1)*_!/5pt/{i_s} ; "t"+a(299);"t" ; **{} ; ?(0)/1mm/;"t" ; **\dir{.} ;
"a4"+a(284);"a4" ; **{} ; ?(0)/4mm/*\dir{.} ; "a4"+a(264);"a4" ; **{} ; ?(0)/4mm/*\dir{.} ; "a4"+a(244);"a4" ; **{} ; ?(0)/4mm/*\dir{.} ; "a4"+a(224);"a4" ; **{} ; ?(0)/4mm/*\dir{.} ;
"a4"+a(284);"a4" ; **{} ; ?(0)/6mm/="t1" ; "a4"+a(264);"a4" ; **{} ; ?(0)/6mm/="t2" ; "a4"+a(244);"a4" ; **{} ; ?(0)/6mm/="t3" ; "a4"+a(224);"a4" ; **{} ; ?(0)/6mm/="t4" ;
"t4";"t1" ; **\crv{"t1"&"t2"&"t3"&"t4"} ; ?(1)*\dir{>} ; ?*_!/6pt/{\gamma_s} ;
"a4"+a(189);"a4" ; **{} ; ?(0)/6mm/="t";"a4" ; **\dir{-} ; ?(1)*_!/6pt/{j_s} ; "t"+a(189);"t" ; **{} ; ?(0)/1mm/;"t" ; **\dir{.} ;
"a4"+a(174);"a4" ; **{} ; ?(0)/4mm/*\dir{.} ; "a4"+a(154);"a4" ; **{} ; ?(0)/4mm/*\dir{.} ;
"a4"+a(99);"a4" ; **{} ; ?(0)/6mm/="t";"a4" ; **\dir{-} ; ?(1)*_!/6pt/{k_s} ; "t"+a(99);"t" ; **{} ; ?(0)/1mm/;"t" ; **\dir{.} ;
"a4"*{\bullet} ; "a4"*\xycircle(10,10){} ;
"a5"+a(342);"a5" ; **{} ; ?(0)/6mm/="t";"a5" ; **\dir{-} ; ?(1)*_!/5pt/{1} ; "t"+a(342);"t" ; **{} ; ?(0)/1mm/;"t" ; **\dir{.} ;
"a5"+a(297);"a5" ; **{} ; ?(0)/6mm/="t";"a5" ; **\dir{-} ; ?(1)*_!/5pt/{2} ; "t"+a(297);"t" ; **{} ; ?(0)/1mm/;"t" ; **\dir{.} ;
"a5"+a(282);"a5" ; **{} ; ?(0)/4mm/*\dir{.} ; "a5"+a(262);"a5" ; **{} ; ?(0)/4mm/*\dir{.} ; "a5"+a(242);"a5" ; **{} ; ?(0)/4mm/*\dir{.} ;
"a5"+a(227);"a5" ; **{} ; ?(0)/6mm/="t";"a5" ; **\dir{-} ; ?(1)*_!/5pt/{i_t} ; "t"+a(227);"t" ; **{} ; ?(0)/1mm/;"t" ; **\dir{.} ;
"a5"+a(212);"a5" ; **{} ; ?(0)/4mm/*\dir{.} ; "a5"+a(192);"a5" ; **{} ; ?(0)/4mm/*\dir{.} ; "a5"+a(172);"a5" ; **{} ; ?(0)/4mm/*\dir{.} ; "a5"+a(152);"a5" ; **{} ; ?(0)/4mm/*\dir{.} ;
"a5"+a(212);"a5" ; **{} ; ?(0)/6mm/="t1" ; "a5"+a(192);"a5" ; **{} ; ?(0)/6mm/="t2" ; "a5"+a(172);"a5" ; **{} ; ?(0)/6mm/="t3" ; "a5"+a(152);"a5" ; **{} ; ?(0)/6mm/="t4" ;
"t4";"t1" ; **\crv{"t1"&"t2"&"t3"&"t4"} ; ?(1)*\dir{>} ; ?*_!/6pt/{\gamma_t} ;
"a5"+a(117);"a5" ; **{} ; ?(0)/6mm/="t";"a5" ; **\dir{-} ; ?(1)/1pt/*_!/4pt/{j_t} ; "t"+a(117);"t" ; **{} ; ?(0)/1mm/;"t" ; **\dir{.} ;
"a5"+a(102);"a5" ; **{} ; ?(0)/4mm/*\dir{.} ; "a5"+a(82);"a5" ; **{} ; ?(0)/4mm/*\dir{.} ;
"a5"+a(27);"a5" ; **{} ; ?(0)/6mm/="t";"a5" ; **\dir{-} ; ?(1)*_!/7pt/{k_t} ; "t"+a(27);"t" ; **{} ; ?(0)/1mm/;"t" ; **\dir{.} ;
"a5"*{\bullet} ; "a5"*\xycircle(10,10){} ;
"a"+a(48);"a" ; **{} ; ?(0)/25mm/*\dir{.} ; "a"+a(51);"a" ; **{} ; ?(0)/25mm/*\dir{.} ; "a"+a(54);"a" ; **{} ; ?(0)/25mm/*\dir{.} ; "a"+a(57);"a" ; **{} ; ?(0)/25mm/*\dir{.} ;"a"+a(60);"a" ; **{} ; ?(0)/25mm/*\dir{.} ;
"a"+a(120);"a" ; **{} ; ?(0)/25mm/*\dir{.} ; "a"+a(123);"a" ; **{} ; ?(0)/25mm/*\dir{.} ; "a"+a(126);"a" ; **{} ; ?(0)/25mm/*\dir{.} ; "a"+a(129);"a" ; **{} ; ?(0)/25mm/*\dir{.} ;"a"+a(132);"a" ; **{} ; ?(0)/25mm/*\dir{.} ;
"a"+a(192);"a" ; **{} ; ?(0)/25mm/*\dir{.} ; "a"+a(195);"a" ; **{} ; ?(0)/25mm/*\dir{.} ; "a"+a(198);"a" ; **{} ; ?(0)/25mm/*\dir{.} ; "a"+a(201);"a" ; **{} ; ?(0)/25mm/*\dir{.} ;"a"+a(204);"a" ; **{} ; ?(0)/25mm/*\dir{.} ;
"a"+a(264);"a" ; **{} ; ?(0)/25mm/*\dir{.} ; "a"+a(267);"a" ; **{} ; ?(0)/25mm/*\dir{.} ; "a"+a(270);"a" ; **{} ; ?(0)/25mm/*\dir{.} ; "a"+a(273);"a" ; **{} ; ?(0)/25mm/*\dir{.} ;"a"+a(276);"a" ; **{} ; ?(0)/25mm/*\dir{.} ;
"a"+a(336);"a" ; **{} ; ?(0)/25mm/*\dir{.} ; "a"+a(339);"a" ; **{} ; ?(0)/25mm/*\dir{.} ; "a"+a(342);"a" ; **{} ; ?(0)/25mm/*\dir{.} ; "a"+a(345);"a" ; **{} ; ?(0)/25mm/*\dir{.} ;"a"+a(348);"a" ; **{} ; ?(0)/25mm/*\dir{.} ;
"a"*[o]=<5mm>{\sigma} *\frm{o} ;
"a"+a(18);"a" ; **{} ; ?(0)/2.5mm/="t" ; "a"+a(18);"a" ; **{} ; ?(0)/10mm/;"t" ; **\dir{-} ; ?*\dir{>} ;
"a"+a(90);"a" ; **{} ; ?(0)/2.5mm/="t" ; "a"+a(90);"a" ; **{} ; ?(0)/10mm/;"t" ; **\dir{-} ; ?*\dir{<} ;
"a"+a(162);"a" ; **{} ; ?(0)/2.5mm/="t" ; "a"+a(162);"a" ; **{} ; ?(0)/10mm/;"t" ; **\dir{-} ; ?*\dir{>} ;
"a"+a(234);"a" ; **{} ; ?(0)/2.5mm/="t" ; "a"+a(234);"a" ; **{} ; ?(0)/10mm/;"t" ; **\dir{-} ; ?*\dir{>} ;
"a"+a(306);"a" ; **{} ; ?(0)/2.5mm/="t" ; "a"+a(306);"a" ; **{} ; ?(0)/10mm/;"t" ; **\dir{-} ; ?*\dir{<} ;
"a"+a(54);"a" ; **{} ; ?(0)/2.5mm/="t1" ; "a"+a(54);"a" ; **{} ; ?(0)/10mm/="t2";"t1" ; **\dir{-} ; ?*\dir{>} ; "t2"+a(54);"t2" ; **{} ; ?(0)/2mm/;"t2" ; **\dir{.} ;
"a"+a(30);"a" ; **{} ; ?(0)/7mm/*\dir{.} ; "a"+a(36);"a" ; **{} ; ?(0)/7mm/*\dir{.} ; "a"+a(42);"a" ; **{} ; ?(0)/7mm/*\dir{.} ;
"a"+a(66);"a" ; **{} ; ?(0)/7mm/*\dir{.} ; "a"+a(72);"a" ; **{} ; ?(0)/7mm/*\dir{.} ; "a"+a(78);"a" ; **{} ; ?(0)/7mm/*\dir{.} ;
"a"+a(126);"a" ; **{} ; ?(0)/2.5mm/="t1" ; "a"+a(126);"a" ; **{} ; ?(0)/10mm/="t2";"t1" ; **\dir{-} ; ?*\dir{<} ; "t2"+a(126);"t2" ; **{} ; ?(0)/2mm/;"t2" ; **\dir{.} ;
"a"+a(102);"a" ; **{} ; ?(0)/7mm/*\dir{.} ; "a"+a(108);"a" ; **{} ; ?(0)/7mm/*\dir{.} ; "a"+a(114);"a" ; **{} ; ?(0)/7mm/*\dir{.} ;
"a"+a(138);"a" ; **{} ; ?(0)/7mm/*\dir{.} ; "a"+a(144);"a" ; **{} ; ?(0)/7mm/*\dir{.} ; "a"+a(150);"a" ; **{} ; ?(0)/7mm/*\dir{.} ;
"a"+a(198);"a" ; **{} ; ?(0)/2.5mm/="t1" ; "a"+a(198);"a" ; **{} ; ?(0)/10mm/="t2";"t1" ; **\dir{-} ; ?*\dir{<} ; "t2"+a(198);"t2" ; **{} ; ?(0)/2mm/;"t2" ; **\dir{.} ;
"a"+a(174);"a" ; **{} ; ?(0)/7mm/*\dir{.} ; "a"+a(180);"a" ; **{} ; ?(0)/7mm/*\dir{.} ; "a"+a(186);"a" ; **{} ; ?(0)/7mm/*\dir{.} ;
"a"+a(210);"a" ; **{} ; ?(0)/7mm/*\dir{.} ; "a"+a(216);"a" ; **{} ; ?(0)/7mm/*\dir{.} ; "a"+a(222);"a" ; **{} ; ?(0)/7mm/*\dir{.} ;
"a"+a(270);"a" ; **{} ; ?(0)/2.5mm/="t1" ; "a"+a(270);"a" ; **{} ; ?(0)/10mm/="t2";"t1" ; **\dir{-} ; ?*\dir{>} ; "t2"+a(270);"t2" ; **{} ; ?(0)/2mm/;"t2" ; **\dir{.} ;
"a"+a(246);"a" ; **{} ; ?(0)/7mm/*\dir{.} ; "a"+a(252);"a" ; **{} ; ?(0)/7mm/*\dir{.} ; "a"+a(258);"a" ; **{} ; ?(0)/7mm/*\dir{.} ;
"a"+a(282);"a" ; **{} ; ?(0)/7mm/*\dir{.} ; "a"+a(288);"a" ; **{} ; ?(0)/7mm/*\dir{.} ; "a"+a(294);"a" ; **{} ; ?(0)/7mm/*\dir{.} ;
"a"+a(342);"a" ; **{} ; ?(0)/2.5mm/="t1" ; "a"+a(342);"a" ; **{} ; ?(0)/10mm/="t2";"t1" ; **\dir{-} ; ?*\dir{>} ; "t2"+a(342);"t2" ; **{} ; ?(0)/2mm/;"t2" ; **\dir{.} ;
"a"+a(318);"a" ; **{} ; ?(0)/7mm/*\dir{.} ; "a"+a(324);"a" ; **{} ; ?(0)/7mm/*\dir{.} ; "a"+a(330);"a" ; **{} ; ?(0)/7mm/*\dir{.} ;
"a"+a(354);"a" ; **{} ; ?(0)/7mm/*\dir{.} ; "a"+a(0);"a" ; **{} ; ?(0)/7mm/*\dir{.} ; "a"+a(6);"a" ; **{} ; ?(0)/7mm/*\dir{.} ;
"a"+(0,+40)*{\Gamma'\to\Gamma':h_i\mapsto h_{(\sigma,\gamma)[i]}} ;
\end{xy}
\end{displaymath}

One soon realises that every automorphism is obtained in this way, that is to say that this map is an isomorphism of groups. Every automorphism of $\Gamma'$ must preserve the vertices of $\Gamma'$, hence it is a bijection on the vertices which is described by some permutation $\sigma\in \tilde{S}_l$. Furthermore, this automorphism must preserve the cyclic ordering at each vertex of $\Gamma'$ and hence is described locally at every vertex $v_i$ by some permutation $\gamma_i\in Z_{k_i}$. Since the automorphism should preserve the edges of $\Gamma'$, it follows from Lemma \ref{lem_graph_chord} that the permutation $(\sigma,\gamma)\in K$. This permutation is the member of $K$ from which the automorphism of $\Gamma'$ is obtained.

It follows from Lemma \ref{lem_graph_chord} that if $c\in\chd{k}$ lies in the orbit of the point $\bar{c}_{\tilde{\Gamma}'}$ under the action of the group $H$, then the ribbon graphs $\grpchd{\hat{c}}$ and $\Gamma'$ are isomorphic (disregarding orientation). For the same reasons as stated in the preceding paragraph, the converse is also true; that is to say that
\[\orb\left(\bar{c}_{\tilde{\Gamma}'}\right)=\{c\in\chd{k}: \grpchd{\hat{c}}\text{ and }\Gamma'\text{ are isomorphic ribbon graphs}\}.\]
Briefly, every ribbon graph isomorphism $\grpchd{\hat{c}}\to\Gamma'$ is described by some permutation $\sigma\in\tilde{S}_l$ of the vertices and a series of cyclic permutations $\gamma_i\in Z_{k_i}$ of the incident half-edges at each vertex. It then follows from Lemma \ref{lem_graph_chord} that $c=(\sigma,\gamma)^{-1}\cdot\bar{c}_{\tilde{\Gamma}'}$.

Let $(\sigma,\gamma)\in H$ and let $c:=(\sigma,\gamma)\cdot\bar{c}_{\tilde{\Gamma}'}$. There exist 2-cycles $\tau_1,\ldots,\tau_k$ such that
\[(\sigma,\gamma)\cdot c_{\tilde{\Gamma}'}=\tau_1\cdots\tau_k\cdot\hat{c}.\]
Since we assumed that $\Gamma'$ had no orientation reversing automorphisms, it follows from the definitions that,
\begin{equation} \label{eqn_Feynman_spinvariantsdummy2}
\langle\Gamma_{k_1,\ldots,k_l}(\hat{c}),\Gamma'\rangle = \sgn(\sigma)\sgn(\tau_1)\cdots\sgn(\tau_k).
\end{equation}

Now it follows from Remark \ref{rem_feynman_equivariance} and equation \eqref{eqn_Feynman_spinvariantsdummy2} that for all $(\sigma,\gamma)\in K$,
\[ \sgn(\sigma)\tilde{F}_{\tilde{\Gamma}'}\left[(\sigma,\gamma)\cdot(x_1\otimes\ldots\otimes x_l)\right] = \tilde{F}_{\tilde{\Gamma}'}\left[x_1\otimes\ldots\otimes x_l\right].\]

Finally, we use all the preceding observations together with the orbit-stabiliser theorem to verify equation \eqref{eqn_Feynman_spinvariants}:
\begin{displaymath}
\begin{split}
\langle\langle x,\Gamma' \rangle\rangle &= \frac{1}{|\Aut(\Gamma')|} \sum_{(\sigma,\gamma)\in H} \sgn(\sigma) \tilde{F}_{\tilde{\Gamma}'}\left[(\sigma,\gamma)\cdot x\right], \\
&= \frac{|K|}{|\Aut(\Gamma')|} \sum_{c\in\orb\left(\bar{c}_{\tilde{\Gamma}'}\right)}\langle\Gamma_{k_1,\ldots,k_l}(\hat{c}),\Gamma'\rangle\tilde{F}_{\grpchd{\hat{c}}}[x], \\
&= \sum_{c\in\chd{k}}\langle\Gamma_{k_1,\ldots,k_l}(\hat{c}),\Gamma'\rangle\beta_{\hat{c}}(x), \\
&= \langle I(x),\Gamma' \rangle.
\end{split}
\end{displaymath}
\end{proof}

Recall that the complex $\stlim$ is the direct limit of the complexes $\cek[\bullet\bullet]$,
\[\stlim=\dilim{n,m}{\cek[\bullet\bullet]}.\]
The bigrading is given as follows: an element $x\in\stlim$ has bidegree $(i,j)$ if it lies in the $i$th exterior power and has order $2j$ (the invariant theory for the Lie algebras $\osp$ implies that there are no elements of odd order, cf. \cite{sergeev} and \cite{hamgraph}, Lemma 4.3).

The maps $I:\cek[\bullet\bullet]\to\gc$ defined in Definition \ref{def_integral} give rise to a map
\begin{equation} \label{eqn_integral}
I:\stlim\to\gc
\end{equation}
on the direct limit. The following theorem is the appropriate super-analogue of Kontsevich's well known theorem on graph homology \cite{kontsympgeom}.

\begin{theorem} \label{thm_kontthm}
The map $I:\stlim\to\gc$ is an isomorphism of complexes.
\end{theorem}

\begin{proof}
It follows from Proposition \ref{prop_adjointness} and theorems \ref{thm_kontsevich} and \ref{thm_Feynman_spinvariants} that $I$ is a morphism of complexes. One can establish that the map $I$ is bijective by constructing an explicit inverse for $I$ using the invariant theory for the Lie algebras $\osp$. This was done in \cite{hamgraph} for the graph complex associated to the commutative operad. The proof for the ribbon graph complex, which is the graph complex associated to the associative operad, proceeds mutatis mutandis.
\end{proof}

\subsection{Hopf algebra structure}

In this section we define the Hopf algebra structures present in both the graph complex and the stable relative Chevalley-Eilenberg complex and show that the map \eqref{eqn_integral} is an isomorphism of Hopf algebras. We begin by describing the Hopf algebra structure on the graph complex.

\begin{defi}
Let $\Gamma$ and $\Gamma'$ be two oriented ribbon graphs, we define their disjoint union $\Gamma\sqcup\Gamma'$ in the obvious way as follows:
\begin{enumerate}
\item The set of half-edges of $\Gamma\sqcup\Gamma'$ is the disjoint union of the half-edges of $\Gamma$ and the half-edges of $\Gamma'$.
\item The set of vertices of $\Gamma\sqcup\Gamma'$ is the union of the vertices of $\Gamma$ and $\Gamma'$. The vertices retain their original cyclic ordering. Likewise the set of edges of $\Gamma\sqcup\Gamma'$ is the union of the edges of $\Gamma$ and $\Gamma'$.
\item The orientation on $\Gamma\sqcup\Gamma'$ is given by retaining the original orientation on the edges of $\Gamma$ and $\Gamma'$ and ordering the vertices of $\Gamma\sqcup\Gamma'$ by concatenating the orderings on the vertices of $\Gamma$ and $\Gamma'$.
\end{enumerate}

The operation of disjoint union described above gives the graph complex $\gc$ the structure of a graded-commutative algebra, where the grading is given by counting the number of vertices. It is easy to check that the differential \eqref{eqn_hdiff} is a derivation with respect to the operation of disjoint union. Furthermore, contracting an edge in a \emph{connected} graph produces another \emph{connected} graph and hence it follows that the subspace generated by \emph{connected} oriented ribbon graphs is a subcomplex of the graph complex.

Every oriented ribbon graph can be uniquely decomposed into the disjoint union of its connected components, hence it follows that the graph complex is a free graded-commutative algebra which is freely generated by the set of all \emph{connected} oriented ribbon graphs (the role of the unit is played by the `empty graph'). As a consequence of this and the fact that the connected oriented ribbon graphs generate a subcomplex of the graph complex, it follows that the graph complex has the canonical structure of a commutative-cocommutative differential graded Hopf algebra in which the subcomplex of primitive elements coincides with the subcomplex generated by all \emph{connected} oriented ribbon graphs.
\end{defi}

Now we will define the Hopf algebra structure on the stable relative Chevalley-Eilenberg complex
\[ \stlim=\dilim{n,m}{\cek} .\]

\begin{defi}
Let $\mathfrak{l}$ be any Lie algebra, then the Chevalley-Eilenberg complex $C_\bullet(\mathfrak{l})$ has the structure of a differential graded cocommutative coalgebra, where the coproduct $\Delta$ is induced by the following map of Lie algebras:
\begin{displaymath}
\begin{array}{ccc}
\mathfrak{l} & \to & \mathfrak{l}\oplus\mathfrak{l}, \\
g & \mapsto & (g,g). \\
\end{array}
\end{displaymath}
Furthermore, for any $g\in\mathfrak{l}$ the adjoint action $\ad g:C_\bullet(\mathfrak{l}) \to C_\bullet(\mathfrak{l})$ is a coderivation of the coproduct $\Delta$.

It follows from the general theory outlined above that the complex $\stlim$ has the structure of a differential graded cocommutative coalgebra.

The complex $\stlim$ may be given the structure of an algebra as follows: there is a canonical embedding of Lie algebras
\begin{equation} \label{eqn_ceprod}
\mathfrak{g}_{2n|m}\oplus\mathfrak{g}_{2n'|m'} \to \mathfrak{g}_{2n+2n'|m+m'}
\end{equation}
in which the left factor is embedded on the left in the standard way and the right factor is embedded on the right by shifting the indices by $2n|m$. This embedding determines a product
\[ \mu:C_{\bullet\bullet}(\mathfrak{g}_{2n|m})\otimes C_{\bullet\bullet}(\mathfrak{g}_{2n'|m'}) \to C_{\bullet\bullet}(\mathfrak{g}_{2n+2n'|m+m'}). \]
Since \eqref{eqn_ceprod} is a map of Lie algebras, it follows that the Chevalley-Eilenberg differential is a derivation with respect to this product and furthermore, that this product descends to the level of coinvariants, or relative homology:
\[ \mu:C_{\bullet\bullet}(\mathfrak{g}_{2n|m},\mathfrak{osp}_{2n|m})\otimes C_{\bullet\bullet}(\mathfrak{g}_{2n'|m'},\mathfrak{osp}_{2n'|m'}) \to C_{\bullet\bullet}(\mathfrak{g}_{2n+2n'|m+m'},\mathfrak{osp}_{2n+2n'|m+m'}). \]

Moreover, this map is \emph{well-defined stably} and is \emph{graded-commutative}; that is to say that it gives rise to a graded-commutative product
\[ \mu:\stlim\otimes\stlim \to \stlim. \]
This map is only well-defined and graded-commutative when one considers \emph{relative} Lie algebra homology. This is because (linear) symplectomorphisms act trivially on $\cek[\bullet\bullet]$ in the stable limit. This, in turn, implies the statements above. A subtle point arises here due to a discrepancy between the invariants of the Lie algebra $\osp$ and its Lie group, cf. \cite{sergeev}. Fortunately, this discrepancy disappears as one stabilises.

One can easily see that the coproduct and product defined above are compatible. All told, they give $\stlim$ the structure of a commutative-cocommutative differential graded Hopf algebra.
\end{defi}

\begin{theorem}
The map \eqref{eqn_integral} $I:\stlim\to\gc$ is an isomorphism of differential graded Hopf algebras.
\end{theorem}

\begin{proof}
It follows from Theorem \ref{thm_kontthm} that all we need to prove is that $I$ is a map of Hopf algebras. First let us prove that $I$ is a map of algebras. Let
\begin{equation} \label{eqn_leftrightinc}
\iota_1:\gf^{2n_1|m_1}\to\gf^{2n_1+2n_2|m_1+m_2} \quad \text{and} \quad \iota_2:\gf^{2n_2|m_2}\to\gf^{2n_1+2n_2|m_1+m_2}
\end{equation}
be the canonical inclusions on the left and on the right respectively. These are the inclusions which give rise to the product on $\stlim$ described above. Given two chord diagrams $c\in\chd{k}$ and $c'\in\chd{k'}$ we define the chord diagram $c\sqcup c'$ to be the chord diagram in $\chd{k+k'}$ given by concatenating $c$ with $c'$.

Let $x:=x_1\otimes\ldots\otimes x_l$ and $y:=y_1\otimes\ldots\otimes y_{l'}$ represent Chevalley-Eilenberg chains in $C_{\bullet\bullet}(\mathfrak{g}_{2n_1|m_1},\mathfrak{osp}_{2n_1|m_1})$ and $C_{\bullet\bullet}(\mathfrak{g}_{2n_2|m_2},\mathfrak{osp}_{2n_2|m_2})$ respectively; where $x_i\in(\gf^{2n_1|m_1})^{\otimes k_i}$, $y_j\in(\gf^{2n_2|m_2})^{\otimes k_j}$ and $k_1+\ldots+k_l=2k$, $k'_1+\ldots+k'_l=2k'$:
\begin{displaymath}
\begin{split}
I(x\cdot y) &= I\left[\iota^*_1(x)\wedge\iota^*_2(y)\right], \\
&= \sum_{c\in\chd{k+k'}} \beta_{\hat{c}}\left[\iota^*_1(x)\otimes\iota^*_2(y)\right]\Gamma_{k_1,\ldots,k_l,k'_1,\ldots,k'_{l'}}(\hat{c}).
\end{split}
\end{displaymath}
Now since $\langle\iota_1(a),\iota_2(b)\rangle=0$ for all $a\in\gf^{2n_1|m_1}$ and $b\in\gf^{2n_2|m_2}$, it follows that this sum factorises as follows:
\begin{displaymath}
\begin{split}
I(x\cdot y) &=
\sum_{(c,c')\in\chd{k}\times\chd{k'}}\beta_{\hat{c}\sqcup\hat{c}'}\left[\iota^*_1(x)\otimes\iota^*_2(y)\right]\Gamma_{k_1,\ldots,k_l,k'_1,\ldots,k'_{l'}} \left(\hat{c}\sqcup\hat{c}'\right), \\
&= \sum_{(c,c')\in\chd{k}\times\chd{k'}}\beta_{\hat{c}}[\iota^*_1(x)]\beta_{\hat{c}'}[\iota^*_2(y)]\Gamma_{k_1,\ldots,k_l}(\hat{c})\sqcup\Gamma_{k'_1,\ldots,k'_{l'}}(\hat{c}'), \\
&= I(x)\cdot I(y).
\end{split}
\end{displaymath}
Hence we have shown that $I$ is a map of algebras as claimed.

Now we must show that $I$ is a map of coalgebras, or equivalently, that $I^{-1}$ is a map of coalgebras. This is equivalent to the statement that $I^{-1}$ maps primitive
elements to primitive elements. To prove this we need to make use of the explicit form for $I^{-1}$ constructed in \cite{hamgraph}, Definition 5.6.

Let $\Gamma$ be a connected oriented ribbon graph and let $\tilde{\Gamma}$ be a fully ordered graph representing it of type $(k_1,\ldots,k_l)$, where $k_1+\ldots+k_l=2k$. Let $u\in V^{\otimes 2k}$ be the vector given by the formula
\[ u:=\sigma_{c_{\tilde{\Gamma}}}^{-1}\cdot(p_1\otimes q_1)\ldots(p_k\otimes q_k);\]
that is to say that if $c_{\tilde{\Gamma}}=(i_1,j_1),\ldots,(i_k,j_k)$; then the $i_1$th factor of $u$ is $p_1$, the $j_1$th factor of $u$ is $q_1$, the $i_2$th factor of $u$ is $p_2$ and so on. The vector $u$ represents $I^{-1}(\Gamma)$, i.e. the image of $u$ under the projection
\[ V^{\otimes 2k}\to (V^{\otimes k_1})_{Z_{k_1}}\otimes\ldots\otimes(V^{\otimes k_l})_{Z_{k_l}}\to \Lambda^l(\tilde{\mathfrak{g}}_{2k|0}) \]
maps to $\Gamma$ under $I$. Let us denote this image by
\[ \tilde{u} = g_1\wedge\ldots\wedge g_l. \]

Now
\[ \Delta(\tilde{u}) = \tilde{u}\otimes 1 + 1\otimes\tilde{u} + O(<l)\otimes O(<l),\]
where $O(<l)\otimes O(<l)$ is a sum of terms of the form
\[ (g_{\tau(1)}\wedge\ldots\wedge g_{\tau(r)}) \otimes (g_{\tau(r+1)}\wedge\ldots\wedge g_{\tau(l)}), \]
for some $\tau\in S_l$ and $1\leq r < l$. Since $\Gamma$ is \emph{connected}, it follows that each of these terms has the form
\[ (g_{\tau(1)}\wedge\ldots\wedge(\cdots p_i\cdots)\wedge\ldots\wedge g_{\tau(r)}) \otimes (g_{\tau(r+1)}\wedge\ldots\wedge(\cdots q_i\cdots)\wedge\ldots\wedge g_{\tau(l)}), \]
\emph{where the coordinates $p_i$ and $q_i$ do not occur in any other factor of the above tensor}. Because all those coordinates which occur elsewhere in the above expression all have the form $p_j$ or $q_j$ for $j\neq i$, it follows that both the left factor and the right factor are zero modulo the action of $\mathfrak{osp}_{2k|0}$. For instance the left factor is the image of
\[ (g_{\tau(1)}\wedge\ldots\wedge(\cdots q_i\cdots)\wedge\ldots\wedge g_{\tau(r)}) \]
under the action of the symplectic vector field $q_i\mapsto p_i$.

It follows from the above argument that $\tilde{u}$ is a primitive element and therefore that $I^{-1}$ maps primitives to primitives, thus our claim that $I$ is an isomorphism of differential graded Hopf algebras is established.
\end{proof}

\section{$\ai$-algebras} \label{sec_Ainfinity}

In this section we review the prerequisite definitions of an $\ai$-algebra and a \emph{symplectic} $\ai$-algebra. The former were introduced by Stasheff in \cite{Stas}. The latter are also known as `cyclic' $\ai$-algebras or $\ai$-algebras with an invariant inner product in the literature. Since we shall only be interested in making use of symplectic $\ai$-algebras in the constructions that follow in the next section, we shall assume that the underlying vector space of our $\ai$-algebra is finite dimensional when stating all definitions.

In order to introduce the notion of a symplectic $\ai$-algebra we will have to make use of \emph{formal} noncommutative geometry (cf. \cite{hamlaz}). This is a kind of completion of the noncommutative geometry described in section \ref{sec_noncomgeom}. Since the $\ai$-algebras that we will be working with are assumed to be finite dimensional, this completion just amounts to the usual $I$-adic completion.

The differences are as follows: for a finite dimensional vector space $V$ the tensor algebra $T(V)$ gets replaced with the \emph{completed tensor algebra}
\[ \widehat{T}(V):=\prod_{i=0}^\infty V^{\otimes i}.\]
A \emph{vector field} on $\widehat{T}(V)$ is a derivation
\[ \xi:\widehat{T}(V)\to\widehat{T}(V) \]
which is \emph{continuous} with respect to the $I$-adic topology. We say that this vector field \emph{vanishes at zero} if for every $v\in V$, $\xi(v)$ is a formal power series with \emph{vanishing constant term}. We call a map
\[\phi:\widehat{T}(V)\to\widehat{T}(W) \]
a \emph{diffeomorphism} if it is an invertible homomorphism of algebras which is \emph{continuous} with respect to the $I$-adic topology.

There exists a formal version of the de Rham complex of $V$ which is denoted by $\cdrham{V}$. This complex is just the $I$-adic completion of the ordinary de Rham complex $\drham{V}$. Given a vector field $\xi$ on $\widehat{T}(V)$ one can define formal analogues
\[ L_\xi,i_\xi:\cdrham{V}\to\cdrham{V} \]
of the Lie derivative and contraction operators described in Definition \ref{def_operators}. Likewise, any diffeomorphism
\[ \phi:\widehat{T}(V)\to\widehat{T}(W) \]
gives rise to a map of complexes $\phi^*:\cdrham{V}\to\cdrham{W}$.

Similarly, there is a formal version of noncommutative \emph{symplectic} geometry. A symplectic form is a closed 2-form $\omega\in\cdrham[2]{V}$ such that the formal analogue of the map \eqref{eqn_nondegenerate} is bijective. Likewise, a vector field is said to be symplectic if its Lie derivative annihilates the symplectic form and a diffeomorphism is called a symplectomorphism if it preserves the symplectic forms, cf. Definition \ref{def_sympmaps}.

\begin{defi}
Let $U$ be a vector space:
\begin{enumerate}
\item
An $\ai$-structure on $U$ is a vector field
\[ m:\ctalg{U} \to \ctalg{U} \]
of degree one and vanishing at zero, such that $m^2=0$.
\item
A \emph{symplectic} $\ai$-structure on $U$ is a symplectic form $\omega \in \cdrham[2]{\Pi U^*}$ together with a \emph{symplectic} vector field
\[ m:\ctalg{U} \to \ctalg{U} \]
of degree one and vanishing at zero, such that $m^2=0$.

To any symplectic $\ai$-structure we can associate a 0-form
\[ h\in\prod_{i=2}^\infty \left((\Pi U^*)^{\otimes i}\right)_{Z_i}\subset\cdrham[0]{\Pi U^*} \]
with vanishing constant and linear terms, by the formula $dh=i_m(\omega)$. We call $h$ the Hamiltonian of the symplectic $\ai$-structure.
\end{enumerate}
\end{defi}

\begin{rem}
Let $m:\ctalg{U} \to \ctalg{U}$ be an $\ai$-structure. It may be written as an infinite sum
\[ m=m_1+m_2+\ldots+m_n+\ldots, \]
where $m_k$ is a vector field of order $k$ and corresponds to a map
\[ m_k:(\Pi U)^{\otimes k}\to\Pi U \]
under the isomorphism $\Der(\ctalg{U})\cong\Hom(T(\Pi U),\Pi U)$. We say that $m$ is a minimal $\ai$-structure if $m_1=0$.

Suppose that $(m,\omega)$ is a symplectic $\ai$-structure whose symplectic form $\omega$ \emph{has order zero} and let
\[ \langle -,- \rangle:\Pi U\otimes \Pi U \to \gf \]
be the inner product associated to $\omega$ by the map $\Upsilon$ of Proposition \ref{prop_closed_twoforms}. It follows from Theorem 10.6 of \cite{hamlaz} that the map
\begin{equation} \label{eqn_hamiltonian}
h_k:(\Pi U)^{\otimes k} \to \gf
\end{equation}
defined by the formula $h_k(x_1,\ldots, x_k):=\langle m_{k-1}(x_1,\ldots,x_{k-1}), x_k \rangle$, is cyclically invariant.
\end{rem}

\begin{defi}
Let $U$ and $W$ be vector spaces:
\begin{enumerate}
\item
Let $m$ and $m'$ be $\ai$-structures on $U$ and $W$ respectively. An $\ai$-morphism from $U$ to $W$ is a continuous algebra homomorphism
\[ \phi:\ctalg{W} \to \ctalg{U} \]
of degree zero such that $\phi \circ m'=m \circ \phi$.
\item
Let $(m,\omega)$ and $(m',\omega')$ be \emph{symplectic} $\ai$-structures on $U$ and $W$ respectively. A \emph{symplectic} $\ai$-morphism from $U$ to $W$ is a continuous algebra homomorphism
\[ \phi:\ctalg{W} \to \ctalg{U} \]
of degree zero such that $\phi \circ m' = m \circ \phi$ and $\phi^*(\omega') = \omega$.
\end{enumerate}
\end{defi}

We will now formulate the Darboux theorem which says that there is a canonical form for any symplectic form. We refer the reader to \cite[\S6]{ginzsympgeom} for a proof:

\begin{theorem} \label{thm_drboux}
Let $V$ be a vector space of dimension $2n|m$ and let $\omega \in \cdrham[2]{V}$ be a homogeneous (either purely even or purely odd) 2-form. The 2-form $\omega$ can be written as an infinite sum
\[ \omega = \omega_0 + \omega_1 + \ldots + \omega_r + \ldots, \]
where $\omega_i$ is a 2-form of order $i$:
\begin{enumerate}
\item
The 2-form $\omega$ is nondegenerate if and only if the 2-form $\omega_0$ is nondegenerate.
\item \label{item_drbouxnonlin}
Suppose that $\omega$ is in fact a symplectic form, then there exists a symplectic vector field $\gamma\in\Der(\widehat{T}(V))$ consisting of terms of order $\geq 2$ and having degree zero, such that $\exp(\gamma)^*[\omega] = \omega_0$.
\item \label{item_drbouxlin}
Furthermore, if $\omega$ is an \emph{even} symplectic form, there exists a linear symplectomorphism $\phi:V\to \gf^{2n|m}$ such that
\[ \phi^*(\omega_0)=\sum_{i=1}^n dp_i dq_i + \frac{1}{2}\sum_{i=1}^m dx_i dx_i.\]
\end{enumerate}
\end{theorem}
\noproof

The following lemma describes how the Hamiltonians of two isomorphic symplectic $\ai$-algebras are related.

\begin{lemma} \label{lem_hamsymplecto}
Let $V$ and $W$ be vector spaces and let $\omega\in\cdrham[2]{V}$ and $\omega'\in\cdrham[2]{W}$ be symplectic forms. For any symplectomorphism $\psi:\widehat{T}(V)\to\widehat{T}(W)$ we have the following commutative diagram:
\begin{displaymath}
\xymatrix{ \cdrham[0]{W} \ar[r]^{d} & \cdrham[1]{W} && \Der(\widehat{T}(W)) \ar[ll]_{i_\xi(\omega')\mapsfrom\xi} \\ \cdrham[0]{V} \ar[r]^{d} \ar[u]^{\psi^*} & \cdrham[1]{V} \ar[u]^{\psi^*} && \Der(\widehat{T}(V)) \ar[ll]_{i_\xi(\omega)\mapsfrom\xi} \ar[u]_{\scriptsize{\begin{array}{c} \psi\xi\psi^{-1} \\ \uparrow \ \ \\ \xi \ \ \end{array}}} }
\end{displaymath}
\end{lemma}

\begin{proof}
We need only prove the commutativity of the right-hand side. This follows from a one line calculation:
\[ i_{\psi\xi\psi^{-1}}(\omega') = \psi^*i_\xi\psi^{*-1}(\omega') = \psi^*i_\xi(\omega).\]
\end{proof}

\begin{example}
Given any two $\ai$-algebras one can define their direct sum in a natural way. Let
\[ A_1:=(U_1,h_1,\omega_1) \quad\text{and}\quad A_2:=(U_2,h_2,\omega_2) \]
be two symplectic $\ai$-algebras, where $h_1\in\cdrham[0]{\Pi U_1^*}$ and $h_2\in\cdrham[0]{\Pi U_2^*}$ are the Hamiltonians of the symplectic $\ai$-structures. Let
\[ \iota_1:\Pi U_1^*\to \Pi U_1^* \oplus \Pi U_2^* \quad\text{and}\quad \iota_2:\Pi U_2^*\to \Pi U_1^* \oplus \Pi U_2^* \]
be the canonical inclusions. The direct sum $A:=(U,h,\omega)$ of $A_1$ and $A_2$ is defined by the formulae:
\begin{enumerate}
\item $U:=U_1\oplus U_2$,
\item $h:=\iota^*_1(h_1)+\iota^*_2(h_2)$,
\item $\omega:=\iota^*_1(\omega_1)+\iota^*_2(\omega_2)$.
\end{enumerate}
\end{example}

\section{Constructing graph homology classes} \label{sec_classes}

In this section we describe a combinatorial construction, due to Kontsevich \cite{kontfeynman}, which produces a class in graph homology from the initial data of a (minimal) symplectic $\ai$-algebra. This construction is known as a \emph{partition function}. An important point here is that the symplectic form for this symplectic $\ai$-algebra should be \emph{even}; we hereafter assume this property of all the symplectic $\ai$-algebras in this section. We also describe an algebraic construction which produces classes in the relative Lie algebra homology of $\mathfrak{g}$ from the same initial data. These classes are called \emph{characteristic classes} and we prove that they are a homotopy invariant of the symplectic $\ai$-algebra. Finally, we show that these two constructions are two sides of the same story, in that they produce the same class when one identifies Lie algebra homology with graph homology using Kontsevich's theorem.

\subsection{Partition functions of $\ai$-algebras}

In this section we define the \emph{partition function} of a symplectic $\ai$-algebra. It will follow from Theorem \ref{thm_main} that this produces a class in graph homology.

Let $A:=(U,m,\omega)$ be a symplectic $\ai$-algebra whose symplectic form has order zero and let $\innprod:=\Upsilon(\omega)$ be the inner product associated to $\omega$ by Lemma \ref{lem_symp_innprod}. Since this
inner product provides an identification of $\Pi U$ with $\Pi U^*$ it gives rise to an inner product $\innprod^{-1}$ on $\Pi U^*$ which is defined by the following commutative
diagram:
\begin{equation} \label{eqn_dualinn}
\xymatrix{ \Pi U \ar@{=}[rr] \ar[rd]_{x\mapsto \langle x,-\rangle} && \Pi U^{**} \\ & \Pi U^* \ar[ru]_{y\mapsto \langle y,-\rangle^{-1}} }
\end{equation}
This gives us a map $\kappa^A:(\Pi U^*)^{\otimes 2k}\to\gf$ defined by the formula
\[ \kappa^A(x_1\otimes\ldots\otimes x_{2k}):= \langle x_1,x_2 \rangle^{-1}\langle x_3,x_4 \rangle^{-1}\ldots\langle x_{2k-1},x_{2k}\rangle^{-1}. \]

Given the initial data of a (minimal) symplectic $\ai$-algebra one can construct a function on graphs which we will call the partition function of this symplectic $\ai$-algebra. It will
follow from subsequent results that this chain will be a cycle in the graph complex.

\begin{defi}\label{partit}
Let $A:=(U,m,\omega)$ be a \emph{minimal} symplectic $\ai$-algebra whose symplectic form is \emph{even} and \emph{has order zero}. Let $\innprod^{-1}$ be the inner product on $\Pi U^*$ given by diagram \eqref{eqn_dualinn} above and let $h_k\in (\Pi U^*)^{\otimes k}$ be the series of tensors defined by \eqref{eqn_hamiltonian}. We define a function $Z_A:\gc\to\gf$ on graphs, called the partition function of $A$, as follows:

Let $\Gamma$ be an oriented ribbon graph and choose some fully ordered graph $\tilde{\Gamma}$ of type $(k_1,\ldots,k_m)$ which represents it. Define the tensor $x_{\tilde{\Gamma}}\in (\Pi U^*)^{\otimes k_1}\otimes\ldots\otimes (\Pi U^*)^{\otimes k_m}$ by the formula
\[ x_{\tilde{\Gamma}}:=h_{k_1}\otimes\ldots\otimes h_{k_m}.\]
The partition function $Z_A$ is given by the formula,
\begin{equation} \label{eqn_partfun}
Z_A(\Gamma):=\frac{\tilde{F}_{\tilde{\Gamma}}^A(x_{\tilde{\Gamma}})}{|\Aut(\Gamma)|}:=\frac{\kappa^A(\sigma_{c_{\tilde{\Gamma}}}\cdot x_{\tilde{\Gamma}})}{|\Aut(\Gamma)|}.
\end{equation}
\end{defi}

\begin{rem}
Since the tensors $h_k\in (\Pi U^*)^{\otimes k}$ are \emph{cyclically invariant} tensors of \emph{degree one}, it follows from Remark \ref{rem_feynman_equivariance} that the partition function $Z_A$ is well defined; that is to say that $Z_A(\Gamma)$ is independent of the choice of fully ordered graph representing $\Gamma$.
\end{rem}

\begin{rem} \label{rem_onlyeven}
Note that since the tensors $h_k\in (\Pi U^*)^{\otimes k}$ are odd and the map $\kappa$ is even, it follows that $Z_A$ only produces classes in bidegree $(2\bullet,\bullet)$; that is to say that $Z_A(\Gamma)$ vanishes on any graph with an odd number of vertices.
\end{rem}

This definition is illustrated graphically as follows:
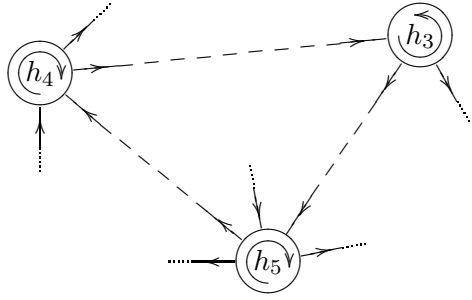
\begin{displaymath}
\begin{xy}
(0,0)="a" ;
"a"*[o]=<9mm>{h_5} *\frm{o} ; "a"*\cir<3mm>{l_d} ; "a"+(2.8,-0.7)*\dir{<} ;
(-30,25)="b" ;
"b"*[o]=<9mm>{h_4} *\frm{o} ; "b"*\cir<3mm>{l_d} ; "b"+(2.8,-0.7)*\dir{<} ;
(20,30)="c" ;
"c"*[o]=<9mm>{h_3} *\frm{o} ; "c"*\cir<3mm>{d^l} ; "c"+(-0.7,2.8)*_\dir{<} ;
"a"+a(10);"a" ; **{} ; ?(0)/4.5mm/="t1" ; "a"+a(10);"a" ; **{} ; ?(0)/10mm/="t2";"t1" ; **\dir{-} ; ?*!/-1mm/\dir{>} ; "t2"+a(10);"t2" ; **{} ; ?(0)/3mm/;"t2" ; **\dir{.} ;
"a"+a(100);"a" ; **{} ; ?(0)/4.5mm/="t1" ; "a"+a(100);"a" ; **{} ; ?(0)/10mm/="t2";"t1" ; **\dir{-} ; ?*!/1mm/\dir{<} ; "t2"+a(100);"t2" ; **{} ; ?(0)/3mm/;"t2" ; **\dir{.} ;
"a"+a(180);"a" ; **{} ; ?(0)/4.5mm/="t1" ; "a"+a(180);"a" ; **{} ; ?(0)/10mm/="t2";"t1" ; **\dir{-} ; ?*!/-1mm/\dir{>} ; "t2"+a(180);"t2" ; **{} ; ?(0)/3mm/;"t2" ; **\dir{.} ;
"b"+a(45);"b" ; **{} ; ?(0)/4.5mm/="t1" ; "b"+a(45);"b" ; **{} ; ?(0)/10mm/="t2";"t1" ; **\dir{-} ; ?*!/-1mm/\dir{>} ; "t2"+a(45);"t2" ; **{} ; ?(0)/3mm/;"t2" ; **\dir{.} ;
"b"+a(270);"b" ; **{} ; ?(0)/4.5mm/="t1" ; "b"+a(270);"b" ; **{} ; ?(0)/10mm/="t2";"t1" ; **\dir{-} ; ?*!/1mm/\dir{<} ; "t2"+a(270);"t2" ; **{} ; ?(0)/3mm/;"t2" ; **\dir{.} ;
"c"+a(300);"c" ; **{} ; ?(0)/4.5mm/="t1" ; "c"+a(300);"c" ; **{} ; ?(0)/10mm/="t2";"t1" ; **\dir{-} ; ?*!/-1mm/\dir{>} ; "t2"+a(300);"t2" ; **{} ; ?(0)/3mm/;"t2" ; **\dir{.} ;
"a";"b" ; **{} ; ?(0)/4.5mm/="t_1" ; ?(0)/10.5mm/="t_2" ; ?(1)/-10.5mm/="t_3" ; ?(1)/-4.5mm/="t_4" ;
"t_1";"t_2" ; **\dir{-} ; ?*\dir{>} ; "t_2";"t_3" ; **\dir{--} ; "t_3";"t_4" ; **\dir{-} ; ?*!/-1mm/\dir{>} ;
"b";"c" ; **{} ; ?(0)/4.5mm/="t_1" ; ?(0)/10.5mm/="t_2" ; ?(1)/-10.5mm/="t_3" ; ?(1)/-4.5mm/="t_4" ;
"t_1";"t_2" ; **\dir{-} ; ?*\dir{>} ; "t_2";"t_3" ; **\dir{--} ; "t_3";"t_4" ; **\dir{-} ; ?*!/-1mm/\dir{>} ;
"c";"a" ; **{} ; ?(0)/4.5mm/="t_1" ; ?(0)/10.5mm/="t_2" ; ?(1)/-10.5mm/="t_3" ; ?(1)/-4.5mm/="t_4" ;
"t_1";"t_2" ; **\dir{-} ; ?*\dir{>} ; "t_2";"t_3" ; **\dir{--} ; "t_3";"t_4" ; **\dir{-} ; ?*!/-1mm/\dir{>} ;
\end{xy}
\end{displaymath}
Given an oriented ribbon graph $\Gamma$ one attaches the tensors $h_k$ to each vertex of $\Gamma$ of valency $k$ using the cyclic ordering of the incident half-edges; this is
possible since the tensors $h_k$ are cyclically invariant. The edges of the graph $\Gamma$ give us a way to contract these tensors using the bilinear form $\langle -,- \rangle^{-1}$ and produce a number $Z_A(\Gamma)\in\gf$, just as we did in Example \ref{exm_feynman_amplitude}.

\begin{defi}
The subspace $\gc^{\mathcal{C}}$ of the graph complex spanned by \emph{connected} graphs is a subcomplex of $\gc$. Given any minimal symplectic $\ai$-algebra $A:=(U,m,\omega)$ whose symplectic form is even and has order zero, we define the \emph{connected level} partition function
\[ Z_A^{\mathcal{C}}:\gc^{\mathcal{C}}\to\gf \]
to be the restriction of the partition function $Z_A$ to $\gc^{\mathcal{C}}$.

Using the pairing $\innprod$ given by Definition \ref{def_graphinn} we may regard the partition function as a member of $\gc$. In order to avoid confusion, we denote this chain by $\tilde{Z}_A$; that is to say that $\tilde{Z}_A$ is the unique member of $\gc$ satisfying the identity,
\[ Z_A(\Gamma)=\langle \tilde{Z}_A,\Gamma \rangle, \quad \text{for all } \Gamma\in\gc.\]
Equivalently, one may write $\tilde{Z}_A:=\sum_\Gamma Z_A(\Gamma)\Gamma$; where the sum is taken over all isomorphism classes of ribbon graphs.

Since the restriction of the inner product $\innprod$ to the subcomplex $\gc^{\mathcal{C}}$ spanned by connected graphs is also nondegenerate, one can similarly define $\tilde{Z}_A^{\mathcal{C}}$ to be the unique member of $\gc^{\mathcal{C}}$ such that
\[ Z_A^{\mathcal{C}}(\Gamma)=\langle \tilde{Z}_A^{\mathcal{C}},\Gamma \rangle, \quad \text{for all } \Gamma\in\gc^{\mathcal{C}}.\]
Equivalently, one may write $\tilde{Z}_A^{\mathcal{C}}:=\sum_\Gamma Z_A(\Gamma)\Gamma$; where the sum is taken over all isomorphism classes of \emph{connected} ribbon graphs.
\end{defi}

A standard result, see e.g. \cite{Manin}, is that the partition function $Z_A$ can be calculated in terms of the connected level partition function $Z_A^{\mathcal{C}}$.

\begin{prop} \label{prop_exppar}
Let $A:=(U,m,\omega)$ be a minimal symplectic $\ai$-algebra whose symplectic form is even and has order zero.
\[ \tilde{Z}_A = \exp\left(\tilde{Z}_A^{\mathcal{C}}\right). \]
\end{prop}

\begin{proof}
Let $\Gamma_i, i\in \mathcal{I}$ be a basis of $\gc^{\mathcal{C}}$ consisting of connected oriented ribbon graphs. Let $\Gamma$ be an oriented ribbon graph, we may assume that  there exists an $\mathbf{i}:=(i_1,\ldots,i_m)\in\mathcal{I}^m$ such that
\[ \Gamma=\Gamma_{i_1}\sqcup\ldots\sqcup\Gamma_{i_m}; \]
where each component $\Gamma_{i_r}$ has an even number of vertices, since otherwise by Remark \ref{rem_onlyeven} the contributions to the partition function would be zero.

$S_m$ acts on the set $\mathcal{I}^m$ by permuting the indices. Let $K\subset S_m$ be the stabiliser of the point $\mathbf{i}$. Now clearly any automorphism of $\Gamma_{i_r}$ determines an automorphism of $\Gamma$. Furthermore any $\sigma\in K$ also determines an automorphism of $\Gamma$ by simply permuting the corresponding connected components of $\Gamma$. Conversely, since the decomposition of $\Gamma$ into connected components is unique and since the image of any connected component of $\Gamma$ under an automorphism is a connected component of $\Gamma$, it follows that the automorphism group of $\Gamma$ is given by the wreath product:
\begin{equation} \label{eqn_exppardummy1}
\Aut(\Gamma) \cong K\ltimes\left(\Aut(\Gamma_{i_1})\times\ldots\times\Aut(\Gamma_{i_m})\right).
\end{equation}
Furthermore, it follows that the orbit of the point $\mathbf{i}$ is
\begin{equation} \label{eqn_exppardummy2}
\orb(\mathbf{i}) = \{(r_1,\ldots,r_m)\in\mathcal{I}^m : \Gamma_{r_1}\sqcup\ldots\sqcup\Gamma_{r_m} \cong \Gamma\}.
\end{equation}

Now the proposition follows from the following calculation:
\begin{displaymath}
\begin{split}
\left\langle \exp\left(\tilde{Z}_A^{\mathcal{C}}\right),\Gamma \right\rangle &= \sum_{n=0}^\infty \frac{1}{n!}\left[\sum_{(r_1,\ldots,r_n)\in\mathcal{I}^n} Z_A(\Gamma_{r_1})\cdots Z_A(\Gamma_{r_n}) \langle\Gamma_{r_1}\sqcup\ldots\sqcup\Gamma_{r_n},\Gamma\rangle\right], \\
&= \frac{1}{m!}\left[\sum_{(r_1,\ldots,r_m)\in\orb(\mathbf{i})} Z_A(\Gamma_{r_1})\cdots Z_A(\Gamma_{r_m})\right], \\
&= \frac{1}{|K|} Z_A(\Gamma_{i_1})\cdots Z_A(\Gamma_{i_m}), \\
&=  Z_A(\Gamma).
\end{split}
\end{displaymath}
Line 2 follows from \eqref{eqn_exppardummy2} and the fact that each $\Gamma_{i_r}$ has an even number of vertices. Line 3 follows from the orbit-stabiliser theorem. Finally, line 4 follows from \eqref{eqn_exppardummy1}.
\end{proof}

\subsection{Characteristic classes of $\ai$-algebras}

In this section we describe how, given any symplectic $\ai$-algebra, one can construct a class in the relative Lie algebra homology of $\mathfrak{g}$ which we will call the characteristic class of this symplectic $\ai$-algebra.

In fact we shall restrict our treatment to \emph{minimal} symplectic $\ai$-algebras. It is possible to include the non-minimal case into our scheme, however that will unnecessarily complicate the exposition while not leading to any new phenomena. The reason for this is that the characteristic class of a non-minimal symplectic $A_\infty$-algebra naturally takes values in a complex isomorphic to the graph complex with bivalent vertices, while the latter could be shown to contain essentially no more information than the usual graph complex. Furthermore, any symplectic $A_\infty$-algebra is weakly equivalent to a minimal one and (as will be shown later) different minimal models have the same characteristic class.

\begin{defi}
Let $A:=(U,m,\omega)$ be a minimal symplectic $\ai$-algebra of dimension $2n|m$, whose symplectic form is \emph{even} and \emph{has order zero}; let $h\in \cdrham[0]{\Pi U^*}$ be its Hamiltonian. By Theorem \ref{thm_drboux} \eqref{item_drbouxlin} there is a linear symplectomorphism
\[ \phi:\Pi U^*\to\gf^{2n|m} \]
transforming $\omega$ to the canonical symplectic form \eqref{eqn_canonicalsymplecticform}. Let $h':=\phi^*(h)\in\tglie$ be the image of $h$ under this symplectomorphism. We define a chain $c_A\in\stlim$, called the \emph{characteristic class} of $A$, by the following formula:
\begin{equation} \label{eqn_charclass}
c_A:=\exp(h'):=1+h'+\frac{1}{2}h'\wedge h'+\ldots+\frac{1}{i!}\underbrace{h'\wedge\ldots\wedge h'}_{\text{$i$ factors}}+\ldots.
\end{equation}
\end{defi}

\begin{rem}
It follows from Lemma 4.3 (ii) of \cite{hamgraph} that the characteristic class $c_A$ is well defined; that is to say that since, as has already been noted, (linear) symplectomorphisms act trivially in the stable limit $\stlim$, it follows that $c_A$ is independent of the choice of linear symplectomorphism $\phi$.

Note also that the chain $c_A$ is inhomogeneous, it is an infinite sum of chains in each bidegree of $\stlim$ and as such lives in the completion of $\stlim$ given by
\[ \widehat{C}_{\bullet\bullet}(\mathfrak{g},\mathfrak{osp}):=\prod_{i,j=0}^\infty \stlim[ij]. \]
Since this obviously causes no problems, we continue to regard it as a member of $\stlim$ by an abuse of notation.
\end{rem}

Let us denote the category of minimal symplectic $\ai$-algebras with even symplectic forms by $\SA$ and let $\ISA$ be the class whose objects are isomorphism classes of these minimal symplectic $\ai$-algebras. We will denote the full subcategory of $\SA$ whose objects consist of minimal symplectic $\ai$-algebras whose symplectic form \emph{has order zero} by $\SAZ$. The following theorem tells us that the characteristic class is a homotopy invariant of the symplectic $\ai$-algebra.

\begin{theorem} \label{thm_invariant}
The characteristic class of any minimal symplectic $\ai$-algebra is a Chevalley-Eilenberg cycle. Furthermore, isomorphic minimal symplectic $\ai$-algebras give rise to homologous cycles. As such the characteristic class map \eqref{eqn_charclass} on $\SAZ$ can be uniquely lifted to a map on $\ISA$ such that the following diagram commutes:
\begin{displaymath}
\xymatrix{\ISA \ar[rr]^c && \hstlim \\ \SAZ \ar[u] \ar[rru]^c}
\end{displaymath}
\end{theorem}

\begin{proof}
Firstly, let us show that for any object $A:=(U,m,\omega)$ of $\SAZ$, the chain $c_A$ is a cycle. Note that the $\ai$-constraint is equivalent to the condition $[m,m]=0$. By  Proposition \ref{prop_hamsympiso} and Lemma \ref{lem_hamsymplecto} it follows that $\{h',h'\}=0$. This immediately implies that the chain \eqref{eqn_charclass} is a Chevalley-Eilenberg cycle.

It follows from the Darboux theorem (Theorem \ref{thm_drboux} \eqref{item_drbouxnonlin}) that any such lifting of \eqref{eqn_charclass} must be unique, since any object in $\SA$ is isomorphic to an object of the subcategory $\SAZ$. Let
\[ A_1:=(U_1,m_1,\omega_1) \quad \text{and} \quad A_2:=(U_2,m_2,\omega_2) \]
be isomorphic objects of $\SAZ$. In order to prove that a lifting of \eqref{eqn_charclass} exists, we must show that $A_1$ and $A_2$ determine the same characteristic class in $\hstlim$.

Let $h_1',h_2'\in\tglie$ be the Hamiltonians corresponding to $A_1$ and $A_2$ respectively. By Lemma \ref{lem_hamsymplecto} there exists a symplectomorphism
\[ \phi:\widehat{T}(\gf^{2n|m})\to\widehat{T}(\gf^{2n|m}) \]
such that $h_2'=\phi^*(h_1')$. Any such symplectomorphism factors as $\phi=\exp(\gamma)\phi_1$, where $\phi_1$ is a linear symplectomorphism and $\gamma$ is an even symplectic vector field consisting of terms of order $\geq 2$. Since any Lie algebra acts trivially on its own homology and linear symplectomorphisms act trivially in the stable limit; it follows that $A_1$ and $A_2$ give rise to homologous characteristic classes. The theorem now follows.
\end{proof}

\begin{rem}
One should, of course, not expect the characteristic class to distinguish inequivalent $\ai$-algebras of the same dimension. First of all, the characteristic class is clearly continuous and so any two $\ai$-algebras, one of which lies in the closure of the other's orbit under the action of the group of formal symplectomorphisms, will have the same characteristic class. Furthermore, $\ai$-algebras with vanishing higher products -- noncommutative Frobenius algebras -- determine classes in $H^0$ of the moduli spaces
of curves. Clearly the only two invariants present -- the genus and the number of marked points -- cannot classify the vast variety of Frobenius algebras.
\end{rem}

\begin{prop}
Let $A_1:=(U_1,h_1,\omega_1)$ and $A_2:=(U_2,h_2,\omega_2)$ be two symplectic $\ai$-algebras.
\[ c_{A_1\oplus A_2} = c_{A_1}\cdot c_{A_2}. \]
\end{prop}

\begin{proof}
We may assume that the vector space $\Pi U_1^*=\gf^{2n_1|m_1}$, that $\Pi U_2^*=\gf^{2n_2|m_2}$ and that the symplectic forms $\omega_1$ and $\omega_2$ are already written in the canonical form \eqref{eqn_canonicalsymplecticform}; then it follows from the definition of the direct sum of two symplectic $\ai$-algebras that
\[ c_{A_1\oplus A_2} = \exp[\iota^*_1(h_1)+\iota^*_2(h_2)], \]
where $\iota_1$ and $\iota_2$ are the inclusions \eqref{eqn_leftrightinc} given by the definition of the Hopf structure on $\stlim$. The proposition now follows from the standard multiplicative properties of the exponential function.
\end{proof}

\subsection{Equivalence of the two constructions}

In this section we show that the partition function and the characteristic class of a symplectic $\ai$-algebra give rise to the same chain in the graph complex via the pairing \eqref{eqn_grphpair} which was used to establish an isomorphism between relative Lie algebra and graph homology.

\begin{theorem} \label{thm_main}
The following diagram commutes:
\begin{displaymath}
\xymatrix{ & \Hom(\gc,\gf) \\ \SAZ \ar[ru]^Z \ar[rd]_c \\ & \stlim \ar[uu]_{\begin{array}{c} \langle\langle x,- \rangle\rangle \\ \uparrow \\ x \end{array}}}
\end{displaymath}
\end{theorem}

\begin{proof}
Let $A:=(U,m,\omega)$ be an object of $\SAZ$ and let $h\in\cdrham[0]{\Pi U^*}$ be its Hamiltonian. Let $V:=\Pi U^*$, we may assume that $V=\gf^{2n|m}$ and that $\omega$ is written in the canonical form \eqref{eqn_canonicalsymplecticform}. A direct calculation shows that the inner product on $V$ determined by diagram \eqref{eqn_dualinn} is just the canonical bilinear form $\innprod$ given by \eqref{eqn_canonical_form}.

Let $\Theta:\prod_{k=1}^\infty V^{\otimes k}\to \cdrham[1]{V}$ be the map defined by Lemma \ref{lem_oneformiso} and let $h_k\in V^{\otimes k}$ be the tensors defined by \eqref{eqn_hamiltonian}. Lemma \ref{lem_normdiff} gives us the following following identity:
\begin{equation} \label{eqn_maindummy2}
\begin{split}
N\cdot h &= \Theta^{-1}d(h) \\
&=\Theta^{-1}[i_m(\omega)] = \sum_{k=3}^\infty h_k,
\end{split}
\end{equation}
where the last line follows by making an explicit calculation.

Now we are ready to make the calculation that will prove the theorem. Let $\Gamma$ be an oriented ribbon graph and let $\tilde{\Gamma}$ be some fully ordered graph of type $(k_1,\ldots,k_r)$ representing it.
\begin{displaymath}
\begin{split}
F_\Gamma(c_A) &= \tilde{F}_{\tilde{\Gamma}}\left(\epsilon\circ T(N)\left[\sum_{i=0}^\infty \frac{1}{i!}h^{\otimes i}\right]\right), \\
&= \tilde{F}_{\tilde{\Gamma}}\left(\epsilon\left[\sum_{i=0}^\infty \frac{1}{i!}\left(\sum_{k=3}^\infty h_k\right)^{\otimes i}\right]\right), \\
&= \tilde{F}_{\tilde{\Gamma}}\left(\sum_{i=0}^\infty \left(\sum_{k=3}^\infty h_k\right)^{\otimes i}\right), \\
&= \tilde{F}_{\tilde{\Gamma}}(h_{k_1}\otimes\ldots\otimes h_{k_r}) = |\Aut(\Gamma)|Z_A(\Gamma).
\end{split}
\end{displaymath}
Line 2 follows from equation \eqref{eqn_maindummy2}. Line 3 follows by virtue of the fact that the Hamiltonian $h$ is odd. Finally, line 4 follows directly from the definition of $\tilde{F}_{\tilde{\Gamma}}$.
\end{proof}

\begin{cor}
The partition function of a minimal symplectic $\ai$-algebra is a cycle in the graph complex. Furthermore, isomorphic minimal symplectic $\ai$-algebras give rise to homologous partition functions.
\end{cor}

\begin{proof}
This is a direct consequence of theorems \ref{thm_invariant} and \ref{thm_main}.
\end{proof}

\section{Topological conformal field theories} \label{sec_TCFT}

In this section we briefly outline one application of our methods to topological conformal field theory (TCFT). We begin by introducing a certain differential graded category $\FT$. Denote by $V$ the infinite dimensional `symplectic' vector space $\dilim{n,m}{\gf^{2n|m}}$.

\begin{defi}
The objects of the category $\FT$ are the symbols $[m]$ where $m$ is a nonnegative integer. The space of morphisms $[m]\rightarrow [n]$ is the relative \emph{cohomological} Chevalley-Eilenberg complex $C^\bullet(\mathfrak{g},\mathfrak{osp};\Hom(V^{\otimes m},V^{\otimes n}))$. Here $\mathfrak{g}$ acts on $\Hom(V^{\otimes m},V^{\otimes n})$ via the projection $\mathfrak{g}\rightarrow \mathfrak{osp}$ and the composition
\begin{equation} \label{comp}
C^\bullet(\mathfrak{g},\mathfrak{osp};\Hom(V^{\otimes m},V^{\otimes
n}))\otimes C^\bullet(\mathfrak{g},\mathfrak{osp};\Hom(V^{\otimes n},V^
{\otimes k}))\rightarrow C^\bullet(\mathfrak{g},\mathfrak{osp};\Hom(V^
{\otimes m},V^{\otimes k}))
\end{equation}
is defined through the obvious map $\Hom(V^{\otimes m},V^{\otimes n})\otimes \Hom(V^{\otimes n},V^{\otimes k})\rightarrow \Hom(V^{\otimes m},V^{\otimes k})$.
\end{defi}

The category $\FT$ has an obvious monoidal structure: $[m]\otimes [n]=[m+n]$. We will call an \emph{open topological conformal field theory} (TCFT) a dg monoidal functor from $\FT$ into the category of $\gf$-vector spaces with monoidal structure given by the usual tensor product of vector spaces and trivial dg structure. Below we will omit the adjective `open' since no other TCFT's will be considered.

To connect this definition of a TCFT with the standard one we introduce a generalisation of the complex of ribbon graphs which allows \emph{graphs with legs}. In other words our ribbon graphs are now allowed to have 1-valent vertices (but not bi-valent vertices), all of which are labelled. These vertices are called \emph{external} and other vertices are called \emph{internal}. The edges attached to external vertices are called \emph{external} edges or \emph{legs} and other edges are called \emph{internal} edges. All legs are partitioned into two classes -- incoming and outgoing legs. The orientation of a graph with legs is the ordering (up to an even permutation) of the set of internal edges and of the sets of endpoints of every edge that is not a leg.

\begin{defi}
The underlying space of the complex $\gc(m,n)$ is the free $\gf$-vector space spanned by isomorphism classes of oriented ribbon graphs with $m$ incoming and $n$ outgoing legs. The bigrading is introduced according to the number of internal vertices and internal edges. The differential $\gc(m,n) \rightarrow \gc[\bullet+1,\bullet+1](m,n)$ is introduced just as in the legless case by expanding internal vertices. The homological ribbon graph complex is introduced similarly.
\end{defi}

\begin{rem}
Suppose that we have a pair of ribbon graphs $\Gamma_1,\Gamma_2$ where the number of incoming legs of $\Gamma_2$ is the same as the number of outgoing legs of $\Gamma_1$. Gluing the corresponding legs of $\Gamma_1$ and $\Gamma_2$ produces another ribbon graph. This operation induces a pairing of complexes
\[ \gc(m,n) \otimes \gc(n,k) \rightarrow \gc(m,k). \]
\end{rem}

Next, employing the same arguments as those used in the proof of Theorem \ref{thm_kontthm} we obtain the following result. An analogous result in the context of commutative graphs was obtained by M. Kapranov, \cite{Kap}.

\begin{theorem} \label{legs}
The complexes $\gc(m,n)$ and $C^\bullet(\mathfrak{g},\mathfrak{osp};\Hom(V^{\otimes m},V^{\otimes n}))$ are isomorphic.
Under this isomorphism the composition \eqref{comp} corresponds to the pairing of graph complexes induced by gluing.
\end{theorem}
\noproof

Note that the graph complex $\gc(m,n)$ is the chain complex corresponding to a certain cell decomposition of the moduli space of Riemann surfaces with $m$ open incoming boundaries and $n$ open outgoing boundaries. This shows that our definition of a TCFT is essentially equivalent to the standard one, found in e.g. \cite{Costello}. The only difference is the dimension shift (for example, graphs having a single vertex correspond to the cells of the moduli spaces of surfaces sitting in
the highest dimension).

The data of a TCFT is equivalent to a (minimal) symplectic $\ai$-algebra. Given an $\ai$-algebra $A:=(U,m,\omega)$ one associates to any object $[m]$ of $\FT$ the vector space $U^{\otimes m}$. For a graph $\Gamma$ with $m$ incoming and $n$ outgoing legs one associates a homomorphism
\[ U^{\otimes m}\rightarrow U^{\otimes n} \]
by the following generalisation of the partition function (which it would be appropriate to call a \emph{correlation function} in this context). To each internal vertex of $\Gamma$ we associate a tensor as in Definition \ref{partit}. After contracting  along the internal edges we are left with a certain tensor whose valency is determined by the number of legs of $\Gamma$. We interpret this tensor as an element in $\Hom(U^{\otimes m}, U^{\otimes n})$

The analogue of the characteristic class construction within this context is as follows. Note that for any Lie algebra $\mathfrak{l}$ and an $\mathfrak{l}$-module $M$ there is a canonical pairing
\[ \cap:C^\bullet(\mathfrak{l},M)\otimes C_\bullet(\mathfrak{l})\rightarrow M \]
and similarly with relative (co)homology. Then the $\ai$-algebra $A$ associates to any morphism $f\in C^\bullet(\mathfrak{g},\mathfrak{osp};\Hom(V^{\otimes m},V^{\otimes n}))$ the morphism $f\cap c_A\in \Hom(V^{\otimes m},V^{\otimes n})$.

It could then be proved similarly to the proof of Theorem \ref{thm_main} that these two constructions are compatible in the sense that the map \[C^\bullet(\mathfrak{g},\mathfrak{osp};\Hom(V^{\otimes m},V^{\otimes n}))\rightarrow \Hom(U^{\otimes m},U^{\otimes n})\]
specified by the correlation function  gives the same relative homology cycle as the characteristic class construction after taking a composition with the stabilisation map
\[ \Hom(U^{\otimes m},U^{\otimes n})\rightarrow \Hom(V^{\otimes m},V^{\otimes n}). \]

Since an open TCFT is the same as a minimal symplectic $\ai$-algebra the natural question is what can be said about TCFT's which correspond to isomorphic $\ai$-algebras. This question can now be answered as follows.

Associated to the dg category $\FT$ is its homology category $H(\FT)$ having the same objects but whose spaces of morphisms are homologies of the corresponding complexes in $\FT$. Then a TCFT $F:\FT\mapsto\Vect$ induces in an obvious manner a functor $H(F):H(\FT)\mapsto\Vect$. It would be natural to call such a functor a `cohomological TCFT' but that seems to be at odds with the notion of a `cohomological field theory' of \cite{KoM} which involves the Deligne-Mumford compactification of the moduli spaces of curves. We say that two TCFT's $F_1$ and $F_2$ are \emph{homologically equivalent} if the corresponding homology functors $H(F_1)$ and $H(F_2)$ are isomorphic. Then arguments similar to those used in the proof of Theorem \ref{thm_invariant} give the following result.

\begin{theorem}
Topological conformal field theories corresponding to isomorphic minimal symplectic $\ai$-algebras are homologically equivalent.
\end{theorem}
\noproof


\begin{thebibliography}{12}
%
\bibitem{CV}
J. Conant, K. Vogtmann; \emph{On a theorem of Kontsevich.} Algebr. Geom. Topol. 3 (2003), 1167--1224 (electronic).
%
\bibitem{Costello}
K. Costello, \emph{Topological conformal field theories and Calabi-Yau categories.} \texttt{math.QA/0412149}.
%
\bibitem{culv}
M. Culler, K. Vogtmann; \emph{Moduli of graphs and automorphisms of free groups.} Invent. Math. 84 (1986), no. 1, 91--119.
%
\bibitem{fuchscohom}
D. Fuchs, \emph{Cohomology of infinite-dimensional Lie algebras.} Contemporary Soviet Mathematics, Consultants Bureau, New York and London, 1986.
%
\bibitem{GK}
E. Getzler, M. Kapranov; \emph{Modular operads.} Compositio Math. 110 (1998), no. 1, 65--126.
%
\bibitem{ginzsympgeom}
V. Ginzburg, \emph{Noncommutative symplectic geometry, quiver varieties, and operads.} Math. Res. Lett. Vol. 8, 2001, No. 3, pp. 377--400.
%
\bibitem{hamgraph}
A. Hamilton, \emph{A super-analogue of Kontsevich's theorem on graph homology.} Lett. Math. Phys. Vol. 76(1), pp. 37--55, (2006).
%
\bibitem{hamlaz}
A. Hamilton, A. Lazarev; \emph{Homotopy algebras and noncommutative geometry.} \texttt{math.QA/0410621}.
%
\bibitem{har}
J. Harer, \emph{The cohomology of the moduli spaces of curves.} Theory of moduli (Montecatini Terme, 1985), 138--221, Lecture Notes in Math., 1337, Springer, Berlin, 1988.
%
\bibitem{igusa}
K. Igusa, \emph{Graph cohomology and Kontsevich cycles.} Topology 43 (2004), No. 6, pp. 1469--1510.
%
\bibitem{Kap}
M. Kapranov, \emph{Rozansky-Witten invariants via Atiyah classes.} Compositio Math. 115 (1999), no. 1, 71--113.
%
\bibitem{kontsympgeom}
M. Kontsevich, \emph{Formal noncommutative symplectic geometry.} The Gelfand Mathematical Seminars, 1990-1992, pp. 173--187, Birkh\"auser Boston, Boston, MA, 1993.
%
\bibitem{kontfeynman}
M. Kontsevich, \emph{Feynman diagrams and low-dimensional topology.} First European Congress of Mathematics, Vol 2 (Paris, 1992), 97-121, Progr. Math., 120, Birkh\"auser, Basel, 1994.
%
\bibitem{KoM}
M. Kontsevich, Yu. Manin; \emph{Gromov-Witten classes, quantum cohomology, and enumerative geometry.} Comm. Math. Phys. 164 (1994), no. 3, 525--562.
%
\bibitem{KS}
M. Kontsevich, Y. Soibelman; \emph{Notes on $\ai$-algebras, $\ai$-categories and noncommutative geometry.} \texttt{math.RA/0606241}.
%
\bibitem{LV}
A. Lazarev, V. Voronov; \emph{Graph homology: Koszul and Verdier duality.} preprint.
%
\bibitem{loday}
J.-L. Loday, \emph{Cyclic homology.} Grundlehren der mathematischen Wissenschaften 301, second edition, Springer 1998.
%
\bibitem{Manin}
Yu. Manin, \emph{Frobenius manifolds, quantum cohomology, and moduli spaces.} American Mathematical Society Colloquium Publications, 47. American Mathematical Society, Providence, RI, 1999.
%
\bibitem{mondello}
G. Mondello, \emph{Combinatorial classes on $\overline{\mathcal{M}}_{g,n}$ are tautological.} Int. Math. Res. Not. 2004, no. 44, 2329--2390.
%
\bibitem{pen}
R. Penner, \emph{The decorated Teichm\"uller space of punctured surfaces.} Comm. Math. Phys. 113 (1987), no. 2, 299--339.
%
\bibitem{sergeev}
A. Sergeev, \emph{An analog of the classical invariant theory for Lie superalgebras. I, II.} Michigan Math. J. 49 (2001), no. 1, 113--146, 147--168.
%
\bibitem{Stas}
J. Stasheff, \emph{Homotopy associativity of $H$-spaces. I, II.} Trans. Amer. Math. Soc. 108 (1963), 275-292; ibid. 108 1963 293--312.
%
\end{thebibliography}
\end{document}